\documentclass[11pt,a4paper]{amsart}
\usepackage[T1]{fontenc}
\usepackage{lmodern}
\usepackage{amsmath,amssymb,mathtools}
\usepackage[margin=28mm]{geometry}
\usepackage{microtype}
\usepackage{enumitem}
\usepackage{xcolor}
\usepackage[
unicode,
colorlinks=true,
linkcolor=blue,
citecolor=blue,
urlcolor=blue
]{hyperref}
\usepackage{bookmark}
\hypersetup{pdftitle={Fundamental gaps of the restricted fractional Laplacian on thin planar convex domains}}
\newtheorem{theorem}{Theorem}[section]
\newtheorem{proposition}[theorem]{Proposition}
\newtheorem{lemma}[theorem]{Lemma}
\newtheorem{corollary}[theorem]{Corollary}
\theoremstyle{remark}
\newtheorem{remark}[theorem]{Remark}
\numberwithin{equation}{section}
\newcommand{\R}{\mathbb R}
\newcommand{\C}{\mathbb C}
\newcommand{\N}{\mathbb N}
\newcommand{\dd}{\,\mathrm d}
\newcommand{\Hs}[2]{\widetilde H^{#1}(#2)}
\newcommand{\Es}{\mathcal E}
\newcommand{\gap}{\Gamma}
\newcommand{\gs}{\mathfrak g}
\newcommand{\dom}{\operatorname{dom}}
\newcommand{\supp}{\operatorname{supp}}
\newcommand{\diam}{\operatorname{diam}}
\newcommand{\Span}{\operatorname{span}}
\newcommand{\econv}{\xrightarrow{\mathrm{epi}}}
\newcommand{\wh}[1]{\widehat{#1}}
\newcommand{\norm}[1]{\lVert #1\rVert}
\newcommand{\ip}[2]{\langle #1,#2\rangle}
\newcommand{\step}[1]{\par\smallskip\noindent\textbf{#1}\enspace}
\setlist[enumerate]{label=\textup{(\roman*)},leftmargin=2em,itemsep=3pt,topsep=4pt}
\title[Fractional fundamental gaps on convex domains]{Fundamental gaps of the restricted fractional Laplacian on thin planar convex domains}

\author{Rui Chen}
\address{
School of Mathematical Sciences, Fudan University,
Shanghai 200433, China
}
\address{
Brandenburg University of Technology Cottbus--Senftenberg,
03046 Cottbus, Germany
}
\email{chenrui23@m.fudan.edu.cn}

\date{}
\begin{document}

\begin{abstract}
In this paper, we study the fundamental gap $\gap_s(\Omega)=\lambda_{2,s}(\Omega)-\lambda_{1,s}(\Omega)$ of the restricted Dirichlet fractional Laplacian for $0<s<1$. Building on the thin-rectangle estimates of Dyda and Kulczycki \cite{DydaKulczycki2007}, we establish compact Mosco convergence of the shifted quadratic forms and convergence of their eigenvalues on collapsing convex strips. As a consequence, we obtain the next-order asymptotic expansion of every fixed eigenvalue of a thin rectangle, including the exact leading coefficients. We then consider the optimal gap over planar convex domains with prescribed area-to-diameter-squared ratio $\theta$,
\[
G_s(\theta):=\inf_{\substack{\Omega\subset\R^2\ \mathrm{bounded,\ open,\ convex}\\
|\Omega|/(\diam\Omega)^2=\theta}}
(\diam\Omega)^{2s}\gap_s(\Omega).
\]
We prove that, as $\theta\downarrow0$,
\[
G_s(\theta)\sim
\begin{cases}
\kappa_s\mathfrak g_{s+1/2}\theta,&0<s<\frac12,\\
3\pi^2\gamma\theta\log(1/\theta),&s=\frac12,\\
3\pi^2\beta_s\theta^{2-2s},&\frac12<s<1.
\end{cases}
\]
Here $\kappa_s$, $\gamma$, and $\beta_s$ are determined by the transverse interval ground state, while
\[
\mathfrak g_t:=
\inf_{\substack{V\in C(\overline I),\,V\text{ convex}}}
\bigl(\lambda_{2,t}(I,V)-\lambda_{1,t}(I,V)\bigr),
\qquad I=(-1/2,1/2),
\]
with $\lambda_{k,t}(I,V)$ denoting the Dirichlet eigenvalues of $(-\Delta)^t+V$. We prove that $\mathfrak g_t>0$ for all $0<t<1$ and that the infimum is attained by a finite convex potential. Thin rectangles attain the optimal leading asymptotic coefficient of $G_s(\theta)$ for $1/2\le s<1$. For $0<s<1/2$, their asymptotic optimality is equivalent to the one-dimensional equality $\mathfrak g_{s+1/2}=\gap_{s+1/2}(I,0)$. Finally, every asymptotically minimizing sequence converges, after diameter normalization and transverse rescaling, to $[-1/2,1/2]\times[0,1]$ in Hausdorff distance.
\end{abstract}

	\maketitle

	\section{Introduction and main results}

Let $0<s<1$ and let $\Omega\subset\R^2$ be a bounded open set. We consider the restricted Dirichlet fractional Laplacian on $\Omega$, associated with the quadratic form
\[
\Es_s[u]=\int_{\R^2}|\xi|^{2s}|\wh u(\xi)|^2\,\dd\xi,
\qquad
u\in\Hs{s}{\Omega}:=
\{u\in H^s(\R^2):u=0\text{ a.e. on }\Omega^c\}.
\]
Here and throughout, $\wh u(\xi)=(2\pi)^{-d/2}\int_{\R^d}e^{-ix\cdot\xi}u(x)\,\dd x$ denotes the unitary Fourier transform in dimension $d$. The Dirichlet eigenvalues, counted with multiplicity, are denoted by
\[
0<\lambda_{1,s}(\Omega)<\lambda_{2,s}(\Omega)
\le\lambda_{3,s}(\Omega)\le\cdots.
\]
We define the fundamental gap by ${\gap_s(\Omega):=\lambda_{2,s}(\Omega)-\lambda_{1,s}(\Omega)}$. The simplicity of the first eigenvalue and the variational results are recalled in Section~\ref{sec:prelim}.

For the classical Dirichlet Laplacian, Andrews and Clutterbuck \cite{AndrewsClutterbuck2011} proved the sharp lower bound $\gap_1(\Omega)\ge3\pi^2/(\diam\Omega)^2$ for convex domains. In contrast, for every $0<s<1$, the diameter-normalized fundamental gap of the restricted fractional Laplacian tends to zero along thin rectangles, with three distinct decay rates in the transverse width $h$ depending on $s$; see \eqref{eq:DK-orders}. Motivated by these estimates, we study not only the precise gap asymptotics for thin rectangles but also the optimal gap among all planar convex domains with a prescribed area-to-diameter-squared ratio $\theta:=|\Omega|/(\diam\Omega)^2$. We determine the sharp asymptotics of the optimal diameter-normalized gap as $\theta\downarrow0$ and the Hausdorff limits of asymptotically minimizing domains after transverse rescaling.

\subsection{Earlier gap estimates and thin-domain limits}

The classical fundamental gap conjecture developed from van den Berg's observations for several classes of convex domains \cite{vandenBerg1983}, Yau's formulation of the sharp diameter problem \cite{Yau1986}, and the explicit convex-domain and convex-potential conjecture of Ashbaugh and Benguria \cite{AshbaughBenguria1989}. Singer, Wong, Yau and Yau \cite{SingerEtAl1985} proved the lower bound $\pi^2/(4D^2)$, which Yu and Zhong \cite{YuZhong1986} improved to $\pi^2/D^2$; see also Smits \cite{Smits1996}. In one dimension, Ashbaugh and Benguria \cite{AshbaughBenguria1989} obtained the sharp bound for symmetric single-well potentials, and Lavine \cite{Lavine1994} proved it for arbitrary convex potentials. Andrews and Clutterbuck \cite{AndrewsClutterbuck2011} established the full conjecture: for a bounded convex domain $\Omega\subset\R^N$ of diameter $D$ and a convex potential $V$,
\[
\lambda_2(-\Delta+V,\Omega)-\lambda_1(-\Delta+V,\Omega)\ge\frac{3\pi^2}{D^2}.
\]
The constant is sharp in every dimension and is approached by collapsing rectangles. Ni \cite{Ni2013} gave a proof using an elliptic maximum principle, and Gong, Li and Luo \cite{GongLiLuo2016,GongLiLuo2016Erratum} gave a proof using reflection coupling. Amato, Bucur and Fragal\`a \cite{AmatoBucurFragala2025} obtained the quantitative improvement
\[
\lambda_2(\Omega)-\lambda_1(\Omega)\ge\frac{3\pi^2}{D^2}+c_N\frac{w_\Omega^6}{D^8},
\]
where $N\ge2$, $w_\Omega$ is the minimum width, and $c_N>0$ depends only on the dimension.

Fundamental gaps and related eigenvalue inequalities for the $p$-Laplacian have been studied in \cite{FleckingerHarrellDeThelin2007,ChengLianWang2014,AndrewsClutterbuckHauer2021,ChenCheng2020,WenZhou2026,YessirkegenovZhangirbayev2026}. More recently, Chen and Hauer \cite{ChenHauer2026Gap} established the sharp one-dimensional fundamental gap inequality for every $p>1$ and  convex potential, with equality if and only if the potential is constant. They also proved that, for zero potential in dimensions $N\ge2$, the optimal diameter-normalized gap vanishes for $1<p<2$ and equals $3\pi^2$ for $p=2$, with the infimum unattained in both cases, whereas for $p>2$ it is positive and attained.

\smallskip

For the restricted fractional Laplacian, Ba\~nuelos and Kulczycki \cite{BanuelosKulczycki2004} related the Cauchy process to a mixed Steklov problem and obtained variational characterizations of its eigenvalues. They subsequently established eigenvalue gap estimates and weighted Poincar\'e inequalities \cite{BanuelosKulczycki2006JFA}, as well as fundamental gap bounds for planar convex domains symmetric with respect to both coordinate axes \cite{BanuelosKulczycki2006CPDE}. Ba\~nuelos, Kulczycki and M\'endez-Hern\'andez \cite{BanuelosKulczyckiMendez2006} proved midconcavity and coordinatewise monotonicity of the ground state on rectangles for symmetric stable processes.
 Dyda and Kulczycki \cite{DydaKulczycki2007} derived a ground-state variational formula and matching estimates for thin rectangles. 
In our notation, their estimates give
\begin{equation}\label{eq:DK-orders}
\gap_s((-a,a)\times(-b,b))\asymp_s
\begin{cases}
b\,a^{-1-2s},&0<s<\frac12,\\
b\,a^{-2}\log(1+a/b),&s=\frac12,\\
b^{2-2s}a^{-2},&\frac12<s<1,
\end{cases}
\qquad b/a\longrightarrow0.
\end{equation}
Here the comparison constants depend only on $s$. These estimates identify the sharp decay rates but do not determine the leading asymptotic coefficients. Since all three expressions tend to zero as $b/a\to0$ with $a$ fixed, there is no constant $C_s>0$ such that
$\gap_s(\Omega)\ge\frac{C_s}{(\diam\Omega)^{2s}}$
for every bounded planar convex domain $\Omega$. 
Thus, unlike the classical Laplacian, the analysis of fundamental gaps for the restricted fractional Laplacian requires additional geometric constraints to obtain positive lower bounds. Kwa\'snicki \cite{Kwasnicki2008} obtained a gap lower bound involving the first eigenvalue and the diameter. Bao, Ruan, Shen and Sheng \cite{BaoRuanShenSheng2019} investigated spectral and restricted fractional gaps and proposed conjectures involving additional geometric quantities, while Bao, Chen, Jiang and Ma \cite{BaoChenJiangMa2020} studied fractional eigenvalue gaps numerically.

The first eigenvalue and ground state on an interval determine the transverse contribution to the thin-domain asymptotics. Kulczycki, Kwa\'snicki, Ma\l ecki and St\'os \cite{KulczyckiKwasnickiMaleckiStos2010} studied the Cauchy process on the half-line and interval, and Kwa\'snicki \cite{Kwasnicki2012} obtained interval eigenvalue asymptotics. Simplicity of interval eigenvalues was established by Fall, Ghimenti, Micheletti and Pistoia \cite{FallGhimentiMichelettiPistoia2023} and by Kwa\'snicki \cite{Kwasnicki2023Simplicity}. Further results include eigenvalue and eigenfunction estimates by Zhang \cite{Zhang2026Interval} and Dirichlet--Neumann bracketing by Kwa\'snicki and Wszo\l a \cite{KwasnickiWszola2026}. For symmetric single-well potentials, Kaleta \cite{Kaleta2012} proved a potential-independent gap lower bound when $1/2<s<1$. The minimization over convex potentials without symmetry, which enters our result below, is a different question.

For the classical Laplacian, thin-domain spectral limits were studied by Friedlander and Solomyak \cite{FriedlanderSolomyak2009}, Borisov and Freitas \cite{BorisovFreitas2009}, and Lu and Rowlett \cite{LuRowlett2014}. Lu and Rowlett \cite{LuRowlett2013Simplices} also studied fundamental gaps on collapsing simplices, proving divergence of the diameter-normalized gap under collapse to a fixed lower-dimensional simplex. Ambrosio, Freddi and Musina \cite{AmbrosioFreddiMusina2020} proved dimension-reduction results for restricted fractional Dirichlet problems on expanding domains. In the present problem, the first two eigenvalues have the same leading transverse term. Determining their difference therefore requires the next-order quadratic form. For the lower bound over convex domains, we also need to control the first shifted eigenvalue along nearly minimizing sequences, as established in Theorem~\ref{thm:window}.

\subsection{Sharp asymptotics and limiting shapes}

For a bounded open convex domain $\Omega\subset\R^2$, set $D:=\diam(\Omega)$ and $\theta:=|\Omega|/D^2$. The scaling law $\gap_s(r\Omega)=r^{-2s}\gap_s(\Omega)$ shows that both $\theta$ and $D^{2s}\gap_s(\Omega)$ are invariant under dilation. For a rectangle of length $L$ and width $h$, we have $\theta=Lh/(L^2+h^2)\sim h/L$ as $h/L\to0$. Thus the thin-rectangle estimates in \eqref{eq:DK-orders} can be expressed in terms of $\theta$. Motivated by these estimates, we consider the optimal diameter-normalized fundamental gap over all planar convex domains with prescribed $\theta$:
\[
G_s(\theta):=
\inf_{\substack{\Omega\subset\R^2\ \mathrm{bounded,\ open,\ convex}\\
|\Omega|/(\diam\Omega)^2=\theta}}
(\diam\Omega)^{2s}\gap_s(\Omega).
\]

The main difficulty in determining the fundamental gap asymptotics for thin rectangles and general planar convex domains is that the first two eigenvalues have the same leading transverse term. Their difference therefore requires a next-order spectral expansion as the domains collapse. Section~\ref{sec:prelim} establishes the transverse spectral estimates, convex compactness results, and one-dimensional gap inequalities needed for this analysis. In Section~\ref{sec:strip}, we consider the more general convex strips
$K_j=\bigl\{(x,z)\in\R^2:x\in I_j,\ \ell_j(x)<z<u_j(x)\bigr\}$,
where $I_j$ is an open interval, $\ell_j$ is convex, and $u_j$ is concave. 
Theorem~\ref{thm:strip-limit} establishes compact convergence of the shifted quadratic forms to an effective one-dimensional form. Combining this result with the abstract spectral convergence lemma (Lemma~\ref{lem:abstract-spectral}), we obtain convergence of every fixed shifted eigenvalue in Corollary~\ref{cor:strip-spectrum}. Applying the corollary to constant-width strips yields the next-order asymptotic expansion of every fixed eigenvalue of a thin rectangle. 

\smallskip

We first introduce the constants appearing in the asymptotic expansions. Let $I=(-1/2,1/2)$, and let $q_s$ be the nonnegative $L^2(I)$-normalized first eigenfunction of the one-dimensional restricted fractional Laplacian, with eigenvalue $\Lambda_s$. Define
\begin{align}
\beta_s&=s\int_\R|\zeta|^{2s-2}|\wh q_s(\zeta)|^2\dd\zeta,\quad\frac12<s<1,\label{eq:beta}\\
\gamma&=|\wh q_{1/2}(0)|^2
=\frac1{2\pi}\left(\int_Iq_{1/2}(x)\dd x\right)^2,
\label{eq:gamma}\\
C_s&=\int_\R\bigl((1+t^2)^s-|t|^{2s}\bigr)\dd t
=\frac{2s\sqrt\pi\,\Gamma(\frac12-s)}
{(1+2s)\Gamma(1-s)},
\quad \kappa_s=C_s|\wh q_s(0)|^2,\quad 0<s<\frac12.\label{eq:kappa}
\end{align}
All these constants are finite and strictly positive in their respective ranges.

For $0<t<1$ and $V\in C(J)$, let $\lambda_{k,t}(J,V)$ denote the $k$-th Dirichlet eigenvalue of the restricted fractional Schr\"odinger operator $(-\Delta)^t+V$ on a bounded open interval $J$, with zero exterior condition, and set $\gap_t(J,V)=\lambda_{2,t}(J,V)-\lambda_{1,t}(J,V)$.

For $L>0$, let $I_L=(-L/2,L/2)$ and $R_{L,h}=I_L\times(-h/2,h/2)$. We obtain the following next-order asymptotic expansion of every fixed eigenvalue.

\begin{theorem}\label{thm:rectangles}
Fix $L>0$, $0<s<1$, and $n\in\N$. As $h\downarrow0$,
\begin{equation}\label{eq:rectangle-expansion}
\lambda_{n,s}(R_{L,h})=\Lambda_s h^{-2s}+
\begin{cases}
\kappa_s h\lambda_{n,s+1/2}(I_L,0)+o(h),&0<s<\frac12,\\
\gamma(n\pi/L)^2h\log(L/h)+o\bigl(h\log(L/h)\bigr),&s=\frac12,\\
\beta_s(n\pi/L)^2h^{2-2s}+o(h^{2-2s}),&\frac12<s<1.
\end{cases}
\end{equation}
\end{theorem}

Subtracting the expansions for $n=1,2$ in \eqref{eq:rectangle-expansion} determines the exact leading coefficients in the gap estimates \eqref{eq:DK-orders}, previously known only up to multiplicative constants.

\smallskip

We next turn to the optimal gap over general planar convex domains. For $0<s<1/2$, the limiting one-dimensional operator involves a convex potential, and its optimal gap is determined by
\begin{equation}\label{eq:convex-gap-infimum}
\gs_t:=\inf\{\gap_t(I,V):V\in C(\overline I)\text{ is convex}\}.
\end{equation}
Kaleta \cite{Kaleta2012} established a uniform positive fundamental gap for symmetric single-well potentials when $1/2<t<1$. We obtain the following relative gap estimate for arbitrary convex potentials, without any symmetry assumption, throughout $0<t<1$. We also prove that the infimum in \eqref{eq:convex-gap-infimum} is attained.

\begin{proposition}\label{prop:relative-gap}
For every $0<t<1$, there exists $c_t>0$ such that every convex $V\in C(\overline I)$ satisfies
\begin{equation}\label{eq:relative-gap}
\gap_t(I,V)\ge c_t\left(\lambda_{1,t}(I,V)-\min_{\overline I}V\right).
\end{equation}
Consequently, $\gs_t\ge c_t\lambda_{1,t}(I,0)>0$. Moreover,
\[
\gs_t=\min\bigl\{\gap_t(I,V):
V:I\to\R\text{ finite and convex}\bigr\}.
\]
In particular, there exists a finite convex potential $V_*:I\to[0,\infty)$ satisfying
\[
\inf_I V_*=0,
\qquad
\gap_t(I,V_*)=\gs_t.
\]
\end{proposition}

The estimate \eqref{eq:relative-gap} gives a positive lower bound independent of the convex potential and, together with convex compactness, yields the existence of a minimizer that need not extend continuously to $\overline I$. 
The relative gap estimate is proved by rescaling a hypothetical counterexample so that the first two eigenvalues both converge to $1$, contradicting the simplicity of the first eigenvalue of the limiting quadratic form. For attainment, the relative estimate bounds the first two eigenvalues of a minimizing sequence. Convex compactness and spectral convergence then yield a limiting potential on a nonempty interval, while the one-dimensional gap bound in Lemma \ref{lem:convex-gap-bound} excludes a limiting interval shorter than $I$.

The relative estimate is specific to one dimension: even for zero potential, no analogous bound with a domain-independent constant holds for planar convex domains. 

\begin{remark}
The relative gap estimate \eqref{eq:relative-gap} does not extend to dimension two. Indeed, let
$\Omega_\varepsilon=I\times\varepsilon I, \enspace 0<t<1$.
By the transverse Poincar\'e inequality, the subadditivity of $r\mapsto r^t$, and the min--max principle applied to products of the transverse ground state and the first two longitudinal eigenfunctions,
\[
\Lambda_t\varepsilon^{-2t}
\le\lambda_{1,t}(\Omega_\varepsilon)
\le\lambda_{2,t}(\Omega_\varepsilon)
\le\Lambda_t\varepsilon^{-2t}
+\lambda_{2,t}(I,0).
\]
Consequently,
\[
0\le
\frac{\gap_t(\Omega_\varepsilon)}
{\lambda_{1,t}(\Omega_\varepsilon)}
\le
\frac{\lambda_{2,t}(I,0)}{\Lambda_t}
\varepsilon^{2t}
\longrightarrow0,
\]
showing that \eqref{eq:relative-gap} fails uniformly in dimension two.
\end{remark}

We next determine the optimal fundamental gap among all planar convex domains. For $1/2\le s<1$, Theorem~\ref{thm:rectangles} provides the required upper bounds. For $0<s<1/2$, we construct concave transverse width functions that realize prescribed convex potentials in the limiting one-dimensional quadratic form; see \eqref{eq:width-realization}. This gives the upper bound with coefficient $\kappa_s\gs_{s+1/2}$. For the lower bounds, Theorem~\ref{thm:window} controls the first shifted eigenvalue of nearly minimizing domains using only an upper bound on their fundamental gaps. We then apply Theorem~\ref{thm:strip-limit} and Corollary~\ref{cor:strip-spectrum} to obtain convergence of the first two shifted eigenvalues. The one-dimensional convex-potential gap estimates in Lemma~\ref{lem:convex-gap-bound} yield the matching lower bounds. We thus obtain our main theorem.

\begin{theorem}\label{thm:main}
For every fixed $0<s<1$, as $\theta\downarrow0$,
\begin{equation}\label{eq:main-asymptotics}
G_s(\theta)=
\begin{cases}
\kappa_s\gs_{s+1/2}\,\theta+o(\theta),&0<s<\frac12,\\
3\pi^2\gamma\,\theta\log(1/\theta)+o\bigl(\theta\log(1/\theta)\bigr),&s=\frac12,\\
3\pi^2\beta_s\,\theta^{2-2s}+o(\theta^{2-2s}),&\frac12<s<1.
\end{cases}
\end{equation}
Thin rectangles in Theorem \ref{thm:rectangles} attain the leading constant for $1/2\le s<1$. For $0<s<1/2$, they attain it if and only if
\begin{equation}\label{eq:remaining-1d}
\gs_{s+1/2}=\gap_{s+1/2}(I,0).
\end{equation}
\end{theorem}

The three asymptotic formulas in \eqref{eq:main-asymptotics} are determined by the corresponding one-dimensional effective operators. For $1/2\le s<1$, the effective operator is the classical Schr\"odinger operator $-a_s(H)d^2/dx^2+W$, with an additional logarithmic factor in the normalization at $s=1/2$. The sharp convex-potential gap inequality in Lemma~\ref{lem:convex-gap-bound}, based on \cite{Lavine1994}, then gives the constant $3\pi^2$ on the unit interval. For $0<s<1/2$, the effective operator is instead the nonlocal Schr\"odinger operator $a_s(H)(-\Delta)^{s+1/2}+W$, where variations in the transverse width produce the convex potential $W$. The optimal coefficient therefore involves $\gs_{s+1/2}$. 

\smallskip

Finally, we determine the limiting shapes of asymptotically minimizing convex domains. We show that, after normalization by the diameter and rescaling in the minimum-width direction, every such sequence converges to a rectangle in Hausdorff distance.

\begin{theorem}\label{thm:shape}
Let $0<s<1$, and let $\Omega_j\subset\R^2$ be bounded open convex domains satisfying
\[
\diam\Omega_j=1,\qquad |\Omega_j|=\theta_j\downarrow0,
\qquad \frac{\gap_s(\Omega_j)}{G_s(\theta_j)}\longrightarrow1.
\]
Orient each domain so that the vertical direction realizes its minimal width, and define
$K_j=\{(x,z):(x,\theta_jz)\in\Omega_j\}$.
Then, after suitable translations,
\begin{equation}\label{eq:shape-limit}
d_H\!\left(\overline K_j,[-\tfrac12,\tfrac12]\times[0,1]\right)\longrightarrow0.
\end{equation}
\end{theorem}

The lower-bound estimates used to prove Theorem~\ref{thm:main} show that asymptotic optimality forces the limiting maximal transverse height $H$ and the effective interval length $L$ to equal one. The area constraint and concavity of the transverse width functions then yield \eqref{eq:shape-limit}. For $0<s<1/2$, however, this Hausdorff convergence does not determine the effective potential, as the following remark shows.

\begin{remark}\label{rem:shape-versus-potential}
For $0<s<1/2$, set $t=s+1/2$. Given any convex $V\in C(\overline I)$, the construction in \eqref{eq:width-realization}, after normalizing $\min_{\overline I}V=0$, produces convex strips $K_\varepsilon$ satisfying
$\overline K_\varepsilon\longrightarrow \overline I\times[0,1]$
in Hausdorff distance, while their normalized fundamental gaps converge to $\kappa_s\gap_t(I,V)$ by \eqref{eq:realization-gap}. Indeed, variations of order $\varepsilon^{1+2s}$ in the transverse width vanish in the Hausdorff limit but produce the effective potential $\kappa_sV$. Thus the rectangular limit in Theorem~\ref{thm:shape} does not imply that thin rectangles attain the optimal subcritical coefficient $\kappa_s\gs_t$.

By the definition of $\gs_t$, the remaining equality \eqref{eq:remaining-1d} is equivalent to
\[
\gap_t(I,V)\ge\gap_t(I,0)
\qquad
\text{for every convex }V\in C(\overline I),
\quad \frac12<t<1.
\]
Proposition~\ref{prop:relative-gap} guarantees the existence of a minimizing convex potential, but whether a constant potential is a minimizer remains open.
\end{remark}

\subsection{The effective one-dimensional form and the gap asymptotics}

The proofs of our main results rely on Sections~\ref{sec:strip} and \ref{sec:window}. Section~\ref{sec:strip} establishes convergence of the shifted quadratic forms and their eigenvalues on collapsing convex strips, while Section~\ref{sec:window} derives the first-eigenvalue bound needed to apply this convergence to nearly minimizing domains. We explain the main steps below.

\smallskip

The first two eigenvalues of a thin domain have the same leading transverse term. To determine their difference, we subtract this term and study the next-order quadratic form. For $\diam\Omega_\varepsilon=1$ and $|\Omega_\varepsilon|=\varepsilon$, rescale the minimum-width direction by setting $K_\varepsilon=\{(x,z):(x,\varepsilon z)\in\Omega_\varepsilon\}$. The corresponding anisotropic form and its normalized remainder are
\[
\mathfrak a_\varepsilon[v]
=\iint_{\R^2}(\zeta^2+\varepsilon^2\xi^2)^s
|\wh v(\xi,\zeta)|^2\,\dd\xi\dd\zeta,
\qquad
F_\varepsilon[v]
=\frac{\mathfrak a_\varepsilon[v]-c_\varepsilon\|v\|_2^2}
{\rho_s(\varepsilon)},
\]
where $c_\varepsilon=\Lambda_sH_\varepsilon^{-2s}$, $H_\varepsilon$ is the maximal transverse section length, and $\rho_s$ is defined in \eqref{eq:rho}.

The key step is to identify the limit of $F_\varepsilon$. It contains two terms. The expansion of the lowest transverse spectral band in Lemma~\ref{lem:band} produces the longitudinal kinetic energy, while variations in the transverse width $h_\varepsilon(x)$ produce the potential
\[
W_\varepsilon(x)
=\frac{\Lambda_s\bigl(h_\varepsilon(x)^{-2s}
-H_\varepsilon^{-2s}\bigr)}{\rho_s(\varepsilon)}.
\]
Since $h_\varepsilon$ is concave, $W_\varepsilon$ is nonnegative and convex. Under the hypotheses of Theorem~\ref{thm:strip-limit}, the limiting one-dimensional form is
\[
F[\psi]
=a_s(H)\int_\R|\xi|^{2\sigma_s}|\wh\psi(\xi)|^2\,\dd\xi
+\int_JW(x)|\psi(x)|^2\,\dd x,
\]
where $\sigma_s=\min\{s+1/2,1\}$, $J=\operatorname{int}\dom W$, and $a_s(H)$ is given in \eqref{eq:effective-coefficient}. Thus the effective operator is nonlocal for $0<s<1/2$ and local for $1/2\le s<1$, with a logarithmic normalization at $s=1/2$.

Theorem~\ref{thm:strip-limit} proves compact Mosco convergence to this form, and Corollary~\ref{cor:strip-spectrum} gives convergence of every fixed shifted eigenvalue. For rectangles, the transverse width is constant and $W=0$, yielding Theorem~\ref{thm:rectangles}. The same convergence also gives the upper bounds for general convex domains. For $s\ge1/2$, thin rectangles provide the optimal coefficients. For $0<s<1/2$, the concave width construction in \eqref{eq:width-realization} realizes any prescribed continuous convex potential $V$, with $W_\varepsilon\to\kappa_sV$. This yields the upper coefficient $\kappa_s\gs_{s+1/2}$.

\smallskip

For the lower bounds, the main difficulty is that a gap bound does not directly control the individual shifted eigenvalues. Theorem~\ref{thm:window} proves that, for nearly minimizing domains,
\[
\nu_{2,j}-\nu_{1,j}=O(\rho_s(\varepsilon_j))
\quad\Longrightarrow\quad
0\le\nu_{1,j}-c_j=O(\rho_s(\varepsilon_j)).
\]
The proof uses a second longitudinal rescaling: if this estimate failed, the first two shifted eigenvalues would converge to the same lowest eigenvalue of a limiting one-dimensional operator, contradicting its simplicity. We can then apply Corollary~\ref{cor:strip-spectrum} and Lemma~\ref{lem:convex-gap-bound} to obtain the matching lower bounds, proving Theorem~\ref{thm:main}. Finally, asymptotic optimality forces the limiting maximal transverse height $H$ and effective interval length $L$ to satisfy $H=L=1$. The area constraint and concavity of the width functions then yield the Hausdorff convergence in Theorem~\ref{thm:shape}.

\subsection{Remarks on the limiting problems and proofs}
\label{remark}
We emphasize four points concerning the assumptions and proofs.

\begin{enumerate}
\item[(i)] \emph{Limiting potentials.} The effective interval $J=\operatorname{int}\dom V$ may be smaller than the original interval or unbounded after rescaling, while $V$ may diverge at its endpoints. Accordingly, Lemma~\ref{lem:convex-compactness} and Proposition~\ref{prop:potential-convergence} allow extended-valued convex potentials, and Lemma~\ref{lem:form-core} also covers unbounded intervals.

\item[(ii)] \emph{Minimizing convex potentials.} Proposition~\ref{prop:relative-gap} proves that $\gs_t>0$ and is attained by a finite convex potential $V_*:I\to[0,\infty)$ with $\inf_I V_*=0$. Such a minimizer is continuous and locally Lipschitz in $I$, but need not be assumed continuous up to the boundary. Its uniqueness and whether it must be constant remain open.

\item[(iii)] \emph{Rescaling and simplicity.} Proposition~\ref{prop:relative-gap} and Theorem~\ref{thm:window} use the same contradiction argument. In the former, we rescale to normalize the first eigenvalue to $1$; in the latter, we normalize its excess above the leading transverse eigenvalue to $1$. If the respective estimates failed, the second normalized eigenvalue would also converge to $1$. Spectral convergence would then yield a limiting operator whose first two eigenvalues coincide, contradicting the simplicity of its first eigenvalue.

\item[(iv)] \emph{Further applications.} The analysis identifies the longitudinal kinetic energy through the lowest transverse spectral band and the effective potential through variations of the transverse width. The same approach may apply to other collapsing spectral problems, provided the corresponding band expansion, compact spectral convergence, and first-eigenvalue bounds can be established.
\end{enumerate}

\medskip

Section~\ref{sec:prelim} collects the transverse spectral estimates, convex compactness results, and one-dimensional gap inequalities. Section~\ref{sec:strip} establishes compact convergence of the shifted quadratic forms on convex strips and derives the eigenvalue expansions for thin rectangles. Section~\ref{sec:window} derives a bound on the first shifted eigenvalue from a fundamental gap bound. Finally, Section~\ref{sec:optimization} proves the sharp gap asymptotics and the Hausdorff convergence of asymptotically minimizing domains.

	\section{Preliminary spectral and variational results}\label{sec:prelim}
	
	In this section, we establish several preliminary results used in the proofs of our main theorems. We first derive the small-frequency asymptotics of the first transverse spectral band. We then establish compactness and spectral convergence results for one-dimensional fractional Schrödinger operators with convex potentials. Finally, we prove uniform lower bounds for their fundamental gaps and the existence of a minimizing convex potential.
	
	\smallskip
	
	Throughout this paper, for $0<t\le1$ and an open set $U\subset\R^d$, we consider
	\[
	\Hs{t}{U}
	:=
	\bigl\{u\in H^t(\R^d):u=0\text{ a.e. on }U^c\bigr\}.
	\]
	Another natural space is the closure of $C_c^\infty(U)$ in the full-space Sobolev norm:
	\[
	H_0^t(U)
	:=
	\overline{C_c^\infty(U)}^{\,H^t(\R^d)}.
	\]
	For general open sets, these spaces satisfy $H_0^t(U)\subset\Hs{t}{U},$
	and the inclusion may be strict. If $U$ has a continuous boundary, then
	$H_0^t(U)=\Hs{t}{U}, \enspace 0<t\le 1$;
	see \cite[Lemma~3.15]{ChandlerWildeHewettMoiola2017}. In particular, this equality holds for bounded open intervals. For the one-dimensional Schrödinger forms with possibly unbounded convex potentials considered below, the corresponding form-core property is established in Lemma~\ref{lem:form-core}.
	
	We define the Dirichlet quadratic form
	\[
	\Es_t[u]
	:=
	\int_{\R^d}|\xi|^{2t}|\wh u(\xi)|^2\dd\xi,
	\qquad u\in\Hs{t}{U}.
	\]
	This form is densely defined, closed, and nonnegative on $L^2(U)$. If $U$ is bounded, the embedding $\Hs{t}{U}\hookrightarrow L^2(U)$
	is compact: for $0<t<1$, apply \cite[Theorem~7.1]{DiNezzaPalatucciValdinoci2012}
	on a smooth bounded ball containing $U$, using the zero extensions of the functions;
	for $t=1$, use the classical Rellich--Kondrachov theorem on that ball. Consequently, the associated self-adjoint operator has compact resolvent. Its eigenvalues, counted with multiplicity, satisfy
	\[
	0<\lambda_{1,t}(U)\le\lambda_{2,t}(U)\le\cdots,
	\qquad
	\lambda_{k,t}(U)\longrightarrow\infty
	\quad\text{as }k\to\infty,
	\]
	where the associated self-adjoint operator is defined by the
	first representation theorem
	\cite[Chapter~VI, Theorem~2.1]{Kato1995}, and its eigenvalues
	satisfy the min--max principle
	\cite[Chapter XIII, Theorem XIII.2]{ReedSimonIV}.
	
	For $0<t<1$ and every real-valued $u\in H^t(\R^d)$, the Gagliardo representation gives
	\begin{equation}\label{eq:absolute-value}
		\Es_t[u]-\Es_t[|u|]
		=
		c_{d,t}\iint_{\R^d\times\R^d}
		\frac{|u(x)u(y)|-u(x)u(y)}
		{|x-y|^{d+2t}}\dd x\dd y,
	\end{equation}
	where $c_{d,t}=4^t\Gamma(d/2+t)/(\pi^{d/2}|\Gamma(-t)|)$; see \cite{Kwasnicki2017}.
	The right-hand side is strictly positive whenever $u$ changes sign on sets of positive measure. Hence every real first eigenfunction has a fixed sign. Moreover, the strong minimum principle \cite[Theorem~1.2]{DelPezzoQuaas2017} implies that a nonnegative first eigenfunction is strictly positive almost everywhere in $U$, even when $U$ is disconnected. If the first eigenvalue were multiple, two orthogonal real first eigenfunctions could both be chosen strictly positive, a contradiction. Consequently,
	$0<\lambda_{1,t}(U)<\lambda_{2,t}(U), \enspace 0<t<1$.
	
	Finally, a change of variables in the Fourier form gives the scaling identity
	\[\lambda_{k,t}(rU)
	=
	r^{-2t}\lambda_{k,t}(U),
	\qquad r>0,\quad k\ge1.\]
	Boundary regularity results for the fractional Laplacian can be found in \cite{Grubb2015,RosOtonSerra2014}. No boundary regularity estimate for planar convex domains is required in the arguments below.

	\subsection{The first transverse band}
	We introduce the three scales arising from the
	small-frequency expansion of the first transverse band.
	For $0<r<e^{-1}$, define
	\begin{equation}\label{eq:rho}
		\rho_s(r)=
		\begin{cases}
			r^{1+2s},&0<s<\tfrac12,\\
			r^2\log(1/r),&s=\tfrac12,\\
			r^2,&\tfrac12<s<1,
		\end{cases}
		\qquad \sigma_s=\min\{s+\tfrac12,1\}.
	\end{equation}
	For $H>0$, define
	\begin{equation}\label{eq:effective-coefficient}
		a_s(H)=
		\begin{cases}
			\kappa_s H, & 0<s<\frac12,\\
			\gamma H, & s=\frac12,\\
			\beta_s H^{2-2s}, & \frac12<s<1,
		\end{cases}
	\end{equation}
	where $\kappa_s$, $\gamma$, and $\beta_s$ are defined in
	\eqref{eq:kappa}, \eqref{eq:gamma}, and \eqref{eq:beta}, respectively.

	Let $q_s$ be the normalized first eigenfunction
	on $I=(-1/2,1/2)$, with eigenvalue $\Lambda_s$.
	By scaling, the normalized first eigenfunction on $(0,H)$ is
	\[
	q_{s,H}(z)
	=
	H^{-1/2}q_s\left(\frac{z-H/2}{H}\right),
	\]
	with eigenvalue $\Lambda_sH^{-2s}$. Its Fourier transform
	satisfies
	$\widehat q_{s,H}(\zeta) = H^{1/2}e^{-iH\zeta/2}\widehat q_s(H\zeta)$.
	
	The following lemma gives the small-frequency
	asymptotics of the first transverse band.

	\begin{lemma}\label{lem:band}
		Define
		\[
		B_s(k)=\inf_{\substack{w\in\Hs{s}{I},\,\norm w_2=1}}
		\int_\R(\zeta^2+k^2)^s|\wh w(\zeta)|^2\dd\zeta,\qquad k\ge0.
		\]
		Then $B_s(0)=\Lambda_s$, $B_s$ is nondecreasing, and
		\begin{equation}\label{eq:band-asymptotic}
			B_s(k)-\Lambda_s\sim
			\begin{cases}
				\kappa_s k^{1+2s},&0<s<\tfrac12,\\
				\gamma k^2\log(1/k),&s=\tfrac12,\\
				\beta_s k^2,&\tfrac12<s<1,
			\end{cases}
			\qquad k\downarrow0.
		\end{equation}
		In particular, $B_s(k)>\Lambda_s$ for every $k>0$. For any interval of length $H>0$, the first transverse
		band is $H^{-2s}B_s(Hk)$. Consequently, for each fixed
		$H>0$,
		\begin{equation}\label{eq:band-scaled}
			\lim_{k\downarrow0}
			\frac{H^{-2s}(B_s(Hk)-\Lambda_s)}{\rho_s(k)}
			=a_s(H).
		\end{equation}
		The convergence is uniform for $H$ in compact subsets
		of $(0,\infty)$.
	\end{lemma}
	
	\begin{proof}
		\step{Step 1. Multiplier estimates.}
		Set $D_{s,k}(\zeta)=(\zeta^2+k^2)^s-|\zeta|^{2s}$. Concavity and subadditivity of $r\mapsto r^s$ give
		\begin{equation}\label{eq:D-bound}
			0\le D_{s,k}(\zeta)\le\min\{k^{2s},sk^2|\zeta|^{2s-2}\},\qquad \zeta\ne0.
		\end{equation}
		We first compute the increment for any family $v_k\in L^2(I)$ such that $v_k\to q_s$ in $L^2(I)$. Put $F_k=|\wh v_k|^2$ and $F=|\wh q_s|^2$. Since the functions are supported in a fixed interval,
		$$\norm{\wh v_k-\wh q_s}_\infty\le(2\pi)^{-1/2}\norm{v_k-q_s}_1\le(2\pi)^{-1/2}\norm{v_k-q_s}_2,$$
		thus, we have
		\begin{equation}\label{eq:band-F-convergence}
			\norm{F_k-F}_\infty\longrightarrow0,
			\qquad \sup_k(\norm{F_k}_\infty+\norm{F_k}_1)<\infty.
		\end{equation}

		\step{Step 2. The three scalar limits.}
		If $0<s<1/2$, the function $\Phi_s(t)=(1+t^2)^s-|t|^{2s}$ belongs to $L^1(\R)$, because $\Phi_s(t)\sim s|t|^{2s-2}$ at infinity. Tonelli's theorem gives
		\begin{align*}
			\int_\R\Phi_s(t)\dd t
			&=s\int_0^1\int_\R(t^2+r)^{s-1}\dd t\dd r\\
			&=s\int_0^1r^{s-1/2}\dd r\int_\R(1+y^2)^{s-1}\dd y
			=\frac{2s\sqrt\pi\,\Gamma(\frac12-s)}{(1+2s)\Gamma(1-s)}=C_s,
		\end{align*}
		where $C_s$ is defined in \eqref{eq:kappa}.
		The substitution $\zeta=kt$, followed by dominated convergence using \eqref{eq:band-F-convergence}, yields
		\[ \frac1{k^{1+2s}}\int_\R D_{s,k}(\zeta)F_k(\zeta)\dd\zeta
		=\int_\R\Phi_s(t)F_k(kt)\dd t
		\longrightarrow C_sF(0)=\kappa_s.\]
		
		If $s=1/2$, rationalization gives
		\[
		D_{1/2,k}(\zeta)
		=
		\frac{k^2}{\sqrt{\zeta^2+k^2}+|\zeta|}.
		\]
		For every fixed $\delta>0$, direct integration yields
		\[  M_\delta(k):=
		\int_{-\delta}^{\delta}
		D_{1/2,k}(\zeta)\,\dd\zeta=
		\delta\sqrt{\delta^2+k^2}-\delta^2
		+k^2\log
		\frac{\delta+\sqrt{\delta^2+k^2}}{k}=
		k^2\log(1/k)+O_\delta(k^2).\]
		By \eqref{eq:band-F-convergence}, we obtain
		\[
		\begin{aligned}
			&\left|
			\int_{-\delta}^{\delta}
			D_{1/2,k}(\zeta)
			\bigl(F_k(\zeta)-F(0)\bigr)\,\dd\zeta
			\right|\le
			M_\delta(k)
			\left(
			\|F_k-F\|_\infty
			+
			\sup_{|\zeta|\le\delta}|F(\zeta)-F(0)|
			\right).
		\end{aligned}
		\]
		On the complement, \eqref{eq:D-bound} gives
		\[
		\int_{|\zeta|>\delta}
		D_{1/2,k}(\zeta)F_k(\zeta)\,\dd\zeta
		\le
		\frac{k^2}{2\delta}\|F_k\|_1.
		\]
		Dividing by $k^2\log(1/k)$, letting $k\downarrow0$,
		and then letting $\delta\downarrow0$, we obtain
		\[\lim_{k\downarrow0}
		\frac1{k^2\log(1/k)}
		\int_{\R}D_{1/2,k}F_k\,\dd\zeta
		=
		F(0)
		=
		|\widehat q_{1/2}(0)|^2
		=\gamma.\]
		
		If $1/2<s<1$, let $w(\zeta)=|\zeta|^{2s-2}$. Since $w$ is integrable near zero and bounded for $|\zeta|\ge1$,
		\[
		\int_\R w|\wh v_k-\wh q_s|^2\dd\zeta
		\le\frac{2}{2s-1}\norm{\wh v_k-\wh q_s}_\infty^2+\norm{v_k-q_s}_2^2\longrightarrow0.
		\]
		Weighted Cauchy--Schwarz inequality therefore gives $\int_\R w|F_k-F|\to0$. As $D_{s,k}/k^2\to sw$ pointwise for every $\zeta\ne 0$ and $D_{s,k}/k^2\le sw$, we obtain
		\[ \frac1{k^2}\int_\R D_{s,k}F_k\dd\zeta
		\longrightarrow s\int_\R|\zeta|^{2s-2}|\wh q_s(\zeta)|^2\dd\zeta=\beta_s.\]
		
		\step{Step 3. Minimization and scaling.}
		Testing with $q_s$ gives
		\[
		0\le B_s(k)-\Lambda_s\le\int_\R D_{s,k}|\wh q_s|^2\dd\zeta\le k^{2s}.
		\]
		For each $k>0$, the direct method and compactness on $I$ give a normalized minimizer $w_k$. Let $P$ be the projection onto $\Span\{q_s\}$. The transverse spectral gap implies
		\[
		(\lambda_{2,s}(I)-\Lambda_s)\norm{(1-P)w_k}_2^2
		\le\int_\R|\zeta|^{2s}|\wh w_k|^2\dd\zeta-\Lambda_s
		\le B_s(k)-\Lambda_s\longrightarrow0.
		\]
		Choose the phase so that $\ip{w_k}{q_s}\ge0$. Then $w_k\to q_s$ in $L^2(I)$, and
		\[
		\int_\R D_{s,k}|\wh w_k|^2\dd\zeta
		\le B_s(k)-\Lambda_s
		\le\int_\R D_{s,k}|\wh q_s|^2\dd\zeta.
		\]
		The limits in Step~2 prove \eqref{eq:band-asymptotic}. Monotonicity follows from the defining variational problem; the positive small-$k$ asymptotics and monotonicity imply strict positivity of the increment for all $k>0$. Dilation gives the band $H^{-2s}B_s(Hk)$. Formula \eqref{eq:band-scaled} follows from \eqref{eq:band-asymptotic}; in the critical case, $\log(1/(Hk))/\log(1/k)\to1$ uniformly on compact $H$-ranges.
	\end{proof}
	
	We next establish the asymptotics of the spectral increment
	for product functions, including the mixed terms. Recall that
	$D_{s,k}(\zeta):=(\zeta^2+k^2)^s-|\zeta|^{2s}$.
	
	\begin{lemma}\label{lem:product-limit}
		Let $\eta_j\downarrow0$, $H>0$, and suppose that $r_j,t_j\to q_{s,H}$ in
		$L^1(\R)\cap L^2(\R)$. For every $f,g\in C_c^\infty(\R)$,
		\begin{equation}\label{eq:product-limit}
			\frac1{\rho_s(\eta_j)}\iint_{\R^2}
			D_{s,\eta_j|\xi|}(\zeta)\wh f(\xi)\overline{\wh g(\xi)}
			\wh r_j(\zeta)\overline{\wh t_j(\zeta)}\dd\xi\dd\zeta
			\longrightarrow a_s(H)\int_\R|\xi|^{2\sigma_s}\wh f(\xi)\overline{\wh g(\xi)}\dd\xi,
		\end{equation}
		where $\rho_s$, $\sigma_s$ and $a_s$ are defined in \eqref{eq:rho} and \eqref{eq:effective-coefficient}, respectively.
	\end{lemma}
	
	\begin{proof}
		\step{Step 1. Fourier convergence.}
		Put $G_j(\zeta)=\wh r_j(\zeta)\overline{\wh t_j(\zeta)}$ and
		$G(\zeta)=|\wh q_{s,H}(\zeta)|^2$. Since the Fourier transform maps
		$L^1(\R)$ continuously into $C_0(\R)$ and by Plancherel formula,
		\begin{equation}\label{eq:product-Fourier-bounds}
			\norm{G_j-G}_{L^\infty(\R)}\to0,
			\qquad \sup_j\bigl(\norm{G_j}_\infty+\norm{G_j}_1\bigr)<\infty.
		\end{equation}
		For $1/2<s<1$, the weight $w(\zeta)=|\zeta|^{2s-2}$ is integrable near zero
		and bounded outside $(-1,1)$. Splitting at $|\zeta|=1$, we obtain
		\[
		\int_\R w|\wh r_j-\wh q_{s,H}|^2\dd\zeta
		\le \frac2{2s-1}\norm{\wh r_j-\wh q_{s,H}}_\infty^2
		+\norm{r_j-q_{s,H}}_2^2\longrightarrow0.
		\]
		The analogous estimate for $t_j$ and weighted Cauchy--Schwarz yield
		\begin{equation}\label{eq:product-weighted}
			\int_\R w|G_j-G|\dd\zeta\longrightarrow0,
			\qquad\sup_j\int_\R w|G_j|\dd\zeta<\infty.
		\end{equation}
		
		For $\xi\in\R$, write
		\[
		I_j(\xi)=\frac1{\rho_s(\eta_j)}
		\int_\R D_{s,\eta_j|\xi|}(\zeta)G_j(\zeta)\dd\zeta.
		\]
		We show that $I_j(\xi)\to a_s(H)|\xi|^{2\sigma_s}$ for every fixed
		$\xi$, and then find a dominating function whose product with
		$|\wh f\,\wh g|$ is integrable.
		
		\step{Step 2. The subcritical case $0<s<1/2$.}
		Let $p=1+2s$, so that $\rho_s(\eta_j)=\eta_j^p$. When $\xi\ne0$,
		$\zeta=\eta_j|\xi|t$ gives
		\[
		I_j(\xi)=|\xi|^p\int_\R\Phi_s(t)G_j(\eta_j|\xi|t)\dd t,
		\qquad \Phi_s(t)=(1+t^2)^s-|t|^{2s}.
		\]
		The integrand converges pointwise to $\Phi_s(t)G(0)$.
		Since $\Phi_s\in L^1(\R)$ and $\sup_j\norm{G_j}_\infty<\infty$,
		dominated convergence gives
		\[
		I_j(\xi)\longrightarrow C_s|\xi|^p|\wh q_{s,H}(0)|^2
		=\kappa_sH|\xi|^p.
		\]
		Furthermore, $D_{s,k}\ge0$ and the exact identity
		$\int_\R D_{s,k}(\zeta)\dd\zeta=C_sk^p$ give
		\[|I_j(\xi)|\le\norm{G_j}_\infty C_s|\xi|^p\le C|\xi|^p.\]
		The conclusion at $\xi=0$ follows because $D_{s,0}=0$. Since
		$|\xi|^p|\widehat f\,\widehat g|\in L^1(\R)$,
		dominated convergence in $\xi$ proves
		\eqref{eq:product-limit} for $0<s<1/2.$
		
		\step{Step 3. The critical case $s=1/2$.}
		Let $L_j=\log(1/\eta_j)$, and $k_j=\eta_j|\xi|$.
		Fix $\xi\ne0$ and $\delta\in(0,1)$. The calculation of Step 2
		in Lemma \ref{lem:band} gives
		\[
		\frac1{\eta_j^2L_j}\int_{|\zeta|\le\delta}
		D_{1/2,k_j}(\zeta)\dd\zeta
		=\xi^2\frac{L_j-\log|\xi|+O_\delta(1)}{L_j}
		\longrightarrow\xi^2.
		\]
		By \eqref{eq:product-Fourier-bounds} and $D_{1/2,k_j}\ge0$,
		\[
		\begin{aligned}
			&\frac1{\eta_j^2L_j}\left|
			\int_{|\zeta|\le\delta}
			D_{1/2,k_j}(\zeta)(G_j(\zeta)-G(0))\,\dd\zeta
			\right|\\
			&\quad\le
			\left(\|G_j-G\|_\infty+
			\sup_{|\zeta|\le\delta}|G(\zeta)-G(0)|\right)
			\frac1{\eta_j^2L_j}\int_{|\zeta|\le\delta}
			D_{1/2,k_j}(\zeta)\,\dd\zeta.
		\end{aligned}
		\]
		Its limit superior is at most
		$\xi^2\sup_{|\zeta|\le\delta}|G(\zeta)-G(0)|$,
		which tends to zero as $\delta\downarrow0$. On $|\zeta|>\delta$, the estimate $D_{1/2,k}(\zeta)
		\le\frac{k^2}{2|\zeta|}$
		and \eqref{eq:product-Fourier-bounds} give
		\[\frac1{\rho_{1/2}(\eta_j)}\int_{|\zeta|>\delta}
		D_{1/2,k_j}|G_j|\,\dd\zeta
		\le\frac{C\xi^2}{\delta L_j}\longrightarrow0.\]
		Taking first $j\to\infty$ and then
		$\delta\downarrow0$, we obtain $I_j(\xi)\to\xi^2G(0)=\gamma H\xi^2.$
		
		\smallskip
		
		We next establish a uniform bound.  Assume $L_j\ge1$ and put $k=\eta_j|\xi|$. If $0<k\le1$, then
		\[
		0\le D_{1/2,k}(\zeta)
		\le\min\left\{k,\frac{k^2}{2|\zeta|}\right\},\qquad\zeta\ne0.
		\]
		Splitting the integral into
		$|\zeta|\le k$, $k<|\zeta|\le1$, and
		$|\zeta|>1$, we obtain from
		\eqref{eq:product-Fourier-bounds}
		\[
		\begin{aligned}
			\int_{|\zeta|\le k}D_{1/2,k}|G_j|\,\dd\zeta
			&\le Ck^2,\quad 
			\int_{|\zeta|>1}D_{1/2,k}|G_j|\,\dd\zeta\le Ck^2,
		\end{aligned}
		\]
		and
		\[\int_{k<|\zeta|\le1}D_{1/2,k}|G_j|\,\dd\zeta
		\le Ck^2\log(1/k).\]
		Hence
		\[
		|I_j(\xi)|
		\le C\xi^2\frac{1+\log(1/(\eta_j|\xi|))}{L_j}
		\le C\xi^2\bigl(1+\log_+(1/|\xi|)\bigr).
		\]
		If $\eta_j|\xi|>1$, using $D_{1/2,k}\le k$ gives $|I_j(\xi)|
		\le\frac{C|\xi|}{\eta_jL_j}
		\le C\xi^2.$
		Thus
		\[|I_j(\xi)|
		\le C\xi^2\bigl(1+\log_+(1/|\xi|)\bigr),
		\qquad \xi\ne0,\]
		for all sufficiently large $j$, with $I_j(0)=0$.
		Since $\xi^2(1+\log_+(1/|\xi|))
		|\widehat f\,\widehat g|\in L^1(\R)$,
		dominated convergence in $\xi$ proves
		\eqref{eq:product-limit} at $s=1/2$.
		
		\step{Step 4. The supercritical case $1/2<s<1$.}
		By \eqref{eq:D-bound},
		\[
		0\le \frac{D_{s,\eta_j|\xi|}(\zeta)}{\eta_j^2}
		\le s\xi^2|\zeta|^{2s-2},
		\qquad
		\frac{D_{s,\eta_j|\xi|}(\zeta)}{\eta_j^2}
		\longrightarrow s\xi^2|\zeta|^{2s-2},
		\quad \zeta\ne0.
		\]
		The weighted convergence \eqref{eq:product-weighted} controls the
		replacement of $G_j$ by $G$; dominated convergence for the fixed
		$wG\in L^1$ then gives
		\[
		I_j(\xi)\longrightarrow
		s\xi^2\int_\R|\zeta|^{2s-2}|\wh q_{s,H}(\zeta)|^2\dd\zeta
		=\beta_sH^{2-2s}\xi^2.
		\]
		The same estimates also give
		\[ |I_j(\xi)|\le s\xi^2\int_\R|\zeta|^{2s-2}|G_j|\dd\zeta
		\le C\xi^2.\]
		As $\xi^2|\widehat f\,\widehat g|\in L^1(\R)$,
		dominated convergence in $\xi$ proves
		\eqref{eq:product-limit} for $1/2<s<1$.
	\end{proof}
	
	\subsection{Minimal width and a full-height chord}
	For a convex
	body $K\subset\R^2$, define
	\[
	w(K)
	:=
	\min_{|e|=1}
	\left(
	\max_{x\in K}x\cdot e-\min_{x\in K}x\cdot e
	\right).
	\]
	
	\begin{lemma}\label{lem:width-chord}
		Let $K\subset\R^2$ be a convex body of area $A>0$
		and diameter $D$. Then
		\[
		\frac{A}{D}\le w(K)\le\frac{2A}{D}.
		\]
		If the vertical direction realizes $w(K)$ and the
		corresponding supporting lines are $z=0$ and $z=w(K)$,
		then there exists $x_0\in\R$ such that
		$\{(x_0,z):0\le z\le w(K)\}\subset K$.
	\end{lemma}
	
	\begin{proof}
		In a direction realizing the minimal width, $K$ is
		contained in a rectangle of height $w(K)$ and horizontal
		side at most $D$. Hence $A\le Dw(K).$
		Conversely, let $[p,q]$ be a diametral segment of $K$,
		so that $|p-q|=D$. Denote by $h_+$ and $h_-$ the maximal
		distances from $K$ to the line through $p,q$ on its
		two sides. By convexity, the two triangles with base
		$[p,q]$ and heights $h_+$ and $h_-$ lie in $K$ and
		have disjoint interiors. Thus
		$A\ge\frac D2(h_++h_-)\ge\frac D2w(K)$.
		The existence of the vertical chord follows from
		\cite[Statement~2.2]{Soltan2005}, which applies to
		convex bodies in arbitrary dimensions.
	\end{proof}

	\subsection{Convex compactness and form cores}
	
	We first establish a compactness result for nonnegative convex potentials. The main assertions follow from \cite[Proposition~7.2 and Theorems~7.6 and~7.17]{RockafellarWets1998}. For $V:\R\to[0,\infty]$, define
	\[
	\operatorname{dom}V:=\{x\in\R:V(x)<\infty\},
	\qquad
	\operatorname{epi}V:=\{(x,r)\in\R^2:r\ge V(x)\}.
	\]
	The function $V$ is proper if $\operatorname{dom}V\ne\varnothing$. Moreover, $V$ is lower semicontinuous if and only if $\operatorname{epi}V$ is closed, and convex if and only if $\operatorname{epi}V$ is convex.
	
	Recall that $V_j$ epi-converges to $V$ if
	\begin{enumerate}
		\item For every sequence $x_j\to x$,
		$V(x)\le\liminf_{j\to\infty}V_j(x_j)$.
		\item For every $x\in\operatorname{dom}V$, there exist $x_j\to x$ such that $V_j(x_j)\rightarrow V(x).$
	\end{enumerate}
	
	\begin{lemma}\label{lem:convex-compactness}
		Let $V_j:\R\to[0,\infty]$ be lower semicontinuous convex functions satisfying $V_j(0)=0$. Up to a subsequence, $V_j$ epi-converges to a proper lower semicontinuous convex function $V:\R\to[0,\infty]$ with $V(0)=0$.
		
		Writing $J:=\operatorname{int}\operatorname{dom}V$, we have
		\begin{equation}\label{eq:epi-local}
			\begin{aligned}
				V_j&\longrightarrow V
				&&\text{locally uniformly on }J,\\
				\inf_A V_j&\longrightarrow\infty
				&&\text{for every nonempty compact }
				A\subset\R\setminus\operatorname{dom}V.
			\end{aligned}
		\end{equation}
		If, for every $M\ge0$, there exists $R_M>0$ such that
		$\{V_j\le M\}\subset[-R_M,R_M], \enspace j\ge1$,
		then $\{V\le M\}\subset[-R_{M+1},R_{M+1}].$
		The two convergence assertions in \eqref{eq:epi-local} remain valid without the normalization $V_j(0)=0$, provided that $V_j$ epi-converges to $V$.
	\end{lemma}

	\begin{proof}
		Since $V_j\ge0$ and $V_j(0)=0$, the epigraphs $\operatorname{epi}V_j$ contain $(0,0)$ and lie in $\R\times[0,\infty)$. By \cite[Theorem~7.6]{RockafellarWets1998}, a subsequence epi-converges to a proper lower semicontinuous function $V\ge0$. By \cite[Theorem~7.17]{RockafellarWets1998}, $V$ is convex and
		\[
		V_j\longrightarrow V
		\qquad\text{locally uniformly on }
		J=\operatorname{int}\operatorname{dom}V.
		\]
		Moreover, \cite[Proposition~7.2]{RockafellarWets1998} gives the lower-limit and recovery-sequence conditions. In particular,
		$0\le V(0)\le\liminf_{j\to\infty}V_j(0)=0$,
		so $V(0)=0$. For every $x\in\operatorname{dom}V$, the recovery property gives $x_j\to x$ such that
		\[
		V(x)\le\liminf_{j\to\infty}V_j(x_j)
		\le\limsup_{j\to\infty}V_j(x_j)\le V(x).
		\]
		Thus $V_j(x_j)\to V(x)$.
		
		\smallskip
		
		We next prove the divergence assertion in \eqref{eq:epi-local}. Let $A\subset\R\setminus\operatorname{dom}V$ be nonempty and compact. If $\inf_A V_j\not\to\infty$, then, along a subsequence, there exist $x_j\in A$ and $C<\infty$ such that $V_j(x_j)\le C$. By compactness, we may assume $x_j\to x\in A$, thus,
		\[
		\infty=V(x)\le\liminf_{j\to\infty}V_j(x_j)\le C,
		\]
		which is a contradiction. Since neither this argument nor Theorem~7.17 requires $V_j(0)=0$ once epi-convergence is known, both assertions in \eqref{eq:epi-local} hold without this normalization.
		
		Finally, suppose that $\{V_j\le M\}\subset[-R_M,R_M]$ for every $j$ and $M\ge0$. For any $x\in\{V\le M\}$, choose recovery points $x_j\to x$ satisfying $V_j(x_j)\to V(x)$. For all sufficiently large $j$,
		$V_j(x_j)\le M+1, \enspace |x_j|\le R_{M+1}$.
		Passing to the limit yields $|x|\le R_{M+1}$, which completes the proof.
	\end{proof}

	Let $0<t\le1$, $a>0$, and let $J\subset\R$ be a nonempty open interval. For a nonnegative finite convex potential $V$ on $J$, define
	\begin{equation}\label{eq:one-d-form}
		Q_{t,a,V}[u]=a\int_\R|\xi|^{2t}|\wh u(\xi)|^2\dd\xi+\int_J V|u|^2\dd x,\quad
		\mathcal D(Q_{t,a,V})=\Hs{t}{J}\cap L^2(J,V\dd x).
	\end{equation}
	The form norm is $(\norm u_2^2+Q_{t,a,V}[u])^{1/2}$. For $t=1$, zero extension identifies $\Hs{1}{J}$ with $H_0^1(J)$, and the kinetic term is $a\int_J|u'|^2$.

\begin{lemma}\label{lem:form-core}
Let $J\subset\R$ be a nonempty open interval, possibly unbounded, and let $V:J\to[0,\infty)$ be a finite convex function. Then $C_c^\infty(J)$ is a form core for $Q_{t,a,V}$.
\end{lemma}

	\begin{proof}
		Fix $u\in \mathcal D(Q_{t,a,V})$ and $x_0\in J$.
		For $0<\rho<1$, define
		\[
		F_\rho(x):=x_0+\rho(x-x_0),
		\qquad
		(T_\rho u)(x):=\rho^{-1/2}u(F_\rho^{-1}(x)),
		\]
		where $u$ is extended by zero outside $J$.
		Since $F_\rho(J)\subset J$, we have
		$T_\rho u=0$ a.e. on $J^c$.
		
		The dilation satisfies
		\[
		\|T_\rho u\|_2=\|u\|_2,
		\qquad
		\int_\R|\xi|^{2t}|\widehat{T_\rho u}(\xi)|^2\,\dd\xi
		=
		\rho^{-2t}
		\int_\R|\xi|^{2t}|\widehat u(\xi)|^2\,\dd\xi.
		\]
		Since $T_\rho\varphi\to\varphi$ in $H^t(\R)$
		for every $\varphi\in C_c^\infty(\R)$,
		density and the uniform boundedness of $T_\rho$
		for $\rho$ near $1$ give, for every $u\in H^t(\R)$,
		\begin{equation}\label{eq:form-core-dilation}
			\|T_\rho u-u\|_{H^t(\R)}
			\longrightarrow0
			\qquad \rho\uparrow1.
		\end{equation}
		
		We next consider the potential term. Set $W=1+V$.
		Since $W\ge1$, convexity gives
		\[
		W(F_\rho(x))
		\le \rho W(x)+(1-\rho)W(x_0)
		\le W(x_0)W(x),
		\qquad x\in J.
		\]
		Therefore,
		\begin{equation}\label{eq:weighted-dilation-bound}
			\begin{aligned}
				\|T_\rho u\|_{L^2(J,W\dd x)}^2
				&=\int_J W(F_\rho(x))|u(x)|^2\,\dd x\le W(x_0)\|u\|_{L^2(J,W\dd x)}^2.
			\end{aligned}
		\end{equation}
		
		Since $V$ is finite and convex on $J$,
		it is continuous and locally bounded. Hence
		$C_c^\infty(J)$ is dense in $L^2(J,W\dd x)$.
		For every $\varphi\in C_c^\infty(J)$,
		the supports of $T_\rho\varphi$ remain in a
		fixed compact subset of $J$ for $\rho$ near $1$.
		Together with \eqref{eq:form-core-dilation},
		this gives
		$\|T_\rho\varphi-\varphi\|_{L^2(J,W\dd x)} \longrightarrow0$.
		For $u\in \mathcal D(Q_{t,a,V})$ and
		$\varphi\in C_c^\infty(J)$,
		\eqref{eq:weighted-dilation-bound} yields
		\[
		\begin{aligned}
			\|T_\rho u-u\|_{L^2(J,W\dd x)}
			&\le
			\|T_\rho(u-\varphi)\|_{L^2(J,W\dd x)}
			+\|T_\rho\varphi-\varphi\|_{L^2(J,W\dd x)}
			+\|\varphi-u\|_{L^2(J,W\dd x)}\\
			&\le
			\bigl(1+\sqrt{W(x_0)}\bigr)
			\|u-\varphi\|_{L^2(J,W\dd x)}
			+\|T_\rho\varphi-\varphi\|_{L^2(J,W\dd x)}.
		\end{aligned}
		\]
		Taking $\rho\uparrow1$ and then approximating
		$u$ by $\varphi\in C_c^\infty(J)$ gives
		$\|T_\rho u-u\|_{L^2(J,W\,\dd x)} \longrightarrow0$.
		In particular,
		\begin{equation}\label{eq:form-core-potential}
			\int_J V|T_\rho u-u|^2\,\dd x
			\longrightarrow0.
		\end{equation}
		Combining \eqref{eq:form-core-dilation}
		and \eqref{eq:form-core-potential}, we obtain
		\begin{equation}\label{eq:form-core-contraction}
			\|T_\rho u-u\|_{Q_{t,a,V}}
			\longrightarrow0.
		\end{equation}
		
		Fix $0<\rho<1$ and put $v=T_\rho u$.
		Since $\operatorname{supp}v\subset F_\rho(\overline J)$,
		and
		$|F_\rho(b)-b| =(1-\rho)|x_0-b|>0$
		for every finite endpoint $b$ of $J$,
		$v$ vanishes near each such endpoint. Choose $\chi\in C_c^\infty(\R)$ satisfying
		\[
		0\le\chi\le1,\qquad
		\chi=1\ \text{on }[-1,1],\qquad
		\chi=0\ \text{outside }[-2,2],
		\]
		and set $\chi_R(x)=\chi(x/R)$.
		The standard cutoff approximation gives $	\|\chi_Rv-v\|_{H^t(\R)}
		\rightarrow0.$
		Moreover, dominated convergence yields
		$\int_J V|(1-\chi_R)v|^2\,\dd x \longrightarrow0$.
		Consequently, $\operatorname{supp}(\chi_Rv)\Subset J$ and
		\begin{equation}\label{eq:form-core-infinity}
			\|\chi_Rv-v\|_{Q_{t,a,V}}
			\longrightarrow0.
		\end{equation}
		
		Finally, set $f=\chi_Rv$ and let
		$\eta_\varepsilon$ be a standard mollifier.
		For sufficiently small $\varepsilon>0$,
		\[
		f_\varepsilon:=\eta_\varepsilon*f
		\in C_c^\infty(J),
		\qquad
		\operatorname{supp}f_\varepsilon\subset K
		\]
		for some fixed compact set $K\Subset J$.
		By the approximation property of mollifiers,
		$\|f_\varepsilon-f\|_{H^t(\R)}\to0$.
		Since $V$ is bounded on $K$,
		\[
		\begin{aligned}
			\int_J V|f_\varepsilon-f|^2\,\dd x
			&\le
			\|V\|_{L^\infty(K)}
			\|f_\varepsilon-f\|_2^2
			\longrightarrow0.
		\end{aligned}
		\]
		Thus, $\|f_\varepsilon-f\|_{Q_{t,a,V}}\to0$.
		Together with \eqref{eq:form-core-contraction}
		and \eqref{eq:form-core-infinity}, this completes the proof.
	\end{proof}

	\subsection{Compact form convergence and eigenvalues}
	
	We first establish a spectral convergence result for nonnegative quadratic forms. 
	For the general theory of compact Mosco convergence and spectral convergence, see \cite[Definition~2.13 and Theorem~2.6]{KuwaeShioya2003}. We give a direct proof in a fixed Hilbert space, allowing the limiting form to be densely defined only on a proper closed subspace.

	\begin{lemma}\label{lem:abstract-spectral}
		Let $Q_j$ and $Q$ be nonnegative closed quadratic forms on a Hilbert space $\mathcal H$, extended by $+\infty$ outside their domains. Suppose that each $\mathcal D(Q_j)$ is compactly embedded in $\mathcal H$ and that the following conditions hold:
		\begin{enumerate}
			\item Every sequence $(u_j)$ satisfying
			$\sup_j\bigl(\|u_j\|_{\mathcal H}^2+Q_j[u_j]\bigr)<\infty$
			is precompact in $\mathcal H$.
			
			\item For every sequence $u_j\rightharpoonup u$ weakly in $\mathcal H$,
			$Q[u]\le\liminf_{j\to\infty}Q_j[u_j]$.
			
			\item For every $u\in \mathcal D(Q)$, there exist $u_j\to u$ strongly in $\mathcal H$ such that
			$Q_j[u_j]\longrightarrow Q[u]$.
		\end{enumerate}
		Then $\mathcal D(Q)$ is compactly embedded in $\mathcal H$. Writing
		\[
		\mathcal H_j:=\overline{\mathcal D(Q_j)}^{\,\mathcal H},
		\qquad
		\mathcal H_\infty:=\overline{\mathcal D(Q)}^{\,\mathcal H}.
		\]
		For every $k\ge1$ with $\dim\mathcal H_\infty\ge k$, their eigenvalues satisfy
		\begin{equation}\label{eq:abstract-eigen-convergence}
			\lambda_k(Q_j)\longrightarrow\lambda_k(Q).
		\end{equation}
		If $\mathcal D(Q)=\{0\}$, $\lambda_1(Q_j)\to+\infty$ whenever the first eigenvalues exist for sufficiently large $j$.
	\end{lemma}
	
	\begin{proof}
		\textbf{Step 1. Compactness of the limiting form.}
		Let $(u_n)\subset \mathcal D(Q)$ satisfy
		${\sup_n\bigl(\|u_n\|_{\mathcal H}^2+Q[u_n]\bigr)<\infty}$.
		By the recovery property, for each $n$, we may choose $j_n>j_{n-1}$ and $v_n\in \mathcal D(Q_{j_n})$ such that
		$\|v_n-u_n\|_{\mathcal H}\le\frac1n, \enspace Q_{j_n}[v_n]\le Q[u_n]+1$.
		Define $w_j=v_n$ when $j=j_n$, and $w_j=0$ otherwise. Then
		$\sup_j\bigl(\|w_j\|_{\mathcal H}^2+Q_j[w_j]\bigr)<\infty$.
		By assumption (i), $(w_j)$ is precompact in $\mathcal H$. In particular, $(v_n)$ has a strongly convergent subsequence. Since $\|v_n-u_n\|_{\mathcal H}\to0$, the same holds for $(u_n)$. Thus $\mathcal D(Q)$ is compactly embedded in $\mathcal H$. Since $Q$ is closed, its restriction to $\mathcal H_\infty$ is a densely defined closed form with compact resolvent. The same conclusion holds for $Q_j$ on $\mathcal H_j$.
		
		\textbf{Step 2. Upper bound for the eigenvalues.}
		Fix $k\ge1$ with $\dim\mathcal H_\infty\ge k$, and let $e_1,\ldots,e_k$ be orthonormal eigenfunctions of $Q$ corresponding to its first $k$ eigenvalues. By assumption (iii), for each $1\le m\le k$ there exist $e_{m,j}\in \mathcal D(Q_j)$ such that
		\[
		e_{m,j}\longrightarrow e_m
		\quad\text{in }\mathcal H,
		\qquad
		Q_j[e_{m,j}]\longrightarrow Q[e_m]=\lambda_m(Q).
		\]
		
		We first establish convergence of the mixed form values. Suppose that $u_j\to u$ and $v_j\to v$ strongly in $\mathcal H$, with $Q_j[u_j]\to Q[u],\,
		Q_j[v_j]\to Q[v].$ Set
		\[
		A_j:=Q_j[u_j],\qquad
		B_j:=Q_j[v_j],\qquad
		C_j:=\operatorname{Re}Q_j[u_j,v_j].
		\]
		The Cauchy--Schwarz inequality for nonnegative quadratic forms gives $|C_j|\le A_j^{1/2}B_j^{1/2},$
		so $(C_j)$ is bounded. Applying assumption (ii) to $u_j+v_j$ and $u_j-v_j$, we obtain
		\[Q[u+v]\le\liminf_{j\to\infty}(A_j+B_j+2C_j)=Q[u]+Q[v]+2\liminf_{j\to\infty}C_j,\]
		and
		\[
		Q[u-v]\le\liminf_{j\to\infty}(A_j+B_j-2C_j)=Q[u]+Q[v]-2\limsup_{j\to\infty}C_j.
		\]
		Hence,
		\[
		\operatorname{Re}Q[u,v]
		\le\liminf_{j\to\infty}C_j
		\le\limsup_{j\to\infty}C_j
		\le\operatorname{Re}Q[u,v].
		\]
		For complex forms, the same argument applied to $u_j\pm iv_j$ gives convergence of the imaginary parts. Consequently, $Q_j[u_j,v_j]\rightarrow Q[u,v].$
		In particular,
		\[
		Q_j[e_{m,j},e_{n,j}]
		\longrightarrow Q[e_m,e_n]
		=\lambda_m(Q)\delta_{mn},
		\qquad 1\le m,n\le k.
		\]
		
		Define the Gram and form matrices
		\[
		G_j:=\bigl(\langle e_{m,j},e_{n,j}\rangle_{\mathcal H}\bigr)_{m,n=1}^k,
		\qquad
		B_j:=\bigl(Q_j[e_{m,j},e_{n,j}]\bigr)_{m,n=1}^k.
		\]
		The preceding convergence and strong convergence of $e_{m,j}$ give
		\begin{equation}\label{eq:abstract-form-matrices}
			\begin{aligned}
				\|G_j-I_k\|_{\mathrm{op}}\longrightarrow0,\\
				\bigl\|B_j-\operatorname{diag}
				(\lambda_1(Q),\ldots,\lambda_k(Q))\bigr\|_{\mathrm{op}}
				\longrightarrow0,
			\end{aligned}
		\end{equation}
		where $\|\cdot\|_{\mathrm{op}}$ denotes the matrix operator norm induced by the Euclidean norm on $\mathbb C^k$.
		
		For all sufficiently large $j$, the matrix $G_j$ is positive definite, and hence
		\[
		E_j:=\operatorname{span}\{e_{1,j},\ldots,e_{k,j}\}
		\]
		has dimension $k$. For every nonzero $c=(c_1,\ldots,c_k)\in\mathbb C^k$, set $u_j=\sum_{m=1}^k c_me_{m,j}$. By \eqref{eq:abstract-form-matrices},
		\[
		\frac{Q_j[u_j]}{\|u_j\|_{\mathcal H}^2}
		=\frac{c^*B_jc}{c^*G_jc}
		\le\frac{\lambda_k(Q)+o(1)}{1-o(1)},
		\]
		uniformly for $c\ne0$. The min--max principle therefore yields
		\begin{equation}\label{eq:abstract-upper}
			\limsup_{j\to\infty}\lambda_k(Q_j)
			\le\lambda_k(Q).
		\end{equation}
		
		\textbf{Step 3. Lower bound for the eigenvalues.}
		Choose a subsequence along which
		\[
		\lambda_k(Q_j)\longrightarrow
		\ell:=\liminf_{j\to\infty}\lambda_k(Q_j)<\infty.
		\]
		Let $v_{1,j},\ldots,v_{k,j}$ be orthonormal eigenfunctions corresponding to the first $k$ eigenvalues of $Q_j$. Then
		\[
		\|v_{m,j}\|_{\mathcal H}=1,
		\qquad
		Q_j[v_{m,j}]
		=\lambda_m(Q_j)\le\lambda_k(Q_j),
		\quad 1\le m\le k.
		\]
		By assumption (i), after passing to a common subsequence,
		\[
		v_{m,j}\longrightarrow v_m
		\quad\text{strongly in }\mathcal H,
		\qquad 1\le m\le k.
		\]
		The limit functions satisfy $\langle v_m,v_n\rangle_{\mathcal H}
		=\delta_{mn}.$
		For any $c=(c_1,\ldots,c_k)\in\mathbb C^k$, assumption (ii) gives
		\[
		\begin{aligned}
			Q\left[\sum_{m=1}^k c_mv_m\right]
			&\le\liminf_{j\to\infty}
			Q_j\left[\sum_{m=1}^k c_mv_{m,j}\right]=\liminf_{j\to\infty}
			\sum_{m=1}^k\lambda_m(Q_j)|c_m|^2\le\ell\sum_{m=1}^k|c_m|^2.
		\end{aligned}
		\]
		Thus $v_1,\ldots,v_k\in \mathcal D(Q)$, and their span has dimension $k$. The min--max principle yields
		\[
		\lambda_k(Q)
		\le
		\sup_{\substack{v\in\operatorname{span}\{v_1,\ldots,v_k\},\,v\ne0}}
		\frac{Q[v]}{\|v\|_{\mathcal H}^2}
		\le\ell.
		\]
		Together with \eqref{eq:abstract-upper}, this proves \eqref{eq:abstract-eigen-convergence}.
		
		\textbf{Step 4. Degenerate limiting form.}
		Suppose that $\mathcal D(Q)=\{0\}$. Let $(u_j)$ satisfy
		$\sup_j\bigl(\|u_j\|_{\mathcal H}^2+Q_j[u_j]\bigr)<\infty$.
		By assumption (i), every subsequence has a strongly convergent further subsequence. If $u$ is such a limit, assumption (ii) gives $Q[u]\le\liminf_{j\to\infty}Q_j[u_j]<\infty.$
		Hence $u\in \mathcal D(Q)=\{0\}$. Since every cluster point is zero, $u_j\to0$ strongly in $\mathcal H$. 
		If $\lambda_1(Q_j)$ did not tend to infinity, there would exist a subsequence and normalized first eigenfunctions $u_j$ satisfying
		$\|u_j\|_{\mathcal H}=1, \enspace Q_j[u_j]=\lambda_1(Q_j)\le C$,
		contradicting $\|u_j\|_{\mathcal H}=1$. This completes the proof.
	\end{proof}
	
	We now apply the preceding lemma to one-dimensional fractional Schrödinger forms with convex potentials. For a proper nonnegative lower semicontinuous convex function $V:\R\to[0,\infty]$, set $J:=\operatorname{int}\operatorname{dom}V$ and use the form $Q_{t,a,V}$ defined in \eqref{eq:one-d-form}, with $\mathcal D(Q_{t,a,V})=\{0\}$ when $J=\varnothing$. Since $\operatorname{dom}V\setminus J$ contains at most two points, the zero-extension condition is equivalent to $u=0$ almost everywhere on $\{V=+\infty\}$. In integrals involving extended-valued potentials, we set $V|u|^2=0$ wherever $u=0$.
	
	\begin{proposition}\label{prop:potential-convergence}
		Let $0<t\le1$, $a_j\to a>0$, and let $V_j,V:\R\to[0,\infty]$ be proper lower semicontinuous convex functions satisfying $V_j\econv V$. Suppose that either
		\begin{enumerate}
			\item all $\operatorname{dom}V_j$ are contained in a fixed bounded interval; or
			\item for every $M>0$, there exists $R_M>0$ such that
			$\{V_j\le M\}\subset[-R_M,R_M], \enspace j\ge1$.
		\end{enumerate}
		Set $Q_j:=Q_{t,a_j,V_j}$ and $Q:=Q_{t,a,V}$. Then the three assumptions of Lemma~\ref{lem:abstract-spectral} hold. If $J:=\operatorname{int}\operatorname{dom}V\ne\varnothing$, the effective intervals of $V_j$ are nonempty for all sufficiently large $j$, and
		\begin{equation}\label{eq:potential-eigen-convergence}
			\lambda_k(Q_j)\longrightarrow\lambda_k(Q),
			\qquad k\ge1.
		\end{equation}
		If $J=\varnothing$, every sequence with uniformly bounded $L^2$ norms and $Q_j$-energies converges strongly to zero, and $\lambda_1(Q_j)\to+\infty$ whenever the first eigenvalues exist for all sufficiently large $j$.
	\end{proposition}
	
	\begin{proof}
		\textbf{Step 1. Compactness.}
		Let $(u_j)$ satisfy
		$\sup_j\bigl(\|u_j\|_2^2+Q_j[u_j]\bigr)<\infty$.
		Since $a_j\to a>0$ and $V_j\ge0$, we have $\sup_j\|u_j\|_{H^t(\R)}<\infty.$
		By Plancherel's identity,
		\[
		\begin{aligned}
			\|u_j(\cdot+h)-u_j\|_2^2=\int_\R|e^{i\xi h}-1|^2|\wh u_j(\xi)|^2\,\dd\xi\le C|h|^{2t}.
		\end{aligned}
		\]
		
		Under assumption (i), all $u_j$ are supported in a fixed bounded interval. Under assumption (ii), for every $M>0$,
		\[
		\int_{|x|>R_M}|u_j(x)|^2\,\dd x
		\le\frac1M\int_\R V_j(x)|u_j(x)|^2\,\dd x
		\le\frac CM.
		\]
		
		The Fr\'echet--Kolmogorov criterion therefore implies that $(u_j)$ is precompact in $L^2(\R)$. The same estimates, with $j$ fixed, give compactness of the individual form embeddings.
		
		\textbf{Step 2. Lower-limit inequality.}
		Let $u_j\rightharpoonup u$ weakly in $L^2(\R)$. We may assume $\liminf_{j\to\infty}Q_j[u_j]<\infty.$
		Passing to a subsequence realizing this lower limit, Step 1 gives
		\[
		u_j\longrightarrow u\quad\text{in }L^2(\R),
		\qquad
		u_j(x)\longrightarrow u(x)\quad\text{a.e.},
		\]
		and bounded kinetic energies imply $u_j\rightharpoonup u$ in $H^t(\R).$ The epi-liminf inequality gives
		$V(x)\le\liminf_{j\to\infty}V_j(x), \enspace x\in\R$.
		Since $V_j\ge0$, Fatou's lemma yields
		\[
		\int_\R V(x)|u(x)|^2\,\dd x
		\le\liminf_{j\to\infty}
		\int_\R V_j(x)|u_j(x)|^2\,\dd x.
		\]
		In particular, $u=0$ almost everywhere on $\{V=+\infty\}$. Since $\operatorname{dom}V\setminus J$ has measure zero, the zero-extension condition holds on $J^c$. Moreover, weak lower semicontinuity gives $\Es_t[u]\le\liminf_{j\to\infty}\Es_t[u_j].$
		Since $\sup_j\Es_t[u_j]<\infty$ and $a_j\to a$,
		$|(a_j-a)\Es_t[u_j]| \le|a_j-a|\sup_j\Es_t[u_j]\longrightarrow0$.
		Combining the kinetic and potential terms, we obtain $Q[u]\le\liminf_{j\to\infty}Q_j[u_j].$
		
		\textbf{Step 3. Recovery sequences.}
		Suppose first that $J\ne\varnothing$. For every $\varphi\in C_c^\infty(J)$, the local uniform convergence of convex potentials in Lemma~\ref{lem:convex-compactness} gives
		$\sup_{\operatorname{supp}\varphi}|V_j-V| \longrightarrow0$.
		In particular, $\varphi\in \mathcal D(Q_j)$ for all sufficiently large $j$, and
		\[
		\begin{aligned}
			|Q_j[\varphi]-Q[\varphi]|
			&\le |a_j-a|\Es_t[\varphi]
			+\int_J|V_j-V||\varphi|^2\,\dd x\\
			&\le |a_j-a|\Es_t[\varphi]
			+\|\varphi\|_2^2
			\sup_{\operatorname{supp}\varphi}|V_j-V|
			\longrightarrow0.
		\end{aligned}
		\]
		
		For arbitrary $u\in \mathcal D(Q)$, Lemma~\ref{lem:form-core} gives $\varphi_n\in C_c^\infty(J)$ such that
		\[
		\|\varphi_n-u\|_2\longrightarrow0,
		\qquad
		Q[\varphi_n]\longrightarrow Q[u].
		\]
		Choose increasing integers $j_n$ such that
		$|Q_j[\varphi_n]-Q[\varphi_n]| \le\frac1n, \enspace j\ge j_n$.
		Define $u_j:=\varphi_n,\,j_n\le j<j_{n+1}.$
		Then $u_j\to u$ strongly in $L^2(\R)$, and
		\[
		\begin{aligned}
			|Q_j[u_j]-Q[u]|
			&\le |Q_j[\varphi_n]-Q[\varphi_n]|
			+|Q[\varphi_n]-Q[u]|\\
			&\le\frac1n+|Q[\varphi_n]-Q[u]|
			\longrightarrow0.
		\end{aligned}
		\]
		If $J=\varnothing$, then $\mathcal D(Q)=\{0\}$, and $u_j=0$ is a recovery sequence.
		
		\textbf{Step 4. Spectral convergence.}
		Steps 1--3 verify the three assumptions of Lemma~\ref{lem:abstract-spectral}. Hence the limiting form is compactly embedded in $L^2(\R)$ on the closure of its domain.
		
		If $J\ne\varnothing$, then $C_c^\infty(J)\subset \mathcal D(Q)$ is dense in $L^2(J)$. Thus the limiting form is densely defined on $L^2(J)$, and its associated operator has compact resolvent. Since $J$ is a nonempty open interval, the form domain is infinite-dimensional. Lemma~\ref{lem:abstract-spectral} therefore gives \eqref{eq:potential-eigen-convergence} for every $k\ge1$. If $J=\varnothing$, the degenerate-limit conclusion of the same lemma gives the asserted strong convergence to zero and divergence of the first eigenvalues.
	\end{proof}

	\begin{corollary}\label{cor:ground-simple}
		Let $0<t\le1$, $a>0$, and let $V:\R\to[0,\infty]$ be a proper lower semicontinuous convex function with $J:=\operatorname{int}\dom V\ne\varnothing$. Suppose that either $\dom V$ is bounded or every sublevel set of $V$ is bounded. Let $Q_{t,a,V}$ be defined in \eqref{eq:one-d-form}, then $Q_{t,a,V}$ has compact resolvent, and its eigenvalues satisfy
		\[
		0<\lambda_1(Q_{t,a,V})<\lambda_2(Q_{t,a,V})\le\lambda_3(Q_{t,a,V})\le\cdots,
		\qquad \lambda_k(Q_{t,a,V})\longrightarrow\infty.
		\]
		Moreover, for any bounded open interval $J$ and finite
		convex potential $V:J\to\R$, the form $Q_{t,a,V}$
		has compact resolvent and a simple first eigenvalue.
	\end{corollary}
	
	\begin{proof}
		Compactness follows from the proof in Proposition~\ref{prop:potential-convergence}. 
		Let $\phi_1$ be an $L^2$-normalized first eigenfunction. Since $V\ge0$ and $a>0$,
		$\lambda_1(Q_{t,a,V})=Q_{t,a,V}[\phi_1] \ge a\Es_t[\phi_1]>0$.
		
		Suppose first that $0<t<1$. If a real first eigenfunction $u$ changes sign, then \eqref{eq:absolute-value} gives
		\[
		Q_{t,a,V}[|u|]<Q_{t,a,V}[u]
		=\lambda_1(Q_{t,a,V})\|u\|_2^2
		=\lambda_1(Q_{t,a,V})\||u|\|_2^2,
		\]
		contradicting the variational characterization of $\lambda_1(Q_{t,a,V})$. Thus every real first eigenfunction has a fixed sign. If $\lambda_1$ were multiple, there would exist two orthogonal nonnegative first eigenfunctions $u,v$, we have
		$uv=0$ a.e. on $J$. As $u,v$ are nonzero, their
		difference changes sign on sets of positive measure.
		By linearity, $u-v$ is also a first eigenfunction,
		a contradiction.

		For $t=1$, let $u$ be a real first eigenfunction. Since $Q_{1,a,V}[|u|]=Q_{1,a,V}[u]$, the function $w:=|u|$ is also a first eigenfunction and satisfies
		$-aw''+(V-\lambda_1)w=0 \enspace\text{in }J$.
		Since $V\in C(J)$, one-dimensional regularity gives $w\in C^2(J)$. If $w(x_0)=0$ for some $x_0\in J$, then
		$w(x_0)=w'(x_0)=0 \quad\Longrightarrow\quad w\equiv0$
		by uniqueness for the initial-value problem, a contradiction. Hence $w>0$ in $J$. Every real first eigenfunction therefore has a fixed sign, and orthogonality proves that $\lambda_1$ is simple.
		
		Finally, let $J$ be a bounded open interval and $V:J\to\R$ be finite and convex. Since $V$ admits a supporting affine lower bound, there exists $c\ge0$ such that $V+c\ge0$ on $J$. Extend $V+c$ to a lower semicontinuous convex function $\widetilde V:\R\to[0,\infty]$ by its one-sided limits at the endpoints and by $+\infty$ outside $\overline J$. Then $\operatorname{int}\operatorname{dom}\widetilde V=J$, and the preceding conclusions apply. Since the form domains coincide and
		$Q_{t,a,\widetilde V}[u] =Q_{t,a,V}[u]+c\|u\|_2^2$,
		the associated operators satisfy $H_{t,a,\widetilde V}=H_{t,a,V}+cI.$
		Thus compactness of the resolvent and simplicity of the first eigenvalue also hold for $V$.
	\end{proof}

	\subsection{Uniform gaps and attainment for convex potentials}

We first prove the following gap estimate for finite convex potentials on arbitrary bounded intervals.

	\begin{lemma}\label{lem:convex-gap-bound}
Let $0<t\le1$, $a>0$, and let $J$ be a bounded open interval of length $L$. For every finite convex potential $V:J\to\R$, the first two eigenvalues of $Q_{t,a,V}$ satisfy
\begin{equation}\label{eq:convex-gap-bound}
\lambda_2(Q_{t,a,V})-\lambda_1(Q_{t,a,V})
\ge
\begin{cases}
aL^{-2t}\gs_t,&0<t<1,\\
3\pi^2aL^{-2},&t=1.
\end{cases}
\end{equation}
\end{lemma}
	
	\begin{proof}
		\step{Step 1. Continuous convex potentials.}
		Suppose first that $V\in C(\overline J)$. Write $J=x_0+LI$ and define
		$W(y):=a^{-1}L^{2t}V(x_0+Ly), \enspace y\in I$.
		The unitary rescaling gives
		$\lambda_k(Q_{t,a,V}) =aL^{-2t}\lambda_{k,t}(I,W), \enspace k\ge1$.
		For $0<t<1$, \eqref{eq:convex-gap-bound} follows from the definition of $\gs_t$. For $t=1$, the result follows from \cite[Theorem~3.1]{Lavine1994}, applied to $W$ with Dirichlet boundary conditions.
		
		\step{Step 2. Finite convex potentials.}
		Let $V:J\to\R$ be finite and convex. Since $J$ is bounded, a supporting affine function gives a lower bound for $V$. Adding a constant, we may assume $V\ge0$. 
Choose a dense sequence $(x_n)\subset J$ and slopes $p_n\in\partial V(x_n)$, where the subdifferential of $V$ at $x\in J$ is defined by
\[
\partial V(x):=
\left\{p\in\R:
V(y)\ge V(x)+p(y-x)
\text{ for all }y\in J
\right\}.
\]
Define
\[
\ell_n(x):=V(x_n)+p_n(x-x_n),
\qquad
V_m(x):=\max\{0,\ell_1(x),\ldots,\ell_m(x)\}.
\]
Each $V_m$ is continuous and convex on $\overline J$, and
		$0\le V_m\le V_{m+1}\le V \enspace\text{on }J$.
		
		We show that $V_m(x)\uparrow V(x)$ for every $x\in J$. Fix $x\in J$ and choose a subsequence $x_{n_j}\to x$. Since $V$ is locally Lipschitz, its subgradients are locally bounded. Hence
		\[
		\begin{aligned}
			|\ell_{n_j}(x)-V(x)|
			&\le |V(x_{n_j})-V(x)|
			+|p_{n_j}|\,|x-x_{n_j}|\longrightarrow0.
		\end{aligned}
		\]
		Since $\ell_n\le V$ and $V_m$ is increasing, we obtain
		$V_m(x)\uparrow V(x), \enspace x\in J$.
		
		We next prove convergence of the corresponding eigenvalues. Fix $k\ge1$. By the min--max principle,
		\[
		\lambda_k(Q_{t,a,V_m})\uparrow b_k
		\le\lambda_k(Q_{t,a,V})<\infty,\quad 
		b_k:=
		\lim_{m\to\infty}\lambda_k(Q_{t,a,V_m}).
		\]
		Let $u_{m,1},\ldots,u_{m,k}$ be orthonormal eigenfunctions corresponding to the first $k$ eigenvalues of $Q_{t,a,V_m}$. Since $V_m\ge0$,
		\[
		a\int_\R|\xi|^{2t}|\wh u_{m,\ell}(\xi)|^2\,\dd\xi
		\le\lambda_\ell(Q_{t,a,V_m})
		\le b_k,
		\qquad 1\le\ell\le k.
		\]
		The compact embedding $\Hs{t}{J}\hookrightarrow L^2(J)$ gives, after passing to a common subsequence,
		\[
		u_{m,\ell}\longrightarrow u_\ell
		\quad\text{strongly in }L^2(J)
		\quad\text{and weakly in }H^t(\R),
		\qquad 1\le\ell\le k.
		\]
		In particular, $\langle u_\ell,u_r\rangle_{L^2(J)}
		=\delta_{\ell r}.$ Fix $n\ge1$. For all $m\ge n$ and every $c=(c_1,\ldots,c_k)\in\C^k$,
		\[
		\begin{aligned}
			Q_{t,a,V_n}\left[\sum_{\ell=1}^kc_\ell u_{m,\ell}\right]
			&\le
			Q_{t,a,V_m}\left[\sum_{\ell=1}^kc_\ell u_{m,\ell}\right]=\sum_{\ell=1}^k
			\lambda_\ell(Q_{t,a,V_m})|c_\ell|^2\le b_k\sum_{\ell=1}^k|c_\ell|^2.
		\end{aligned}
		\]
		By weak lower semicontinuity of $Q_{t,a,V_n}$ and $V_n\uparrow V$,
		\[
		Q_{t,a,V}\left[\sum_{\ell=1}^kc_\ell u_\ell\right]
		\le b_k\sum_{\ell=1}^k|c_\ell|^2.
		\]
		Thus $\operatorname{span}\{u_1,\ldots,u_k\}$ is a $k$-dimensional subspace of $\mathcal D(Q_{t,a,V})$. The min--max principle yields
		$\lambda_k(Q_{t,a,V}) \le b_k =\lim_{m\to\infty}\lambda_k(Q_{t,a,V_m})$.
		Combining this with the reverse inequality gives
		\begin{equation}\label{eq:monotone-potential-spectrum}
			\lambda_k(Q_{t,a,V_m})
			\longrightarrow\lambda_k(Q_{t,a,V}),
			\qquad k\ge1.
		\end{equation}
		
		Applying Step~1 to $V_m$ and passing to the limit proves \eqref{eq:convex-gap-bound} for every finite convex potential. Moreover, \eqref{eq:monotone-potential-spectrum} shows that the infimum defining $\gs_t$ is unchanged when continuous convex potentials on $\overline I$ are replaced by all finite convex potentials on $I$.
	\end{proof}

        We next establish a uniform lower bound for the fundamental gap in terms of the first eigenvalue and the minimum of the convex potential. For every $0<t<1,$ recall that
	\[ \gs_t=\inf\{\gap_t(I,V):V\in C(\overline I)\text{ convex}\}.\]

\begin{proof}[\textbf{Proof of
			Proposition~\ref{prop:relative-gap}}.]
		Subtracting the minimum of $V$ does not change the gap. Suppose that no positive constant satisfies \eqref{eq:relative-gap} for nonnegative convex potentials with minimum zero. Then there exist such $V_j\in C(\overline I)$ satisfying
		\[
		e_j:=\lambda_{1,t}(I,V_j)>0,
		\qquad
		\frac{\gap_t(I,V_j)}{e_j}\longrightarrow0.
		\]
		Choose $x_j\in\overline I$ with $V_j(x_j)=0$, and set
		\[
		r_j:=e_j^{-1/(2t)},\qquad
		J_j:=r_j^{-1}(I-x_j),\qquad
		U_j(y):=e_j^{-1}V_j(x_j+r_jy)
		\]
		for $y\in\overline J_j$, with $U_j=+\infty$ elsewhere. Then $U_j$ is nonnegative, lower semicontinuous, and convex, with $U_j(0)=0$. The unitary rescaling
		$v(y)=r_j^{1/2}u(x_j+r_jy)$ gives
		\begin{equation}\label{eq:relative-rescaled-eigenvalues}
			\lambda_1(Q_{t,1,U_j})=1,
			\qquad
			\lambda_2(Q_{t,1,U_j})
			=1+\frac{\gap_t(I,V_j)}{e_j}
			\longrightarrow1.
		\end{equation}
		
		\step{Step 1. Uniform confinement.}
		Fix $\tau=1/2$ and $\varphi\in C_c^\infty((0,1))$ with $\|\varphi\|_2=1$. Set
		$A_\varphi:=\int_\R|\xi|^{2t}|\wh\varphi(\xi)|^2\,\dd\xi$.
		Choose $R>0$ sufficiently large that $R^{-2t}A_\varphi+\tau<1.$
		The sublevel set
		$S_j:=\{y\in\R:U_j(y)\le\tau\}$
		is a closed interval containing zero, possibly a singleton. We claim that $|S_j|\le R$. Otherwise, there exists $(b,b+R)\subset S_j$, and the normalized function
		\[
		\varphi_{R,b}(y)
		:=R^{-1/2}\varphi\left(\frac{y-b}{R}\right)
		\]
		belongs to $\mathcal D(Q_{t,1,U_j})$. By the min--max principle,
		\[
		\begin{aligned}
			1=\lambda_1(Q_{t,1,U_j})
			&\le Q_{t,1,U_j}[\varphi_{R,b}]\le R^{-2t}A_\varphi+\tau<1,
		\end{aligned}
		\]
		a contradiction. Hence $S_j\subset[-R,R]$.
		
		\smallskip
		
		For $|y|\ge2R$ with $U_j(y)<\infty$, convexity and $U_j(0)=0$ yield
		$\tau<U_j(2R\operatorname{sgn}y) \le\frac{2R}{|y|}U_j(y)$.
		Since $U_j\ge0$, we obtain
		\begin{equation}\label{eq:relative-confinement}
			U_j(y)\ge\frac{\tau}{2R}|y|-\tau,
			\qquad y\in\R.
		\end{equation}
		Consequently, for every $M\ge0$,
		\[
		\{U_j\le M\}
		\subset
		\left[
		-\frac{2R(M+\tau)}{\tau},
		\frac{2R(M+\tau)}{\tau}
		\right].
		\]
		Thus the sublevel sets of $U_j$ are uniformly bounded.
		
		\step{Step 2. Contradiction by spectral convergence.}
		By Lemma~\ref{lem:convex-compactness}, after passing to a subsequence,
		$U_j\econv U$ for some proper nonnegative lower semicontinuous convex function $U$ with $U(0)=0$. For every $y\in\dom U$, choose recovery points $y_j\to y$ with $U_j(y_j)\to U(y)$. Passing to the limit in \eqref{eq:relative-confinement} gives
		$U(y)\ge\frac{\tau}{2R}|y|-\tau, \enspace y\in\R$.
		In particular, every sublevel set of $U$ is bounded.
		
		Set $J:=\operatorname{int}\dom U$. If $J=\varnothing$, Proposition~\ref{prop:potential-convergence} gives $\lambda_1(Q_{t,1,U_j})\rightarrow\infty,$
		contradicting \eqref{eq:relative-rescaled-eigenvalues}. Hence $J\ne\varnothing$. Applying the same proposition with $a_j=a=1$, we obtain
		\[
		\lambda_1(Q_{t,1,U})
		=\lim_{j\to\infty}\lambda_1(Q_{t,1,U_j})=1,\quad
		\lambda_2(Q_{t,1,U})
		=\lim_{j\to\infty}\lambda_2(Q_{t,1,U_j})=1.
		\]
		This contradicts the simplicity of the first eigenvalue established in Corollary~\ref{cor:ground-simple}. Therefore, the infimum of the ratios in \eqref{eq:relative-gap} is strictly positive.
		
		\smallskip
		
		For any convex $V\in C(\overline I)$, set
		$W:=V-\min_{\overline I}V$. Then
		\[
		\gap_t(I,V)=\gap_t(I,W),
		\qquad
		\lambda_{1,t}(I,V)-\min_{\overline I}V
		=\lambda_{1,t}(I,W).
		\]
		Since $W\ge0$,
		\[
		\gap_t(I,V)
		\ge c_t\lambda_{1,t}(I,W)
		\ge c_t\lambda_{1,t}(I,0).
		\]
		Consequently, $\gs_t\ge c_t\lambda_{1,t}(I,0)>0.$

        \step{Step 3. Attainment of the minimum.}
		Choose convex potentials $V_j\in C(\overline I)$ satisfying
		\[
		V_j\ge0,\qquad
		\min_{\overline I}V_j=0,\qquad
		\gap_t(I,V_j)\longrightarrow\gs_t.
		\]
		By \eqref{eq:relative-gap}, we obtain
		\begin{equation}\label{eq:minimizing-eigen-bound}
			\sup_j\lambda_{k,t}(I,V_j)<\infty,
			\qquad k=1,2.
		\end{equation}
		
		Choose $x_j\in\overline I$ with $V_j(x_j)=0$. After passing to a subsequence, assume $x_j\to x_*\in\overline I$. Define
		\[
		U_j(y):=
		\begin{cases}
			V_j(x_j+y), & y\in\overline I-x_j,\\
			+\infty, & y\notin\overline I-x_j.
		\end{cases}
		\]
		Then $U_j$ is nonnegative, lower semicontinuous, and convex, with
		$U_j(0)=0, \enspace \dom U_j\subset[-1,1]$.
		By Lemma~\ref{lem:convex-compactness}, after a further subsequence, $U_j\econv U,$
		where $U:\R\to[0,\infty]$ is proper, lower semicontinuous, and convex, with $U(0)=0$.
		
		Set
		$a_*:=-\frac12-x_*, \enspace b_*:=\frac12-x_*$.
		For every $y\in\dom U$, choose recovery points $y_j\to y$ satisfying $U_j(y_j)\to U(y)$. Since $y_j\in\overline I-x_j$ for all sufficiently large $j$, we obtain
		${\dom U\subset[a_*,b_*], \enspace b_*-a_*=1}$.
		
		Let $J_*:=\operatorname{int}\dom U$. If $J_*=\varnothing$, Proposition~\ref{prop:potential-convergence} gives $\lambda_{1,t}(I,V_j)
		=\lambda_1(Q_{t,1,U_j})
		\rightarrow\infty,$
		contradicting \eqref{eq:minimizing-eigen-bound}. Hence $J_*\ne\varnothing$. 
		By Proposition~\ref{prop:potential-convergence},
		$\lambda_k(Q_{t,1,U_j})\to\lambda_k(Q_{t,1,U})$
		for $k=1,2$. Since translation preserves the
		eigenvalues, we obtain
		\begin{equation}\label{eq:attainment-limit-gap}
			\gap_t(J_*,U|_{J_*})
			=\lim_{j\to\infty}\gap_t(I,V_j)
			=\gs_t.
		\end{equation}

		Write $L_*:=|J_*|$. Since $J_*\subset(a_*,b_*)$, we have $0<L_*\le1$. Applying Lemma \ref{lem:convex-gap-bound} to the finite convex potential $U|_{J_*}$ gives
		\[
		\gs_t
		=\gap_t(J_*,U|_{J_*})
		\ge L_*^{-2t}\gs_t.
		\]
		Since $\gs_t>0$ by Proposition~\ref{prop:relative-gap}, we obtain $L_*\ge1$. Consequently,
		\[
		L_*=1,
		\qquad
		J_*=(a_*,b_*)=I-x_*.
		\]
		
		Define
		$V_*(x):=U(x-x_*), \enspace x\in I$.
		Then $V_*$ is finite, nonnegative, and convex on $I$, and \eqref{eq:attainment-limit-gap} yields $\gap_t(I,V_*)=\gs_t.$ 
		
		\smallskip
		
		It remains to verify $\inf_I V_*=0$. If $0\in J_*$, this follows from $U(0)=0$. Otherwise, $0$ is an endpoint of $J_*$. For any $y\in J_*$ and $0<\rho<1$, convexity gives
		\[
		0\le U(\rho y)
		\le\rho U(y)+(1-\rho)U(0)
		=\rho U(y).
		\]
		Letting $\rho\downarrow0$ yields $\inf_{J_*}U=0$, and hence $\inf_I V_*=0$. 
	\end{proof}

	\section{Compact convergence on convex strips}\label{sec:strip}
	
	In this section, we prove compact Mosco convergence of the shifted quadratic forms on convex strips. We show that every sequence with uniformly bounded norms and shifted energies has a strongly convergent subsequence whose limit is of the form $\psi(x)q_{s,H}(z)$. We identify the effective one-dimensional quadratic form in each of the three ranges of $s$ and construct recovery sequences. The resulting convergence yields asymptotic expansions for every fixed eigenvalue, which we then apply to thin rectangles.

	\subsection{Convex strips and the limiting form}

	Let $\eta_j\downarrow0$ and consider the convex strips
	\[K_j
	=
	\bigl\{(x,z)\in\R^2:
	x\in I_j,\ \ell_j(x)<z<u_j(x)\bigr\},\]
	where $I_j\subset\R$ is an open interval, $\ell_j$ is convex, $u_j$ is concave, and
	\[
	0\le\ell_j(x)<u_j(x)\le H_j,
	\quad x\in I_j,
	\qquad
	H_j\longrightarrow H\in(0,\infty).
	\]
	Set
	\[
	h_j:=u_j-\ell_j,\qquad
	\rho_j:=\rho_s(\eta_j),\qquad
	c_j:=\Lambda_sH_j^{-2s},
	\]
	where $\rho_s$ is defined in \eqref{eq:rho}.
	Since $h_j$ is positive and concave, the function
	\begin{equation}\label{eq:strip-potentials}
		W_j(x):=
		\frac{\Lambda_sh_j(x)^{-2s}-c_j}{\rho_j},
		\qquad x\in I_j,
	\end{equation}
	is nonnegative and convex. We extend $W_j$ to $\R$ by its lower semicontinuous endpoint values and set $W_j=+\infty$ outside $\overline I_j$.
	
	Assume that $W_j$ epi-converges to a proper nonnegative lower semicontinuous convex function $W:\R\to[0,\infty]$. Writing $J:=\operatorname{int}\operatorname{dom}W$, Lemma~\ref{lem:convex-compactness} gives
	\[
	W_j\longrightarrow W
	\qquad\text{locally uniformly on }J,
	\]
	and $\inf_AW_j\rightarrow+\infty$
	for every nonempty compact set $A\subset\R\setminus\operatorname{dom}W$. 
	We further assume either that all $I_j$ lie in
	a fixed bounded interval or that, for every
	$M\ge0$, there exists $R_M<\infty$ such that
	\begin{equation}\label{eq:strip-confinement}
		\{W_j\le M\}\subset[-R_M,R_M],
		\qquad j\ge1.
	\end{equation}

	For $v\in\Hs{s}{K_j}$, define the anisotropic Dirichlet quadratic form
	\[
	\mathfrak a_j[v]
	:=
	\iint_{\R^2}
	(\zeta^2+\eta_j^2\xi^2)^s
	|\widehat v(\xi,\zeta)|^2
	\,\dd\xi\,\dd\zeta.
	\]
	The shifted form is
	\begin{equation}\label{eq:strip-forms}
		F_j[v]
		:=
		\frac{\mathfrak a_j[v]-c_j\|v\|_2^2}{\rho_j},
		\qquad v\in\Hs{s}{K_j},
	\end{equation}
	with $F_j[v]=+\infty$ for
	$v\in L^2(\R^2)\setminus\Hs{s}{K_j}$. The effective one-dimensional form is
	\begin{equation}\label{eq:effective-form}
		F[\psi]
		:=
		a_s(H)\int_\R
		|\xi|^{2\sigma_s}|\widehat\psi(\xi)|^2
		\,\dd\xi
		+
		\int_JW(x)|\psi(x)|^2\,\dd x,
	\end{equation}
	where $a_s(H)$ and $\sigma_s$ are defined in
	\eqref{eq:effective-coefficient} and
	\eqref{eq:rho}, respectively. Its domain is
	\[ \mathcal D(F)
	=
	\begin{cases}
		\Hs{s+1/2}{J}\cap L^2(J,W\,\dd x),
		&0<s<\frac12,\\[1mm]
		H_0^1(J)\cap L^2(J,W\,\dd x),
		&\frac12\le s<1.
	\end{cases}\]
	Here all functions are extended by zero to $\R$, and
	\[
	L^2(J,W\,\dd x)
	:=
	\left\{\psi\in L^2(J):
	\int_JW(x)|\psi(x)|^2\,\dd x<\infty\right\}.
	\]
	If $J=\varnothing$, we set $\mathcal D(F)=\{0\}$.
	
	\smallskip
	
	Let $q_{s,H}$ be the positive $L^2$-normalized first Dirichlet eigenfunction on $(0,H)$, extended by zero to $\R$. We define the limiting form on $L^2(\R^2)$ by
	\[F_\infty[v]
	:=
	\begin{cases}
		F[\psi],
		&v(x,z)=\psi(x)q_{s,H}(z),
		\quad\psi\in\mathcal D(F),\\[1mm]
		+\infty,
		&\text{otherwise}.
	\end{cases}\]
	Thus $\mathcal D(F_\infty)$ is contained in the proper closed subspace
	\[
	\bigl\{\psi(x)q_{s,H}(z):
	\psi\in L^2(\R)\bigr\}
	\subset L^2(\R^2),
	\]
	and hence $F_\infty$ is not densely defined on $L^2(\R^2)$.

	We now establish compact Mosco convergence of the shifted forms $F_j$ to the limiting form $F_\infty$ defined above.

	\begin{theorem}\label{thm:strip-limit}
		Under the preceding assumptions, the following statements hold.
		\begin{enumerate}
			\item \textbf{Compactness.}
			Every sequence $(v_j)\subset L^2(\R^2)$ satisfying
			\[
			\sup_j\bigl(\norm{v_j}_2^2+F_j[v_j]\bigr)<\infty
			\]
			admits a subsequence converging strongly in $L^2(\R^2)$ to $\psi(x)q_{s,H}(z)$ for some $\psi\in\mathcal D(F)$.
			
			\item \textbf{Lower bound.}
			For every sequence $v_j\rightharpoonup v$ weakly in $L^2(\R^2)$,
			\[
			F_\infty[v]\le\liminf_{j\to\infty}F_j[v_j].
			\]
			
			\item \textbf{Recovery.}
			For every $\psi\in\mathcal D(F)$, there exist $v_j\in\Hs{s}{K_j}$ such that
			\[
			v_j\longrightarrow\psi q_{s,H}
			\quad\text{strongly in }L^2(\R^2),
			\qquad
			F_j[v_j]\longrightarrow F[\psi].
			\]
		\end{enumerate}
	\end{theorem}

	\subsection{Compactness and the lower bound}
	
	We first prove the compactness and lower-bound assertions of Theorem~\ref{thm:strip-limit}. The recovery sequences will be constructed in the next subsection.
	
	\begin{proof}[\textbf{Proof of
			Theorem~\ref{thm:strip-limit} (i)--(ii)}.]
		\step{Step 1. Transverse estimates.}
		For $v\in\Hs{s}{K_j}$, define
		\[
		\begin{aligned}
			\mathfrak a_0[v]
			&:=\iint_{\R^2}|\zeta|^{2s}
			|\wh v(\xi,\zeta)|^2\,\dd\xi\dd\zeta,\quad
			R_j[v]:=\iint_{\R^2}D_{s,\eta_j|\xi|}(\zeta)
			|\wh v(\xi,\zeta)|^2\,\dd\xi\dd\zeta.
		\end{aligned}
		\]
		Then $\mathfrak a_j=\mathfrak a_0+R_j$, with $R_j\ge0$. Set
		$g_v(x):=\|v(x,\cdot)\|_{L_z^2}, \enspace S_j(x):=(\ell_j(x),u_j(x))$.
		By partial Plancherel formula,
		\[
		\mathfrak a_0[v]
		=\int_{I_j}\int_\R|\zeta|^{2s}
		|\mathcal F_zv(x,\zeta)|^2\,\dd\zeta\dd x.
		\]
		For almost every $x\in I_j$, the function
		$v(x,\cdot)$ belongs to $H^s(\R)$ and vanishes
		almost everywhere outside $S_j(x)$. Since
		$|S_j(x)|=h_j(x)$, the interval Poincar\'e
		inequality gives
		\[
		\begin{aligned}
			&\int_\R|\zeta|^{2s}
			|\mathcal F_zv(x,\zeta)|^2\,\dd\zeta\ge
			\Lambda_sh_j(x)^{-2s}
			\int_\R|v(x,z)|^2\,\dd z.
		\end{aligned}
		\]
		Integrating in $x$ and using
		\eqref{eq:strip-potentials}, we obtain
		\begin{equation}\label{eq:transverse-potential}
			\begin{aligned}
				\mathfrak a_0[v]-c_j\|v\|_2^2
				&\ge\int_{I_j}
				(\Lambda_sh_j^{-2s}-c_j)g_v^2\,\dd x=\rho_j\int_{I_j}W_jg_v^2\,\dd x
				\ge0.
			\end{aligned}
		\end{equation}
		
		We also use the spectral gap on the fixed interval
		$(0,H_j)$. Let $P_j$ be the orthogonal projection
		in the $z$ variable onto
		$\operatorname{span}\{q_{s,H_j}\}$:
		\[
		(P_jv)(x,z)
		=q_{s,H_j}(z)
		\int_\R v(x,y)q_{s,H_j}(y)\,\dd y.
		\]
		Since every transverse slice is supported in
		$(0,H_j)$, the min--max principle gives
		\[
		\begin{aligned}
			&\int_\R|\zeta|^{2s}
			|\mathcal F_zv(x,\zeta)|^2\,\dd\zeta
			-c_j\|v(x,\cdot)\|_{L_z^2}^2\ge
			(\lambda_{2,s}(I,0)-\Lambda_s)H_j^{-2s}
			\|(1-P_j)v(x,\cdot)\|_{L_z^2}^2.
		\end{aligned}
		\]
		Consequently,
		\begin{equation}\label{eq:transverse-projection}
			\mathfrak a_0[v]-c_j\|v\|_2^2
			\ge
			(\lambda_{2,s}(I,0)-\Lambda_s)H_j^{-2s}
			\|(1-P_j)v\|_2^2.
		\end{equation}
		
		Suppose that $\sup_j\bigl(\|v_j\|_2^2+F_j[v_j]\bigr)<\infty.$
		Since
		$\rho_jF_j[v_j] =\mathfrak a_0[v_j]-c_j\|v_j\|_2^2+R_j[v_j]$,
		and both terms on the right are nonnegative,
		\eqref{eq:transverse-potential} and
		\eqref{eq:transverse-projection} imply
		\begin{equation}\label{eq:strip-energy-bounds}
			\int_{I_j}W_jg_{v_j}^2\,\dd x\le C,\quad
			\|(1-P_j)v_j\|_2^2\le C\rho_j,\quad
			R_j[v_j]\le C\rho_j,\quad
			\mathfrak a_0[v_j]\le C.
		\end{equation}
		Here and below, $C$ denotes a constant independent
		of $j$, which may change from line to line.
		
		\step{Step 2. Strong $L^2$ compactness.}
		The last estimate in \eqref{eq:strip-energy-bounds}
		gives, for every $M>0$,
		$\iint_{|\zeta|>M} |\wh v_j(\xi,\zeta)|^2\,\dd\xi\dd\zeta \le CM^{-2s}$.
		
		To control the horizontal frequencies, take the
		Fourier transform only in $x$ and set $\mathcal V_j(\xi,z):=\mathcal F_xv_j(\cdot,z).$
		For almost every $\xi$, the function
		$\mathcal V_j(\xi,\cdot)$ is supported in
		$(0,H_j)$. The variational definition and scaling
		of $B_s$ therefore give (see Lemma \ref{lem:band})
		\begin{equation}\label{eq:horizontal-band-bound}
			\begin{aligned}
				&\mathfrak a_j[v_j]-c_j\|v_j\|_2^2\ge
				\int_\R H_j^{-2s}
				\bigl(B_s(H_j\eta_j|\xi|)-\Lambda_s\bigr)
				\|\mathcal V_j(\xi,\cdot)\|_{L_z^2}^2\,\dd\xi.
			\end{aligned}
		\end{equation}
		Fix $R>0$. By the monotonicity of $B_s$,
		\[
		\begin{aligned}
			&H_j^{-2s}
			\bigl(B_s(H_j\eta_jR)-\Lambda_s\bigr)
			\int_{|\xi|>R}
			\|\mathcal V_j(\xi,\cdot)\|_{L_z^2}^2\,\dd\xi\le
			\mathfrak a_j[v_j]-c_j\|v_j\|_2^2
			\le C\rho_j.
		\end{aligned}
		\]
		By the transverse-band asymptotics in \eqref{eq:band-scaled},
		\[\frac{H_j^{-2s}
			\bigl(B_s(H_j\eta_jR)-\Lambda_s\bigr)}
		{\rho_j}
		\longrightarrow a_s(H)R^{2\sigma_s}>0.\]
		Dividing by the positive band increment and using
		partial Plancherel formula, we obtain
		\[\limsup_{j\to\infty}
		\iint_{|\xi|>R}|\wh v_j|^2\,\dd\xi\dd\zeta
		\le\frac{C}{a_s(H)R^{2\sigma_s}}.\]
		
		Combining the transverse and horizontal estimates
		gives
		\begin{equation}\label{eq:two-frequency-tail}
			\lim_{M\to\infty}\limsup_{j\to\infty}
			\iint_{\xi^2+\zeta^2>M^2}
			|\wh v_j(\xi,\zeta)|^2\,\dd\xi\dd\zeta=0.
		\end{equation}
		Indeed,
		$\{\xi^2+\zeta^2>M^2\} \subset \{|\xi|>M/\sqrt2\} \cup\{|\zeta|>M/\sqrt2\}$.
		
		For $y\in\R^2$, Plancherel's identity gives
		\[
		\|v_j(\cdot+y)-v_j\|_2^2
		=
		\iint_{\R^2}
		|e^{iy\cdot(\xi,\zeta)}-1|^2
		|\wh v_j(\xi,\zeta)|^2\,\dd\xi\dd\zeta.
		\]
		On $\xi^2+\zeta^2\le M^2$, the multiplier is
		bounded by $M^2|y|^2$, and on its complement
		it is bounded by $4$. Hence
		\[
		\begin{aligned}
			\limsup_{j\to\infty}
			\|v_j(\cdot+y)-v_j\|_2^2
			&\le CM^2|y|^2+4\limsup_{j\to\infty}
			\iint_{\xi^2+\zeta^2>M^2}|\wh v_j|^2.
		\end{aligned}
		\]
		First choosing $M$ sufficiently large by
		\eqref{eq:two-frequency-tail} and then taking
		$|y|$ sufficiently small gives asymptotically
		uniform translation continuity. The finitely many
		remaining indices are handled by continuity of
		translations in $L^2(\R^2)$.
		
		The transverse supports lie in a fixed bounded
		interval. If all $I_j$ lie in a fixed bounded
		interval, spatial tightness follows immediately.
		Otherwise, \eqref{eq:strip-confinement} and
		\eqref{eq:strip-energy-bounds} give
		\[
		\begin{aligned}
			&\int_{|x|>R_M}\int_\R
			|v_j(x,z)|^2\,\dd z\dd x\le
			\frac1M\int_{I_j}W_j(x)g_{v_j}(x)^2\,\dd x
			\le\frac CM.
		\end{aligned}
		\]
		The Fr\'echet--Kolmogorov criterion therefore
		gives, after passing to a subsequence,
		\begin{equation}\label{eq:strip-strong-limit}
			v_j\longrightarrow v
			\qquad\text{strongly in }L^2(\R^2).
		\end{equation}
		
		\step{Step 3. Identification of the limit and potential energy.}
		Since $H_j\to H$, continuity of translations
		and dilations in $L^2(\R)$ gives
		$q_{s,H_j}\to q_{s,H}$ strongly in $L^2(\R)$.
		Define
		\[
		\begin{aligned}
			\psi_j(x)&:=\int_\R
			v_j(x,z)q_{s,H_j}(z)\,\dd z,\quad
			\psi(x):=\int_\R
			v(x,z)q_{s,H}(z)\,\dd z.
		\end{aligned}
		\]
		By Cauchy--Schwarz inequality,
		\[
		\begin{aligned}
			\|\psi_j-\psi\|_{L_x^2}
			&\le\|v_j-v\|_2+\|v\|_2
			\|q_{s,H_j}-q_{s,H}\|_{L_z^2}
			\longrightarrow0.
		\end{aligned}
		\]
		Since $P_jv_j=\psi_jq_{s,H_j}$,
		\eqref{eq:strip-energy-bounds} yields
		\[
		\begin{aligned}
			\|v-\psi q_{s,H}\|_2
			&\le\|v-v_j\|_2
			+\|(1-P_j)v_j\|_2+\|\psi_jq_{s,H_j}-\psi q_{s,H}\|_2
			\longrightarrow0.
		\end{aligned}
		\]
		Thus $v=\psi q_{s,H}$. Moreover, $\|g_{v_j}-|\psi|\|_{L_x^2}
		\le\|v_j-v\|_2\rightarrow0.$
		By Cauchy--Schwarz inequality,
		\[
		\|g_{v_j}^2-|\psi|^2\|_{L_x^1}
		\le
		\|g_{v_j}-|\psi|\|_{L_x^2}
		\bigl(\|g_{v_j}\|_{L_x^2}+\|\psi\|_{L_x^2}\bigr)
		\longrightarrow0.
		\]
		Consequently,
		\begin{equation}\label{eq:product-limit-identification}
			v=\psi q_{s,H},\quad
			g_{v_j}\longrightarrow|\psi|
			\quad\text{in }L^2(\R),\quad
			g_{v_j}^2\longrightarrow|\psi|^2
			\quad\text{in }L^1(\R).
		\end{equation}
		
		After a further subsequence,
		$g_{v_j}\to|\psi|$ almost everywhere.
		By epi-convergence, $W_j(x)\rightarrow+\infty$ for every $x\notin\dom W.$
		Set $E:=\R\setminus\dom W$. For every $M>0$,
		Fatou's lemma and \eqref{eq:strip-energy-bounds}
		give
		\[
		\begin{aligned}
			M\int_E|\psi|^2\,\dd x
			&\le\liminf_{j\to\infty}
			\int_E\min\{W_j,M\}g_{v_j}^2\,\dd x\le C.
		\end{aligned}
		\]
		Letting $M\to\infty$ yields $\psi=0$ almost
		everywhere on $E$. Since $\dom W$ differs
		from $J=\operatorname{int}\dom W$ by at most
		two points, we have $\psi=0$ a.e. on $J^c.$
		
		For every compact interval $A\Subset J$,
		the local uniform convergence $W_j\to W$
		and \eqref{eq:product-limit-identification} imply
		\[
		\int_AW_jg_{v_j}^2\,\dd x
		\longrightarrow
		\int_AW|\psi|^2\,\dd x.
		\]
		By exhaustion of $J$ and nonnegativity,
		\begin{equation}\label{eq:potential-lower-limit}
			\liminf_{j\to\infty}
			\int_{I_j}W_jg_{v_j}^2\,\dd x
			\ge\int_JW|\psi|^2\,\dd x.
		\end{equation}
		
		We also record a Fourier-slice estimate used
		below. Choose $H_*>0$ such that $H_j,H\le H_*$
		for all sufficiently large $j$.
		Since $v_j-v$ vanishes outside $(0,H_*)$
		in the $z$ variable, Cauchy--Schwarz gives
		\[
		\begin{aligned}
			|\mathcal F_z(v_j-v)(x,\zeta)|
			&\le(2\pi)^{-1/2}
			\int_0^{H_*}|v_j(x,z)-v(x,z)|\,\dd z\le\left(\frac{H_*}{2\pi}\right)^{1/2}
			\|v_j(x,\cdot)-v(x,\cdot)\|_{L_z^2}.
		\end{aligned}
		\]
		Taking $L_x^2$ norms and applying Plancherel
		in $x$, we obtain
		\begin{equation}\label{eq:mixed-fourier-convergence}
			\sup_{\zeta\in\R}
			\|\wh v_j(\cdot,\zeta)-\wh v(\cdot,\zeta)\|_{L_\xi^2}
			\le C\|v_j-v\|_2
			\longrightarrow0.
		\end{equation}
		The same estimate gives
		$\sup_j\sup_{\zeta\in\R} \|\wh v_j(\cdot,\zeta)\|_{L_\xi^2}<\infty$.
		Here the Fourier slices are understood through
		their continuous $L_\xi^2$-valued representatives.
		
		\step{Step 4. Kinetic lower bound for $1/2<s<1$.}
		In this range, $\rho_j=\eta_j^2$. Fix $R>1$
		and $0<\delta<1$, and set
		$A_{R,\delta} :=\{(\xi,\zeta):|\xi|\le R,\ \delta\le|\zeta|\le R\}$.
		On $A_{R,\delta}$, the following convergence is uniform as $j\to\infty$:
		\[
		\frac{D_{s,\eta_j|\xi|}(\zeta)}{\eta_j^2}
		\longrightarrow s\xi^2|\zeta|^{2s-2}.
		\]
	Moreover,
		\eqref{eq:strip-strong-limit} and Plancherel formula
		imply $|\wh v_j|^2\rightarrow|\wh v|^2$ in $L^1(\R^2).$
		Since the multipliers are uniformly bounded on
		$A_{R,\delta}$,
		\[
		\begin{aligned}
			&\lim_{j\to\infty}
			\int_{A_{R,\delta}}
			\frac{D_{s,\eta_j|\xi|}(\zeta)}{\eta_j^2}
			|\wh v_j|^2\,\dd\xi\dd\zeta=
			s\int_{A_{R,\delta}}
			\xi^2|\zeta|^{2s-2}|\wh v|^2
			\,\dd\xi\dd\zeta.
		\end{aligned}
		\]
		Using $D_{s,k}\ge0$, restricting the full
		energy to $A_{R,\delta}$, and then letting
		$R\uparrow\infty$ and $\delta\downarrow0$,
		we obtain by monotone convergence
		\begin{equation}\label{eq:kinetic-supercritical}
			\begin{aligned}
				\liminf_{j\to\infty}\frac{R_j[v_j]}{\rho_j}
				&\ge s\iint_{\R^2}
				\xi^2|\zeta|^{2s-2}
				|\wh\psi(\xi)|^2
				|\wh q_{s,H}(\zeta)|^2\,\dd\xi\dd\zeta=\beta_sH^{2-2s}
				\int_\R\xi^2|\wh\psi(\xi)|^2\,\dd\xi.
			\end{aligned}
		\end{equation}
		The transverse coefficient is finite and
		strictly positive since $2s-2>-1$ and
		$\wh q_{s,H}\in L^\infty(\R)\cap L^2(\R)$.
		As the left-hand side is bounded,
		\eqref{eq:kinetic-supercritical} also gives
		$\psi\in H^1(\R)$.
		
		\step{Step 5. Kinetic lower bound for $s=1/2$.}
		Set $L_j:=\log(1/\eta_j),\,\rho_j=\eta_j^2L_j.$
		Fix $0<\alpha<1$, $R>0$, and $0<\delta<1$.
		Consider
		$A_j:=\{(\xi,\zeta):|\xi|\le R,\ \eta_j^\alpha\le|\zeta|\le\delta\}$.
		Since $\sup_{A_j}\frac{\eta_j|\xi|}{|\zeta|}
		\le R\eta_j^{1-\alpha}\rightarrow0,$
		rationalization gives, uniformly on $A_j$,
		\[
		\begin{aligned}
			D_{1/2,\eta_j|\xi|}(\zeta)
			&=\frac{\eta_j^2\xi^2}
			{\sqrt{\zeta^2+\eta_j^2\xi^2}+|\zeta|}\ge(1-o_j(1))
			\frac{\eta_j^2\xi^2}{2|\zeta|}.
		\end{aligned}
		\]
		Consequently,
		\begin{equation}\label{eq:critical-annulus}
			\begin{aligned}
				\frac{R_j[v_j]}{\rho_j}
				&\ge\frac{1-o_j(1)}{2L_j}
				\int_{\eta_j^\alpha\le|\zeta|\le\delta}
				\frac1{|\zeta|}\times
				\int_{|\xi|\le R}
				\xi^2|\wh v_j(\xi,\zeta)|^2
				\,\dd\xi\dd\zeta.
			\end{aligned}
		\end{equation}
		
		By \eqref{eq:mixed-fourier-convergence},
		uniformly in $\zeta$,
		\[
		\begin{aligned}
			&\int_{|\xi|\le R}\xi^2
			\bigl||\wh v_j(\xi,\zeta)|^2
			-|\wh v(\xi,\zeta)|^2\bigr|\,\dd\xi\le C_R\|v_j-v\|_2.
		\end{aligned}
		\]
		Furthermore,
		\[
		\frac1{2L_j}
		\int_{\eta_j^\alpha\le|\zeta|\le\delta}
		\frac{\dd\zeta}{|\zeta|}
		=\alpha+\frac{\log\delta}{L_j}
		\longrightarrow\alpha.
		\]
		Hence the error in replacing
		$|\wh v_j|^2$ by $|\wh v|^2$ in
		\eqref{eq:critical-annulus}, after normalization,
		is bounded by
		$C_{R,\alpha}\|v_j-v\|_2\longrightarrow0$.
		
		Set
		\[
		G(\zeta):=|\wh q_{1/2,H}(\zeta)|^2,
		\qquad
		\omega(\delta):=
		\sup_{|\zeta|\le\delta}|G(\zeta)-G(0)|.
		\]
		Since $G$ is continuous, $\omega(\delta)\to0$
		as $\delta\downarrow0$. Using
		$\wh v(\xi,\zeta)=\wh\psi(\xi)\wh q_{1/2,H}(\zeta)$,
		we obtain
		\[
		\liminf_{j\to\infty}\frac{R_j[v_j]}{\rho_j}
		\ge
		\alpha\bigl(G(0)-\omega(\delta)\bigr)
		\int_{|\xi|\le R}\xi^2|\wh\psi(\xi)|^2\,\dd\xi.
		\]
		Letting $\delta\downarrow0$, then
		$\alpha\uparrow1$, and finally $R\uparrow\infty$,
		and using $G(0)=\gamma H$, we conclude that
		\begin{equation}\label{eq:kinetic-critical}
			\liminf_{j\to\infty}\frac{R_j[v_j]}{\rho_j}
			\ge\gamma H
			\int_\R\xi^2|\wh\psi(\xi)|^2\,\dd\xi.
		\end{equation}
		In particular, $\psi\in H^1(\R)$.
		
		\step{Step 6. Kinetic lower bound for $0<s<1/2$.}
		Here $\rho_j=\eta_j^{1+2s}$. The substitution
		$\zeta=\eta_jt$ gives the exact identity
		\begin{equation}\label{eq:subcritical-blowup-integral}
			\begin{aligned}
				\frac{R_j[v_j]}{\rho_j}
				&=\iint_{\R^2}
				\bigl((t^2+\xi^2)^s-|t|^{2s}\bigr)\times
				|\wh v_j(\xi,\eta_jt)|^2\,\dd\xi\dd t.
			\end{aligned}
		\end{equation}
		
		Fix $R,M>0$. By
		\eqref{eq:mixed-fourier-convergence} and
		continuity of $\wh q_{s,H}$,
		\[
		\wh v_j(\xi,\eta_jt)
		\longrightarrow
		\wh\psi(\xi)\wh q_{s,H}(0)
		\quad\text{in }L^2([-R,R]\times[-M,M]).
		\]
		Indeed, uniformly for $|t|\le M$,
		\[
		\begin{aligned}
			&\|\wh v_j(\cdot,\eta_jt)
			-\wh\psi\,\wh q_{s,H}(0)\|_{L_\xi^2}\le C\|v_j-v\|_2
			+\|\psi\|_2
			\sup_{|t|\le M}
			|\wh q_{s,H}(\eta_jt)-\wh q_{s,H}(0)|
			\longrightarrow0.
		\end{aligned}
		\]
		
		Set $Q_{R,M}:=[-R,R]\times[-M,M]$ and
		\[
		\Phi(\xi,t):=(t^2+\xi^2)^s-|t|^{2s},\quad
		U_j(\xi,t):=\wh v_j(\xi,\eta_jt),\quad
		U(\xi,t):=\wh\psi(\xi)\wh q_{s,H}(0).
		\]
		Since $\Phi$ is nonnegative and bounded on
		$Q_{R,M}$,
		\[
		\left|\int_{Q_{R,M}}
		\Phi\bigl(|U_j|^2-|U|^2\bigr)
		\,\dd\xi\dd t\right|
		\rightarrow0.
		\]
		Thus the truncated energies converge.
		Since $\Phi\ge0$, we may restrict the energy
		in \eqref{eq:subcritical-blowup-integral}
		to $Q_{R,M}$ before passing to the limit.
		Letting $R,M\uparrow\infty$ and applying
		monotone convergence, we obtain
		\begin{equation}\label{eq:kinetic-subcritical}
			\begin{aligned}
				\liminf_{j\to\infty}\frac{R_j[v_j]}{\rho_j}
				&\ge|\wh q_{s,H}(0)|^2
				\int_\R|\wh\psi(\xi)|^2\times
				\int_\R\bigl((t^2+\xi^2)^s-|t|^{2s}\bigr)
				\,\dd t\dd\xi\\
				&=\kappa_sH
				\int_\R|\xi|^{1+2s}
				|\wh\psi(\xi)|^2\,\dd\xi.
			\end{aligned}
		\end{equation}
		Here, for $\xi\ne0$, the substitution
		$t=|\xi|r$ gives
		\[
		\int_\R\bigl((t^2+\xi^2)^s-|t|^{2s}\bigr)\,\dd t
		=C_s|\xi|^{1+2s},
		\]
		while the integral vanishes at $\xi=0$.
		Moreover,
		$|\wh q_{s,H}(0)|^2=H|\wh q_s(0)|^2$.
		
		\step{Step 7. Completion of the lower bound.}
		By \eqref{eq:strip-energy-bounds}, the left-hand
		sides of \eqref{eq:kinetic-supercritical},
		\eqref{eq:kinetic-critical}, and
		\eqref{eq:kinetic-subcritical} are bounded.
		Their right-hand sides therefore give
		$\psi\in H^{\sigma_s}(\R)$.
		
		Since $\psi=0$ almost everywhere on $J^c$
		and \eqref{eq:potential-lower-limit} gives
		$\int_JW|\psi|^2\,\dd x<\infty$, we obtain
		$\psi\in\mathcal D(F)$. In particular, when
		$\sigma_s=1$, the zero-extension condition
		and $J$ being an interval give
		$\psi\in H_0^1(J)$. By \eqref{eq:transverse-potential} and
		\eqref{eq:strip-forms},
		\[
		F_j[v_j]
		\ge\frac{R_j[v_j]}{\rho_j}
		+\int_{I_j}W_jg_{v_j}^2\,\dd x.
		\]
		Combining \eqref{eq:potential-lower-limit}
		with the corresponding kinetic lower bound,
		we conclude that
		\[
		F[\psi]\le\liminf_{j\to\infty}F_j[v_j].
		\]
		This proves the compactness assertion (i).
		
		\smallskip
		
		To prove (ii), let $v_j\rightharpoonup v$
		weakly in $L^2(\R^2)$, and set $\ell:=\liminf_{j\to\infty}F_j[v_j].$
		If $\ell=+\infty$, the claim is immediate.
		Otherwise, choose a subsequence along which
		$F_j[v_j]\to\ell<\infty$. Its $L^2$ norms are
		bounded, so the compactness argument gives
		a further subsequence converging strongly
		to $\psi q_{s,H}$ for some
		$\psi\in\mathcal D(F)$. By uniqueness of the
		weak limit, $v=\psi q_{s,H}$, and
		$F_\infty[v]=F[\psi] \le\ell=\liminf_{j\to\infty}F_j[v_j]$.
		This proves the lower-bound assertion (ii).
	\end{proof}

	\subsection{Recovery sequences on common transverse intervals}
	
	We complete the proof of Theorem~\ref{thm:strip-limit} by constructing recovery sequences. Since the transverse sections depend on $x$, the ground state $q_{s,H_j}$ need not be supported in each section. We instead use fixed transverse intervals contained in all sections over small horizontal cells. The resulting trial functions are finite sums of products, to which Lemma~\ref{lem:product-limit} applies.
	
	\begin{proof}[\textbf{Proof of
			Theorem~\ref{thm:strip-limit} (iii)}.]
		\step{Step 1. Control of the transverse boundaries.}
		First let $\psi\in C_c^\infty(J)$, where $J=\operatorname{int}\dom W\ne\varnothing$. Choose compact intervals $A_0,A_1$ such that
		\[
		\supp\psi\subset\operatorname{int}A_0,
		\qquad
		A_0\Subset\operatorname{int}A_1,
		\qquad
		A_1\Subset J.
		\]
		By the local uniform convergence $W_j\to W$, we know that $0\le W_j(x)\le C,\,x\in A_1,$
		for all sufficiently large $j$. Since $h_j\le H_j$ and $r\mapsto r^{-2s}$ is convex,
		\[
		\begin{aligned}
			\rho_jW_j(x)
			&=\Lambda_s\bigl(h_j(x)^{-2s}-H_j^{-2s}\bigr)\ge 2s\Lambda_sH_j^{-2s-1}
			\bigl(H_j-h_j(x)\bigr).
		\end{aligned}
		\]
		As $H_j\to H>0$, we obtain
		\begin{equation}\label{eq:height-close}
			0\le H_j-h_j(x)\le C\rho_j,
			\qquad x\in A_1.
		\end{equation}
		
		Define
		\[
		A_j(x):=\frac{\ell_j(x)}{\rho_j},
		\qquad
		B_j(x):=\frac{H_j-u_j(x)}{\rho_j}.
		\]
		Both functions are nonnegative and convex, and
		\[
		0\le A_j(x)+B_j(x)
		=\frac{H_j-h_j(x)}{\rho_j}\le C,
		\qquad x\in A_1.
		\]
		Let $d:=\operatorname{dist}(A_0,\partial A_1)>0$. The secant-slope inequalities for convex functions bounded on $A_1$ give
		\begin{equation}\label{eq:endpoint-Lipschitz}
			|A_j(x)-A_j(y)|+|B_j(x)-B_j(y)|
			\le C_0|x-y|,
			\qquad x,y\in A_0,
		\end{equation}
		where $C_0$ is independent of $j$.
		
		\step{Step 2. Common transverse intervals.}
		Fix $r>0$ sufficiently small. Cover a neighborhood of $\supp\psi$ by finitely many open intervals $U_m\Subset\operatorname{int}A_0$, $1\le m\le M$, with $\operatorname{diam}U_m\le r$. Choose $\chi_m\in C_c^\infty(U_m)$ such that
		\[
		\chi_m\ge0,\qquad
		\sum_{m=1}^M\chi_m=1
		\quad\text{near }\supp\psi.
		\]
		Set $a_{j,m}:=\sup_{U_m}A_j,\,
		b_{j,m}:=\sup_{U_m}B_j,$
		and define
		\[
		S_{j,m}:=(\rho_ja_{j,m},H_j-\rho_jb_{j,m}),
		\qquad
		L_{j,m}:=|S_{j,m}|.
		\]
		By \eqref{eq:height-close}, $a_{j,m}$ and $b_{j,m}$ are uniformly bounded. Hence
		\[
		L_{j,m}
		=H_j-\rho_j(a_{j,m}+b_{j,m})
		\longrightarrow H>0
		\]
		for each fixed $r$ and $m$, so $S_{j,m}$ is nonempty for all sufficiently large $j$.
		
		For every $x\in U_m$,
		$\ell_j(x)\le\rho_ja_{j,m}, \enspace H_j-\rho_jb_{j,m}\le u_j(x)$.
		Moreover, \eqref{eq:endpoint-Lipschitz} gives
		\[
		\begin{aligned}
			h_j(x)-L_{j,m}
			&=\rho_j\bigl(a_{j,m}-A_j(x)
			+b_{j,m}-B_j(x)\bigr)\\
			&\le C_0\rho_j\operatorname{diam}U_m
			\le C_0\rho_jr.
		\end{aligned}
		\]
		Consequently,
		\begin{equation}\label{eq:common-interval}
			S_{j,m}\subset S_j(x),
			\qquad
			0\le h_j(x)-L_{j,m}\le C_0\rho_jr,
			\quad x\in U_m,
		\end{equation}
		where $S_j(x)=(\ell_j(x),u_j(x))$.
		
		\step{Step 3. Trial functions and their $L^2$ limits.}
		Let $q_{j,m}$ be the positive $L^2$-normalized first Dirichlet eigenfunction on $S_{j,m}$, extended by zero to $\R$. Define
		\[r_j^r(x,z):=\sum_{m=1}^M\chi_m(x)q_{j,m}(z),
			\qquad
			v_j^r(x,z):=\psi(x)r_j^r(x,z).\]
		
		We first verify that $v_j^r\in\Hs{s}{K_j}$. Set $f_m:=\psi\chi_m\in C_c^\infty(U_m)$. Since $0<s<1$, $(\xi^2+\zeta^2)^s
		\le|\xi|^{2s}+|\zeta|^{2s}.$
		Thus
		\[
		\begin{aligned}
			&\iint_{\R^2}(\xi^2+\zeta^2)^s
			|\wh f_m(\xi)|^2|\wh q_{j,m}(\zeta)|^2
			\,\dd\xi\dd\zeta\le
			\Es_s[f_m]\|q_{j,m}\|_2^2
			+\|f_m\|_2^2\Es_s[q_{j,m}]
			<\infty.
		\end{aligned}
		\]
		Moreover, $f_mq_{j,m}$ vanishes almost everywhere outside $U_m\times S_{j,m}\subset K_j.$
		Hence $f_mq_{j,m}\in\Hs{s}{K_j}$, and therefore $v_j^r\in\Hs{s}{K_j}$.
		
		\smallskip
		
		Let $c_{j,m}$ denote the midpoint of $S_{j,m}$. By scaling,
		\[
		q_{j,m}(z)
		=L_{j,m}^{-1/2}
		q_s\left(\frac{z-c_{j,m}}{L_{j,m}}\right),
		\qquad
		c_{j,m}\to\frac H2,\quad L_{j,m}\to H.
		\]
		Since $q_s\in L^1(\R)\cap L^2(\R)$, continuity of translations and dilations gives, for each fixed $r,m$,
		\begin{equation}\label{eq:common-ground-convergence}
			q_{j,m}\longrightarrow q_{s,H}
			\qquad\text{in }L^1(\R)\cap L^2(\R).
		\end{equation}
		Using $\sum_mf_m=\psi$ and the finiteness of $M$, we obtain
		\[
		\begin{aligned}
			\|v_j^r-\psi q_{s,H}\|_2
			&=\left\|
			\sum_{m=1}^Mf_m(x)
			\bigl(q_{j,m}(z)-q_{s,H}(z)\bigr)
			\right\|_2\\
			&\le\sum_{m=1}^M
			\|f_m\|_{L_x^2}
			\|q_{j,m}-q_{s,H}\|_{L_z^2}
			\longrightarrow0.
		\end{aligned}
		\]
		
		\step{Step 4. Transverse energy estimate.}
		For $w\in\Hs{s}{(0,H_j)}$, define
		\[
		\mathcal T_j[w]
		:=
		\int_\R|\zeta|^{2s}|\wh w(\zeta)|^2
		\,\dd\zeta-c_j\|w\|_2^2.
		\]
		Since $c_j=\Lambda_sH_j^{-2s}$, the interval Poincar\'e inequality gives $\mathcal T_j[w]\ge0$. For each $q_{j,m}$,
		${\mathcal T_j[q_{j,m}] =\Lambda_sL_{j,m}^{-2s}-c_j}$.
		For $x\in U_m$, the definition of $W_j$ gives
		\[
		\frac{\mathcal T_j[q_{j,m}]}{\rho_j}
		=
		W_j(x)
		+\frac{\Lambda_s}{\rho_j}
		\bigl(L_{j,m}^{-2s}-h_j(x)^{-2s}\bigr).
		\]
		Since $L_{j,m}\le h_j(x)$ and both lengths converge uniformly to $H>0$ on the relevant cells, the mean value theorem and \eqref{eq:common-interval} yield
		\[
		\begin{aligned}
			0\le L_{j,m}^{-2s}-h_j(x)^{-2s}
			&\le C\bigl(h_j(x)-L_{j,m}\bigr)\le C\rho_jr.
		\end{aligned}
		\]
		Consequently,
		\begin{equation}\label{eq:cell-transverse}
			\frac{\mathcal T_j[q_{j,m}]}{\rho_j}
			\le W_j(x)+Cr,
			\qquad x\in U_m,
		\end{equation}
		where $C$ is independent of $j,m$, and sufficiently small $r$.
		
		\smallskip
		
		For any $w_1,\ldots,w_M\in\Hs{s}{(0,H_j)}$ and $\theta_m\ge0$ with $\sum_m\theta_m=1$, polarization gives
		\[
		\begin{aligned}
			&\sum_{m=1}^M\theta_m\mathcal T_j[w_m]
			-\mathcal T_j\left[\sum_{m=1}^M\theta_mw_m\right]=
			\frac12\sum_{m,n=1}^M
			\theta_m\theta_n\mathcal T_j[w_m-w_n]
			\ge0.
		\end{aligned}
		\]
		For $x\in\supp\psi$, take $\theta_m=\chi_m(x)$ and $w_m=q_{j,m}$. By \eqref{eq:cell-transverse},
		\[
		\begin{aligned}
			\mathcal T_j[r_j^r(x,\cdot)]
			&\le\sum_{m=1}^M\chi_m(x)\mathcal T_j[q_{j,m}]\le\rho_j\bigl(W_j(x)+Cr\bigr).
		\end{aligned}
		\]
		Using partial Plancherel formula,
		\[
		\mathfrak a_0[v_j^r]-c_j\|v_j^r\|_2^2
		=
		\int_\R|\psi(x)|^2
		\mathcal T_j[r_j^r(x,\cdot)]\,\dd x.
		\]
		Since $W_j\to W$ uniformly on $\supp\psi$, we obtain
		\begin{equation}\label{eq:recovery-transverse}
			\limsup_{j\to\infty}
			\frac{\mathfrak a_0[v_j^r]-c_j\|v_j^r\|_2^2}{\rho_j}
			\le
			\int_JW(x)|\psi(x)|^2\,\dd x
			+Cr\|\psi\|_2^2.
		\end{equation}
		
		\step{Step 5. Longitudinal energy and mixed terms.}
		For fixed $r$, write
		\[
		v_j^r(x,z)=\sum_{m=1}^Mf_m(x)q_{j,m}(z),
		\qquad f_m=\psi\chi_m.
		\]
		Expanding the Fourier energy gives
		\[
		\begin{aligned}
			\frac{R_j[v_j^r]}{\rho_j}
			&=\sum_{m,n=1}^M\frac1{\rho_j}
			\iint_{\R^2}D_{s,\eta_j|\xi|}(\zeta)
			\wh f_m(\xi)\overline{\wh f_n(\xi)}\times
			\wh q_{j,m}(\zeta)\overline{\wh q_{j,n}(\zeta)}
			\,\dd\xi\dd\zeta.
		\end{aligned}
		\]
		By \eqref{eq:common-ground-convergence}, Lemma~\ref{lem:product-limit} applies to every pair $(m,n)$. Since $M$ is fixed,
		\[
		\begin{aligned}
			\lim_{j\to\infty}\frac{R_j[v_j^r]}{\rho_j}=a_s(H)\int_\R|\xi|^{2\sigma_s}
			\left|\sum_{m=1}^M\wh f_m(\xi)\right|^2
			\,\dd\xi.
		\end{aligned}
		\]
		Since $\sum_mf_m=\psi$, we conclude that
		\begin{equation}\label{eq:recovery-longitudinal}
			\lim_{j\to\infty}\frac{R_j[v_j^r]}{\rho_j}
			=
			a_s(H)\int_\R|\xi|^{2\sigma_s}
			|\wh\psi(\xi)|^2\,\dd\xi.
		\end{equation}
		Combining \eqref{eq:recovery-transverse} and \eqref{eq:recovery-longitudinal}, we obtain, for every fixed sufficiently small $r>0$,
		\[\limsup_{j\to\infty}F_j[v_j^r]
			\le F[\psi]+Cr\|\psi\|_2^2,\]
		where $C$ is independent of $r$.
		
		\step{Step 6. Diagonal choice of the mesh size.}
		Choose $r_n\downarrow0$. For each fixed $n$, Steps 3--5 give $v_j^{r_n}\rightarrow\psi q_{s,H}$ in $L^2(\R^2),$
		and
		\[
		\limsup_{j\to\infty}F_j[v_j^{r_n}]
		\le F[\psi]+Cr_n\|\psi\|_2^2.
		\]
		Choose strictly increasing integers $j_n$ such that, for every $j\ge j_n$,
		\[
		\begin{aligned}
			\|v_j^{r_n}-\psi q_{s,H}\|_2&<\frac1n,\quad
			F_j[v_j^{r_n}]
			\le F[\psi]+Cr_n\|\psi\|_2^2+\frac1n.
		\end{aligned}
		\]
		Define $v_j:=v_j^{r_n}$ for $j_n\le j<j_{n+1}$, and set $v_j=0$ for $j<j_1$. Then
		\[
		v_j\longrightarrow\psi q_{s,H}
		\quad\text{in }L^2(\R^2),
		\qquad
		\limsup_{j\to\infty}F_j[v_j]\le F[\psi].
		\]
		The Mosco lower bound established in the preceding subsection gives
		\[
		F[\psi]\le\liminf_{j\to\infty}F_j[v_j]
		\le\limsup_{j\to\infty}F_j[v_j]\le F[\psi].
		\]
		Thus, $	F_j[v_j]\longrightarrow F[\psi].$
		This proves the recovery assertion for every $\psi\in C_c^\infty(J)$.
		
		\step{Step 7. Extension to the full form domain.}
		Let $\psi\in\mathcal D(F)$. By Lemma~\ref{lem:form-core}, there exist $\psi_n\in C_c^\infty(J)$ such that
		$\|\psi_n-\psi\|_2\longrightarrow0, \enspace F[\psi_n-\psi]\longrightarrow0$.
		Since $F$ is a nonnegative quadratic form, the Cauchy--Schwarz inequality gives
		\[
		\begin{aligned}
			|F[\psi_n]-F[\psi]|
			&\le F[\psi_n-\psi]
			+2F[\psi]^{1/2}F[\psi_n-\psi]^{1/2}\rightarrow0.
		\end{aligned}
		\]
		
		For each $n$, as $j\to\infty$, Step 6 provides $v_{j,n}\in\Hs{s}{K_j}$ satisfying
		\[
		v_{j,n}\longrightarrow\psi_nq_{s,H}
		\quad\text{in }L^2(\R^2),
		\qquad
		F_j[v_{j,n}]\longrightarrow F[\psi_n].
		\]
		Choose strictly increasing integers $J_n$ such that, for every $j\ge J_n$,
		\[
		\begin{aligned}
			\|v_{j,n}-\psi_nq_{s,H}\|_2&<\frac1n,\quad
			|F_j[v_{j,n}]-F[\psi_n]|<\frac1n.
		\end{aligned}
		\]
		Define $v_j:=v_{j,n}$ for $J_n\le j<J_{n+1}$, and set $v_j=0$ for $j<J_1$. Then
		\[
		\begin{aligned}
			\|v_j-\psi q_{s,H}\|_2
			&\le\frac1n+\|\psi_n-\psi\|_2
			\longrightarrow0,
		\end{aligned}
		\]
		and
		\[
		\begin{aligned}
			|F_j[v_j]-F[\psi]|
			&\le\frac1n+|F[\psi_n]-F[\psi]|
			\longrightarrow0.
		\end{aligned}
		\]
		If $J=\varnothing$, then $\mathcal D(F)=\{0\}$, and the zero sequence is a recovery sequence.
		
		This proves part (iii). Together with parts (i)--(ii), it completes the proof of Theorem~\ref{thm:strip-limit}.
	\end{proof}

\subsection{Eigenvalue convergence and thin rectangles}

The compact Mosco convergence established in Theorem~\ref{thm:strip-limit} yields convergence of the eigenvalues of the shifted forms.

\begin{corollary}\label{cor:strip-spectrum}
Suppose that $J\ne\varnothing$ in Theorem~\ref{thm:strip-limit}. Let $\nu_{k,j}$ be the eigenvalues of $\mathfrak a_j$ and $\mu_k$ those of $F$, counted with multiplicity. The associated operators have compact resolvent on the $L^2$ closures of their form domains. For every fixed $k\ge1$,
\begin{equation}\label{eq:strip-spectrum}
\frac{\nu_{k,j}-c_j}{\rho_s(\eta_j)}
\longrightarrow\mu_k,
\qquad
\frac{\nu_{2,j}-\nu_{1,j}}{\rho_s(\eta_j)}
\longrightarrow\mu_2-\mu_1.
\end{equation}
\end{corollary}

\begin{proof}
We apply Lemma~\ref{lem:abstract-spectral} with
$\mathcal H=L^2(\R^2)$, $Q_j=F_j$, and $Q=F_\infty$.

We first verify the compact embedding of $\mathcal D(F_j)$ into $\mathcal H$. For each fixed $j$, since $\eta_j>0$,
\[
(\zeta^2+\eta_j^2\xi^2)^s
\ge d_j(\xi^2+\zeta^2)^s,
\qquad d_j:=\min\{1,\eta_j^{2s}\}>0.
\]
Hence $\Es_s[v]\le d_j^{-1}\mathfrak a_j[v]$ for every $v\in\Hs{s}{K_j}$. Using $\mathfrak a_j[v]=\rho_jF_j[v]+c_j\|v\|_2^2$, we obtain
\[
\|v\|_{H^s(\R^2)}^2
\le C_j\bigl(\|v\|_2^2+F_j[v]\bigr),
\quad v\in\Hs{s}{K_j}.
\]
Thus boundedness in the form norm of $F_j$ implies boundedness in $H^s(\R^2)$.

Since $K_j\subset\R\times(0,H_j)$, the transverse supports are bounded. If $I_j$ is bounded, spatial tightness follows directly. Otherwise, the confinement assumption and \eqref{eq:transverse-potential} give, for every $M>0$,
\[
\begin{aligned}
\int_{|x|>R_M}\int_\R|v(x,z)|^2\,\dd z\dd x
&\le\frac1M\int_{I_j}W_j(x)g_v(x)^2\,\dd x\le\frac{F_j[v]}{M}.
\end{aligned}
\]
The $H^s$ bound gives uniform translation continuity, while the preceding estimate gives spatial tightness. 
For every sufficiently large $j$, the Fr\'echet--Kolmogorov criterion therefore yields the compact embedding
$\mathcal D(F_j)=\Hs{s}{K_j} \hookrightarrow\hookrightarrow L^2(\R^2)$.

We next identify the limiting form. The kinetic form is closed on its zero-extension Sobolev domain, and the nonnegative potential form is closed on $L^2(J,(1+W)\,\dd x)$. Hence $F$ is a nonnegative closed form on $L^2(J)$. Since $W$ is finite and locally bounded on $J$, we have $C_c^\infty(J)\subset\mathcal D(F)$, so $\mathcal D(F)$ is dense in $L^2(J)$. Moreover, $\|q_{s,H}\|_2=1$, and thus
$T:L^2(J)\longrightarrow L^2(\R^2), \enspace T\psi:=\psi q_{s,H}$,
is an isometry. It follows that
\[
\overline{\mathcal D(F_\infty)}^{\,L^2(\R^2)}
=T(L^2(J)),
\qquad
F_\infty[T\psi]=F[\psi],
\quad \psi\in\mathcal D(F).
\]
Thus $F_\infty$, restricted to the closure of its form domain, is unitarily equivalent to $F$ on $L^2(J)$.

Finally, the compactness, lower-bound, and recovery assertions of Theorem~\ref{thm:strip-limit}, together with the compact embedding proved above, verify the hypotheses of Lemma~\ref{lem:abstract-spectral}. Consequently, the limiting operator has compact resolvent, and for every fixed $k\ge1$,
$\lambda_k(F_j)\longrightarrow\lambda_k(F)=\mu_k$.
Since
\[
\lambda_k(F_j)=\frac{\nu_{k,j}-c_j}{\rho_j},
\qquad \rho_j=\rho_s(\eta_j),
\]
the first convergence in \eqref{eq:strip-spectrum} follows. Subtracting the limits for $k=1,2$ gives the second.
\end{proof}

We apply this result to thin rectangles, for which the effective potential vanishes on the longitudinal interval.

\begin{proof}[\textbf{Proof of
			Theorem~\ref{thm:rectangles}}.]
Fix $L>0$ and translate $R_{L,h}$ to $I_L\times(0,h)$, where $I_L=(-L/2,L/2)$ and $I=(-1/2,1/2)$. For $u\in\Hs{s}{I_L\times(0,h)}$, define
$U(X,z):=(Lh)^{1/2}u(LX,hz), \enspace \eta:=\frac hL$.
Then $U\in\Hs{s}{I\times(0,1)}$ and $\|U\|_2=\|u\|_2$. The Fourier scaling formula gives
${\wh u(\xi,\zeta) =(Lh)^{1/2}\wh U(L\xi,h\zeta)}$.
Consequently,
\[
\begin{aligned}
\Es_s[u]
&=\iint_{\R^2}
(\xi^2+\zeta^2)^s|\wh u(\xi,\zeta)|^2
\,\dd\xi\dd\zeta=h^{-2s}\iint_{\R^2}
(\zeta^2+\eta^2\xi^2)^s
|\wh U(\xi,\zeta)|^2\,\dd\xi\dd\zeta.
\end{aligned}
\]
Let $\nu_n(\eta)$ denote the $n$th eigenvalue of the anisotropic form on $I\times(0,1)$. By the min--max principle,
$\lambda_{n,s}(R_{L,h}) =h^{-2s}\nu_n(h/L)$.

For the fixed strip $K=I\times(0,1)$, we have $H=1$, $c=\Lambda_s$, and
\[
W(x)=
\begin{cases}
0,&x\in\overline I,\\
+\infty,&x\notin\overline I.
\end{cases}
\]
Thus $J=I$. Let $\mu_{n,s}$ denote the $n$th eigenvalue of the effective form $F$ on $L^2(I)$. Since $W=0$ on $I$, the effective form is a multiple of the restricted fractional Dirichlet form of order $s+1/2$ when $0<s<1/2$, and of the classical Dirichlet form when $1/2\le s<1$. Using the coefficients $a_s(1)$ in \eqref{eq:effective-coefficient} and the Dirichlet eigenvalues $(n\pi)^2$ on $I$, we obtain
\[
\mu_{n,s}=
\begin{cases}
\kappa_s\lambda_{n,s+1/2}(I,0),
&0<s<\frac12,\\[1mm]
\gamma(n\pi)^2,
&s=\frac12,\\[1mm]
\beta_s(n\pi)^2,
&\frac12<s<1.
\end{cases}
\]
Corollary~\ref{cor:strip-spectrum} therefore gives, for every fixed $n\ge1$,
\[
\nu_n(\eta)
=\Lambda_s+\rho_s(\eta)\mu_{n,s}
+o(\rho_s(\eta)),
\qquad \eta\downarrow0.
\]
Therefore,
\[
\lambda_{n,s}(R_{L,h})
=\Lambda_sh^{-2s}
+h^{-2s}\rho_s(h/L)\mu_{n,s}
+o\bigl(h^{-2s}\rho_s(h/L)\bigr).
\]
The scaling factors satisfy
\[
h^{-2s}\rho_s(h/L)
=
\begin{cases}
hL^{-1-2s},
&0<s<\frac12,\\[1mm]
hL^{-2}\log(L/h),
&s=\frac12,\\[1mm]
h^{2-2s}L^{-2},
&\frac12<s<1.
\end{cases}
\]
For $0<s<\frac12$, the one-dimensional scaling identity gives
\[
\lambda_{n,s+1/2}(I_L,0)
=L^{-1-2s}\lambda_{n,s+1/2}(I,0).
\]
Substitution into the preceding expansion yields \eqref{eq:rectangle-expansion}, completing the proof.
\end{proof}

\section{A gap bound implies a first-eigenvalue bound}\label{sec:window}

In this section, we prove that if the fundamental gap of a thin convex domain is bounded by $C\rho_s(\varepsilon_j)$, then the first eigenvalue differs from the transverse eigenvalue $c_j=\Lambda_sH_j^{-2s}$ by at most $C'\rho_s(\varepsilon_j)$. To prove this, we suppose that $(\nu_{1,j}-c_j)/\rho_s(\varepsilon_j)\to\infty$ and introduce a new horizontal scale determined by $\rho_s(\eta_j)=\nu_{1,j}-c_j$. After rescaling, the first two eigenvalues converge to $1$. We then show that the rescaled transverse potentials have uniformly bounded sublevel sets, allowing us to pass to a one-dimensional limit whose first eigenvalue is simple. This contradicts the convergence of the first two eigenvalues to the same value.

\smallskip

Let $\Omega_j\subset\R^2$ be bounded open convex sets satisfying $\diam\Omega_j=1,\,
|\Omega_j|=\varepsilon_j\downarrow0.$
By Lemma~\ref{lem:width-chord}, after a rigid motion, the vertical direction realizes the minimal width $w_j$, the corresponding supporting lines are $y=0$ and $y=w_j$, and $\{0\}\times[0,w_j]\subset\overline\Omega_j.$
Define
\begin{equation}\label{eq:geometric-normalization}
K_j:=\{(x,z):(x,\varepsilon_jz)\in\Omega_j\},
\qquad H_j:=\frac{w_j}{\varepsilon_j}.
\end{equation}
Since $\varepsilon_j\le w_j\le2\varepsilon_j$ by Lemma~\ref{lem:width-chord}, the anisotropic scaling gives
\begin{equation}\label{eq:normalized-geometry}
\begin{gathered}
|K_j|=1,\quad 1\le H_j\le2,\quad
K_j\subset\R\times(0,H_j),\quad
\{0\}\times[0,H_j]\subset\overline K_j,
\quad |I_j|\le1,
\end{gathered}
\end{equation}
where $I_j$ is the horizontal projection of $K_j$. Indeed, horizontal projection does not increase diameter, and $0\in\overline I_j$, so $I_j\subset[-1,1]$. Write
\[
S_j(x)=(\ell_j(x),u_j(x)),
\qquad
h_j(x)=u_j(x)-\ell_j(x),
\qquad x\in I_j.
\]
Convexity of $K_j$ implies that $\ell_j$ is convex and $u_j$ is concave. With endpoint values understood by passage to the closure, the full-height chord gives $h_j(0)=H_j=\max_{\overline I_j}h_j.$
Thus the normalized transverse potential is
\[
W_j(x)=
\frac{\Lambda_sh_j(x)^{-2s}-\Lambda_sH_j^{-2s}}
{\rho_s(\varepsilon_j)},
\qquad W_j(0)=0.
\]

Let $\nu_{k,j}$ denote the eigenvalues associated with $\mathfrak a_j$ on $\Hs{s}{K_j}$, with $\eta_j=\varepsilon_j$, and set $c_j:=\Lambda_sH_j^{-2s}$. The unitary transverse rescaling $v(x,z)=\varepsilon_j^{1/2}u(x,\varepsilon_jz)$
preserves the $L^2$ norm and gives
\[
\begin{aligned}
\Es_s[u]
&=\iint_{\R^2}
\left(\xi^2+\varepsilon_j^{-2}\zeta^2\right)^s
|\wh v(\xi,\zeta)|^2\,\dd\xi\dd\zeta=\varepsilon_j^{-2s}\mathfrak a_j[v],
\end{aligned}
\]
where
\begin{equation}\label{eq:ajv}
    \mathfrak a_j[v]=\iint_{\R^2}
(\zeta^2+\varepsilon_j^2\xi^2)^s
|\wh v(\xi,\zeta)|^2\,\dd\xi\dd\zeta.
\end{equation}
The min--max principle therefore yields
\begin{equation}\label{eq:physical-anisotropic-spectrum}
\lambda_{k,s}(\Omega_j)=\varepsilon_j^{-2s}\nu_{k,j}.
\end{equation}

\begin{theorem}\label{thm:window}
Fix $0<s<1$ and $C>0$. Suppose that the domains satisfy \eqref{eq:geometric-normalization}--\eqref{eq:normalized-geometry} and
\[
\nu_{2,j}-\nu_{1,j}\le C\rho_s(\varepsilon_j).
\]
Then, for all sufficiently large $j$,
\begin{equation}\label{eq:window-bound}
0\le\nu_{1,j}-c_j\le C'\rho_s(\varepsilon_j),
\end{equation}
where $C'>0$ depends only on $s$ and $C$.
\end{theorem}

\begin{proof}
Set $\delta_j:=\nu_{1,j}-c_j$. By the transverse Poincar\'e inequality in \eqref{eq:transverse-potential},
for every admissible $v$, we obtain $\mathfrak a_j[v]\ge c_j\|v\|_2^2$, so $\delta_j\ge0$. We first prove $\delta_j\to0$, independently of the gap bound. We then rule out the possibility that $\delta_j/\rho_s(\varepsilon_j)$ is unbounded.

\step{Step 1. Convergence to the transverse eigenvalue.}
Since $|K_j|=1$ and $H_j\le2$, we have $|I_j|\ge1/2$. After a horizontal reflection, we may assume that the right endpoint $d_j$ of $I_j$ satisfies $d_j\ge1/4$. Choose $(d_j,z_j)\in\overline K_j$, where $0\le z_j\le H_j$.

By convexity, the triangle with vertices $(0,0)$, $(0,H_j)$, and $(d_j,z_j)$ lies in $\overline K_j$. Its vertical section at $0<x<d_j$ is
\[
\left(\frac{xz_j}{d_j},
H_j-\frac{x(H_j-z_j)}{d_j}\right).
\]
Since $H_j/d_j\le8$ and $K_j=\operatorname{int}\overline K_j$, we obtain, for every $0<r<1/32$ with $H_j-16r>1/2$,
\begin{equation}\label{eq:window-inner-product}
(0,r)\times(8r,H_j-8r)\subset K_j.
\end{equation}

Fix $f\in C_c^\infty(0,1)$ with $\|f\|_2=1$, and put $f_r(x):=r^{-1/2}f(x/r)$. Let $q_j^r$ be the normalized first Dirichlet eigenfunction on $(8r,H_j-8r)$. By \eqref{eq:window-inner-product}, the product $f_rq_j^r$ belongs to $\Hs{s}{K_j}$ and has unit $L^2$ norm. Using $(a+b)^s\le a^s+b^s$ for $a,b\ge0$, we obtain
\[
\begin{aligned}
\mathfrak a_j[f_rq_j^r]
&\le
\int_\R|\zeta|^{2s}|\wh q_j^r(\zeta)|^2\,\dd\zeta
+\varepsilon_j^{2s}
\int_\R|\xi|^{2s}|\wh f_r(\xi)|^2\,\dd\xi\\
&=\Lambda_s(H_j-16r)^{-2s}
+\left(\frac{\varepsilon_j}{r}\right)^{2s}
\Es_s[f].
\end{aligned}
\]
The min--max principle and the mean value theorem therefore give
\[
\begin{aligned}
0\le\delta_j
&\le\Lambda_s\bigl((H_j-16r)^{-2s}-H_j^{-2s}\bigr)
+\left(\frac{\varepsilon_j}{r}\right)^{2s}\Es_s[f]\le C_{s,f}\left(r+\left(\frac{\varepsilon_j}{r}\right)^{2s}\right).
\end{aligned}
\]
Taking $r=\sqrt{\varepsilon_j}$ for sufficiently large $j$, we obtain
\[0\le\delta_j
\le C_{s,f}\bigl(\varepsilon_j^{1/2}+\varepsilon_j^s\bigr)
\longrightarrow0.\]
Moreover, $\delta_j>0$. Indeed, \eqref{eq:horizontal-band-bound} and $B_s(k)>\Lambda_s$ for $k>0$ imply that $\mathfrak a_j[v]>c_j\|v\|_2^2$ for every nonzero admissible $v$. Since $\nu_{1,j}$ is attained, the inequality is strict at its first eigenfunction.

\step{Step 2. Rescaling at the first-eigenvalue excess.}
Suppose, along a subsequence, that $\frac{\delta_j}{\rho_s(\varepsilon_j)}
\rightarrow\infty.$
Since $\delta_j>0$ and $\delta_j\to0$, and $\rho_s$ is strictly increasing near zero, we may define $\eta_j$ and $r_j$ by $\rho_s(\eta_j)=\delta_j,\,
r_j=\frac{\varepsilon_j}{\eta_j}.$
Then $\eta_j\to0$ and $\eta_j>\varepsilon_j$ for all large $j$.

For $0<s<1/2$, we have $r_j=
\left(\frac{\rho_s(\varepsilon_j)}{\delta_j}\right)^{1/(1+2s)}
\rightarrow0,$
while for $1/2<s<1$, $r_j=
\left(\frac{\rho_s(\varepsilon_j)}{\delta_j}\right)^{1/2}
\rightarrow0.$
At $s=1/2$, the identity
$\delta_j=\eta_j^2\log(1/\eta_j)$ gives
\[
\begin{aligned}
0<r_j^2
&=\frac{\varepsilon_j^2\log(1/\eta_j)}{\delta_j}\le\frac{\varepsilon_j^2\log(1/\varepsilon_j)}{\delta_j}
=\frac{\rho_{1/2}(\varepsilon_j)}{\delta_j}
\longrightarrow0.
\end{aligned}
\]
Thus $r_j\to0$. Passing to a further subsequence, we may assume $H_j\to H\in[1,2]$.

Define
\[
\widetilde K_j
:=\{(t,z):(r_jt,z)\in K_j\},
\qquad
w(t,z):=r_j^{1/2}v(r_jt,z).
\]
This is a unitary change of variables. Since $\varepsilon_j/r_j=\eta_j$, by \eqref{eq:ajv}, Fourier scaling gives
\[
\begin{aligned}
\mathfrak a_j[v]
&=\iint_{\R^2}
(\zeta^2+\eta_j^2\xi^2)^s
|\wh w(\xi,\zeta)|^2\,\dd\xi\dd\zeta=:\widetilde{\mathfrak a}_j[w].
\end{aligned}
\]
Set
\[
\widetilde F_j[w]
:=\frac{\widetilde{\mathfrak a}_j[w]-c_j\|w\|_2^2}
{\delta_j}.
\]
The first two eigenvalues of $\widetilde F_j$ satisfy
\begin{equation}\label{eq:clustered-rescaled-spectrum}
\mu_{1,j}=1,
\qquad
\mu_{2,j}
=1+\frac{\nu_{2,j}-\nu_{1,j}}{\delta_j}
\longrightarrow1,
\end{equation}
where the convergence follows from the assumed gap bound.

The transverse section of $\widetilde K_j$ at $t$ has length $h_j(r_jt)$. Define the corresponding normalized potential by
\[
U_j(t)
:=\frac{\Lambda_sh_j(r_jt)^{-2s}-c_j}{\delta_j}
\]
on the rescaled projection, with lower semicontinuous endpoint values and value $+\infty$ outside its closure. Since $h_j$ is positive and concave and $r\mapsto r^{-2s}$ is convex and decreasing, $U_j$ is convex. Moreover, $U_j\ge0,\, U_j(0)=0.$

\step{Step 3. Uniform confinement after rescaling.}
Fix $\tau>0$, to be chosen below, and define
$E_j:=\{t\in\R:U_j(t)\le\tau\}$.
Each $E_j$ is an interval containing zero. For $t\in E_j$, we have
$0\le \Lambda_s\bigl(h_j(r_jt)^{-2s}-H_j^{-2s}\bigr) \le\tau\delta_j$.
By convexity of $r\mapsto r^{-2s}$ and $1\le H_j\le2$,
\[
\begin{aligned}
h_j(r_jt)^{-2s}-H_j^{-2s}
&\ge2sH_j^{-2s-1}
\bigl(H_j-h_j(r_jt)\bigr),
\end{aligned}
\]
and hence
\[
0\le H_j-h_j(r_jt)\le A_s\tau\delta_j,
\qquad
A_s:=\frac{2^{2s+1}}{2s\Lambda_s}.
\]
Since
$H_j-h_j(r_jt) =\ell_j(r_jt)+H_j-u_j(r_jt)$
and both terms on the right are nonnegative, we obtain
\[
0\le\ell_j(r_jt)\le A_s\tau\delta_j,
\qquad
0\le H_j-u_j(r_jt)\le A_s\tau\delta_j.
\]

Put $\vartheta_j:=A_s\tau\delta_j$ and
$S_j^*:=(\vartheta_j,H_j-\vartheta_j)$.
For all sufficiently large $j$, the interval $S_j^*$ is nonempty and lies in every transverse section over $\operatorname{int}E_j$. Consequently,
\begin{equation}\label{eq:sublevel-product}
\operatorname{int}E_j\times S_j^*
\subset\widetilde K_j.
\end{equation}

We claim that $|E_j|$ is uniformly bounded for a suitable fixed $\tau$. Suppose otherwise. After passing to a subsequence, $|E_j|\to\infty$. For any fixed $L>0$, we may then choose $a_j\in\R$ such that $(a_j-L/2,a_j+L/2)\subset\operatorname{int}E_j.$
Fix $f\in C_c^\infty(I)$ with $\|f\|_2=1$, and define
\[
f_{L,a_j}(t):=L^{-1/2}
f\left(\frac{t-a_j}{L}\right).
\]
Let $q_j^*$ be the normalized first Dirichlet eigenfunction on $S_j^*$. By \eqref{eq:sublevel-product}, the product $f_{L,a_j}q_j^*$ is admissible for $\widetilde F_j$ and has unit $L^2$ norm. Its normalized transverse energy satisfies
\[
\begin{aligned}
&\frac{\mathfrak a_0[f_{L,a_j}q_j^*]
-c_j\|f_{L,a_j}q_j^*\|_2^2}{\delta_j}=
\frac{\Lambda_s(H_j-2\vartheta_j)^{-2s}
-\Lambda_sH_j^{-2s}}{\delta_j}
\le \widetilde C_s\tau,
\end{aligned}
\]
where the constant is independent of $j,\tau$, and $L$ for all sufficiently large $j$.

Since $H_j\to H$ and $\vartheta_j\to0$, scaling of the interval ground states gives $q_j^*\rightarrow q_{s,H}$ in $L^1(\R)\cap L^2(\R).$
Moreover,
\[
\wh f_{L,a_j}(\xi)
=L^{1/2}e^{-ia_j\xi}\wh f(L\xi),
\qquad
|\wh f_{L,a_j}(\xi)|^2
=L|\wh f(L\xi)|^2.
\]
The Fourier modulus is therefore independent of the translation $a_j$. Applying Lemma~\ref{lem:product-limit}, with $\rho_s(\eta_j)=\delta_j$, we obtain
\[
\begin{aligned}
&\frac1{\delta_j}
\iint_{\R^2}D_{s,\eta_j|\xi|}(\zeta)
|\wh f_{L,a_j}(\xi)|^2
|\wh q_j^*(\zeta)|^2\,\dd\xi\dd\zeta\longrightarrow
a_s(H)L^{-2\sigma_s}
\int_\R|\xi|^{2\sigma_s}|\wh f(\xi)|^2\,\dd\xi.
\end{aligned}
\]

Since $a_s(H)$ is bounded uniformly for $H\in[1,2]$, first choose $\tau>0$ so small that $\widetilde C_s\tau<1/4$, and then choose $L>0$ sufficiently large that
\[
\sup_{1\le H\le2}a_s(H)L^{-2\sigma_s}
\int_\R|\xi|^{2\sigma_s}|\wh f(\xi)|^2\,\dd\xi
<\frac14.
\]
For these choices, the Rayleigh quotient of $f_{L,a_j}q_j^*$ with respect to $\widetilde F_j$ satisfies
\[
\widetilde F_j[f_{L,a_j}q_j^*]
\le \widetilde C_s\tau+\frac14+o(1)<1
\]
for all sufficiently large $j$. The min--max principle then gives $\mu_{1,j}<1$, contradicting \eqref{eq:clustered-rescaled-spectrum}. Thus there exists $L_0>0$, independent of $j$, such that $E_j\subset[-L_0,L_0]$
for all sufficiently large $j$, where we also used $0\in E_j$. Set $R:=L_0+1$. For $|t|>R$ with $U_j(t)<\infty$, convexity and $U_j(0)=0$ give
$\tau<U_j(R\operatorname{sgn}t) \le\frac R{|t|}U_j(t)$.
Together with $U_j\ge0$ on $[-R,R]$, this yields
\begin{equation}\label{eq:second-confinement}
U_j(t)\ge\frac{\tau}{R}|t|-\tau,
\qquad t\in\R,
\end{equation}
for all sufficiently large $j$. In particular, the sublevel sets of $U_j$ are uniformly bounded.

\step{Step 4. The limiting contradiction.}
By Lemma~\ref{lem:convex-compactness}, after passing to a further subsequence, $U_j\econv U$, where $U$ is proper, lower semicontinuous, convex, and satisfies $U(0)=0$. The confinement estimate \eqref{eq:second-confinement} passes to the epi-limit. Thus the rescaled domains $\widetilde K_j$, the parameters $\eta_j\to0$, and the potentials $U_j$ satisfy the hypotheses of Theorem~\ref{thm:strip-limit}.

Let $w_{1,j},w_{2,j}$ be orthonormal eigenfunctions of $\widetilde F_j$ corresponding to $\mu_{1,j},\mu_{2,j}$. By \eqref{eq:clustered-rescaled-spectrum}, their energies are uniformly bounded:
\[
\widetilde F_j[w_{1,j}]=1,
\qquad
\widetilde F_j[w_{2,j}]=\mu_{2,j}\le2
\]
for all sufficiently large $j$. The compactness assertion of Theorem~\ref{thm:strip-limit} gives, along a common subsequence,
\[
w_{k,j}\longrightarrow\psi_kq_{s,H}
\quad\text{strongly in }L^2(\R^2),
\qquad k=1,2.
\]
Since $\|q_{s,H}\|_2=1$, strong convergence preserves orthonormality:
\[
\langle\psi_k,\psi_l\rangle_{L^2(\R)}
=\lim_{j\to\infty}
\langle w_{k,j},w_{l,j}\rangle_{L^2(\R^2)}
=\delta_{kl}.
\]
In particular, both $\psi_k$ are nonzero. Hence
$J:=\operatorname{int}\dom U$ cannot be empty, since otherwise Theorem~\ref{thm:strip-limit} would force the limiting functions to vanish.

Corollary~\ref{cor:strip-spectrum} now applies. By \eqref{eq:clustered-rescaled-spectrum}, the first two eigenvalues of the limiting one-dimensional form satisfy
\[
\lambda_1(Q_{\sigma_s,a_s(H),U})
=\lambda_2(Q_{\sigma_s,a_s(H),U})=1.
\]
However, $0<\sigma_s\le1$, $a_s(H)>0$, and $U$ is a convex potential on $J$ with bounded sublevel sets. Corollary~\ref{cor:ground-simple} implies that its first eigenvalue is simple, a contradiction.

We have shown that $\delta_j/\rho_s(\varepsilon_j)$ remains bounded along every sequence satisfying the hypotheses. We now prove that the bound is uniform over all admissible convex domains. Suppose that no constant $C'=C'(s,C)$ gives such a uniform bound for sufficiently small areas. Then, for each $j\ge1$, we could choose an admissible convex domain with $\varepsilon_j<1/j$ such that
\[
\nu_{2,j}-\nu_{1,j}\le C\rho_s(\varepsilon_j),
\qquad
\frac{\nu_{1,j}-c_j}{\rho_s(\varepsilon_j)}>j.
\]
This would give $\varepsilon_j\to0$ and $\frac{\nu_{1,j}-c_j}{\rho_s(\varepsilon_j)}
\rightarrow\infty,$
contradicting the preceding argument. Hence \eqref{eq:window-bound} holds uniformly over all admissible convex domains with sufficiently small area, with $C'$ depending only on $s$ and $C$.

\end{proof}

\section{Optimization over planar convex domains}
\label{sec:optimization}

We now turn to the optimization of the fundamental gap over planar convex domains of prescribed diameter and area. Using convexity and the compact spectral convergence established in the preceding section, we reduce the asymptotic optimization problem to a one-dimensional eigenvalue gap problem with a convex potential. We then establish the sharp asymptotic bounds by constructing matching competitors and characterize the geometry of asymptotically minimizing domains.

For later use, define
\begin{equation}\label{eq:physical-gap-scale}
 f_s(\theta):=\theta^{-2s}\rho_s(\theta)
=
\begin{cases}
\theta,&0<s<\frac12,\\
\theta\log(1/\theta),&s=\frac12,\\
\theta^{2-2s},&\frac12<s<1.
\end{cases}   
\end{equation}

\subsection{Upper bounds and prescribed effective potentials}

\step{Step 1. Thin rectangles.}
Consider $R_{1,h}=I\times(-h/2,h/2)$, with $0<h<1$. Its diameter and area ratio are
\[
D_h=\sqrt{1+h^2},\qquad
\theta_h:=\frac{|R_{1,h}|}{D_h^2}=\frac{h}{1+h^2}.
\]
By Theorem~\ref{thm:rectangles}, the common transverse term $\Lambda_sh^{-2s}$ cancels in the fundamental gap, giving
\[
\gap_s(R_{1,h})=
\begin{cases}
\kappa_s\gap_{s+1/2}(I,0)h+o(h),
&0<s<\frac12,\\[1mm]
3\pi^2\gamma h\log(1/h)
+o\bigl(h\log(1/h)\bigr),
&s=\frac12,\\[1mm]
3\pi^2\beta_s h^{2-2s}
+o(h^{2-2s}),
&\frac12<s<1.
\end{cases}
\]
Since $D_h\to1$, $\theta_h/h\to1$, and
\[
\frac{\log(1/\theta_h)}{\log(1/h)}
=
1+\frac{\log(1+h^2)}{\log(1/h)}
\longrightarrow1,
\]
we obtain
\begin{equation}\label{eq:rectangle-normalized-gap}
\frac{D_h^{2s}\gap_s(R_{1,h})}{f_s(\theta_h)}
\longrightarrow
\begin{cases}
\kappa_s\gap_{s+1/2}(I,0),
&0<s<\frac12,\\[1mm]
3\pi^2\gamma,&s=\frac12,\\
3\pi^2\beta_s,&\frac12<s<1.
\end{cases}
\end{equation}

The map $h\mapsto h/(1+h^2)$ is strictly increasing from $(0,1)$ onto $(0,1/2)$,
thus every sufficiently small area ratio is realized by a unique such rectangle. Since rectangles are admissible in the definition of $G_s$, \eqref{eq:rectangle-normalized-gap} gives
\[
\begin{aligned}
\limsup_{\theta\downarrow0}
\frac{G_{1/2}(\theta)}
{\theta\log(1/\theta)}
&\le3\pi^2\gamma,\qquad
\limsup_{\theta\downarrow0}
\frac{G_s(\theta)}{\theta^{2-2s}}\le3\pi^2\beta_s,
\quad \frac12<s<1.
\end{aligned}
\]
For $0<s<1/2$, rectangles also give an upper bound of order $\theta$, but their coefficient need not be optimal among convex domains with variable widths.

\step{Step 2. Realization of convex potentials in the subcritical range.}
Fix $0<s<1/2$ and set $p:=1+2s$. Let $V\in C(\overline I)$ be convex. Since adding a constant to $V$ does not change the one-dimensional spectral gap, we may assume $V\ge0,\,\min_{\overline I}V=0.$
Define
\begin{equation}\label{eq:width-realization}
b_s:=\frac{\kappa_s}{2s\Lambda_s},
\qquad
H_\varepsilon:=1+b_s\varepsilon^p\int_I V(x)\,\dd x,
\qquad
h_\varepsilon(x):=H_\varepsilon-b_s\varepsilon^pV(x),
\end{equation}
and let
$K_\varepsilon :=\{(x,z):x\in I,\ 0<z<h_\varepsilon(x)\}$.
Since $V$ is bounded and convex, $h_\varepsilon$ is concave, strictly positive for sufficiently small $\varepsilon$, and converges uniformly to $1$. Thus $K_\varepsilon$ is a bounded open convex domain. Moreover,
\[
\begin{aligned}
|K_\varepsilon|
&=\int_Ih_\varepsilon(x)\,\dd x
=1,\qquad
\sup_{x\in I}h_\varepsilon(x)=H_\varepsilon\rightarrow1.
\end{aligned}
\]
The second identity uses $\inf_I V=0$, with the supremum understood in the usual sense when the minimum is attained only at an endpoint.

We next identify the effective potential. Taylor's formula for $t\mapsto t^{-2s}$, applied uniformly on $\overline I$, gives
\[
\begin{aligned}
h_\varepsilon(x)^{-2s}-H_\varepsilon^{-2s}
&=2sb_sH_\varepsilon^{-2s-1}
\varepsilon^pV(x)+O\bigl(\varepsilon^{2p}\|V\|_\infty^2\bigr).
\end{aligned}
\]
Consequently,
\[\begin{aligned}
W_\varepsilon(x)
&:=\frac{\Lambda_s
\bigl(h_\varepsilon(x)^{-2s}
-H_\varepsilon^{-2s}\bigr)}
{\varepsilon^p}=2s\Lambda_sb_sH_\varepsilon^{-2s-1}V(x)
+O(\varepsilon^p)\rightarrow\kappa_sV(x)
\end{aligned}\]
uniformly on $\overline I$, since $2s\Lambda_sb_s=\kappa_s$ and $H_\varepsilon\to1$.

Extend $W_\varepsilon$ to $\R$ by its continuous endpoint values on $\overline I$ and by $+\infty$ outside $\overline I$. The uniform convergence on $\overline I$ implies epi-convergence to
\[
W(x)=
\begin{cases}
\kappa_sV(x),&x\in\overline I,\\
+\infty,&x\notin\overline I.
\end{cases}
\]
The projection interval is fixed and bounded, so Theorem~\ref{thm:strip-limit} applies with $\eta=\varepsilon$, $H=1$, and $J=I$. Since $a_s(1)=\kappa_s$ and $\sigma_s=s+1/2$, for $\psi\in\Hs{s+1/2}{I}\cap L^2(I,V\,\dd x)$, the limiting form is
\[
\begin{aligned}
F[\psi]
&=\kappa_s\int_\R|\xi|^{1+2s}
|\wh\psi(\xi)|^2\,\dd\xi
+\kappa_s\int_I V(x)|\psi(x)|^2\,\dd x=\kappa_sQ_{s+1/2,1,V}[\psi].
\end{aligned}
\]

Let $\nu_k(\varepsilon)$ denote the $k$th eigenvalue of the anisotropic form on $K_\varepsilon$. By Corollary~\ref{cor:strip-spectrum}, for every fixed $k\ge1$,
\[
\frac{\nu_k(\varepsilon)
-\Lambda_sH_\varepsilon^{-2s}}
{\varepsilon^p}
\longrightarrow
\kappa_s\lambda_{k,s+1/2}(I,V).
\]
Subtracting the limits for $k=1,2$, we obtain
\[
\nu_2(\varepsilon)-\nu_1(\varepsilon)
=
\kappa_s\gap_{s+1/2}(I,V)
\varepsilon^p+o(\varepsilon^p).
\]

\step{Step 3. Physical scaling and exact area ratios.}
Define the physical thin domain
\[
\Omega_\varepsilon
:=\{(x,\varepsilon z):(x,z)\in K_\varepsilon\},
\qquad
D_\varepsilon:=\diam\Omega_\varepsilon.
\]
Since $|K_\varepsilon|=1$ and $H_\varepsilon\to1$,
\[
|\Omega_\varepsilon|=\varepsilon,
\qquad
1\le D_\varepsilon
\le\sqrt{1+\varepsilon^2H_\varepsilon^2}.
\]
Consequently, its area ratio satisfies
\[
\theta_\varepsilon
:=\frac{|\Omega_\varepsilon|}{D_\varepsilon^2}
=\frac{\varepsilon}{D_\varepsilon^2}
=\varepsilon+O(\varepsilon^3).
\]

The anisotropic scaling gives $\lambda_{k,s}(\Omega_\varepsilon)
=\varepsilon^{-2s}\nu_k(\varepsilon).$
Hence, using $p=1+2s$,
\[
\begin{aligned}
\gap_s(\Omega_\varepsilon)
&=\varepsilon^{-2s}
\bigl(\nu_2(\varepsilon)-\nu_1(\varepsilon)\bigr)=\kappa_s\gap_{s+1/2}(I,V)
\varepsilon^{p-2s}
+o(\varepsilon^{p-2s})\\
&=\kappa_s\gap_{s+1/2}(I,V)
\varepsilon+o(\varepsilon).
\end{aligned}
\]
Since $\theta_\varepsilon=\varepsilon/D_\varepsilon^2$, we obtain
\begin{equation}\label{eq:realization-gap}
\begin{aligned}
\frac{D_\varepsilon^{2s}
\gap_s(\Omega_\varepsilon)}
{\theta_\varepsilon}
&=D_\varepsilon^{2s+2}
\bigl(\kappa_s\gap_{s+1/2}(I,V)+o(1)\bigr)\rightarrow\kappa_s\gap_{s+1/2}(I,V).
\end{aligned}
\end{equation}

It remains to realize every sufficiently small area ratio. Define
\[
\Phi_\varepsilon(x,t)
:=(x,\varepsilon t h_\varepsilon(x)),
\qquad (x,t)\in\overline I\times[0,1].
\]
Then $\overline\Omega_\varepsilon
=\Phi_\varepsilon(\overline I\times[0,1]).$
Since $\Phi_\varepsilon$ depends continuously on $\varepsilon$ in the uniform norm,
\[
|D_\varepsilon-D_{\varepsilon'}|
\le2\|\Phi_\varepsilon-\Phi_{\varepsilon'}\|_\infty,
\]
and therefore $\varepsilon\mapsto\theta_\varepsilon$ is continuous.

For every sufficiently small $\theta>0$, the diameter bounds give
\[
\theta_{\theta/2}
\le\frac{\theta}{2}<\theta,
\qquad
\theta_{2\theta}
\ge\frac{2\theta}
{1+4\theta^2H_{2\theta}^2}>\theta.
\]
By the intermediate value theorem, there exists
$\varepsilon(\theta)\in[\theta/2,2\theta]$
such that $\theta_{\varepsilon(\theta)}=\theta.$
Since $\varepsilon(\theta)\to0$ as $\theta\downarrow0$,
\eqref{eq:realization-gap} and the definition of $G_s$ imply
\[
\limsup_{\theta\downarrow0}
\frac{G_s(\theta)}{\theta}
\le\kappa_s\gap_{s+1/2}(I,V).
\]
Taking the infimum over all convex
$V\in C(\overline I)$, we obtain
\[\limsup_{\theta\downarrow0}
\frac{G_s(\theta)}{\theta}
\le\kappa_s\gs_{s+1/2}.\]

\subsection{Lower bounds in the three ranges}

We now establish the matching lower bounds by applying the compact strip convergence to nearly minimizing convex domains.

\begin{proof}[\textbf{Completion of the proof of Theorem~\ref{thm:main}}.]
\step{Step 1. Nearly minimizing domains.}
Choose $\theta_j\downarrow0$ such that
\[
\frac{G_s(\theta_j)}{f_s(\theta_j)}
\longrightarrow
\liminf_{\theta\downarrow0}
\frac{G_s(\theta)}{f_s(\theta)}.
\]
The preceding upper bounds ensure that this limit is finite. By the definition of $G_s$ and dilation invariance, we may choose bounded open convex domains $\Omega_j$ satisfying
\begin{equation}\label{eq:near-minimizers}
\diam\Omega_j=1,\qquad
|\Omega_j|=\theta_j,\qquad
G_s(\theta_j)\le\gap_s(\Omega_j)
\le G_s(\theta_j)+\frac{f_s(\theta_j)}j.
\end{equation}
Consequently,
\[
0\le
\frac{\gap_s(\Omega_j)-G_s(\theta_j)}
{f_s(\theta_j)}
\le\frac1j,
\]
so the normalized gaps of $\Omega_j$ converge to the chosen lower limit.

Apply the minimal-width normalization
\eqref{eq:geometric-normalization} with
$\varepsilon_j=\theta_j$, and write
$\rho_j:=\rho_s(\theta_j)$. By
\eqref{eq:physical-anisotropic-spectrum}
and \eqref{eq:physical-gap-scale},
\[
\frac{\nu_{2,j}-\nu_{1,j}}{\rho_j}
=
\frac{\gap_s(\Omega_j)}{f_s(\theta_j)}
\le C.
\]
Theorem~\ref{thm:window} therefore gives $0\le\nu_{1,j}-c_j\le C'\rho_j.$
Combining this with the spectral-gap bound, we obtain
\begin{equation}\label{eq:two-eigenvalue-window}
0\le\nu_{k,j}-c_j\le C''\rho_j,
\qquad k=1,2.
\end{equation}
Thus the first two normalized eigenfunctions have
uniformly bounded shifted energies, as required
for the compactness result in
Theorem~\ref{thm:strip-limit}.

\step{Step 2. A nonempty limiting interval.}
Write $I_j=(a_j,b_j)$. By
\eqref{eq:normalized-geometry},
\[
1\le H_j\le2,\qquad
I_j\subset[-1,1],\qquad |I_j|\le1.
\]
Moreover, $W_j\ge0$ and $W_j(0)=0$.
By Lemma~\ref{lem:convex-compactness}, after passing to a subsequence,
\[
H_j\to H\in[1,2],\qquad
a_j\to a,\qquad b_j\to b,\qquad
W_j\econv W.
\]
Since $W_j=+\infty$ outside $\overline I_j$,
the recovery-point property of epi-convergence gives
${\dom W\subset[a,b],\enspace b-a\le1}$.

Let $v_{1,j},v_{2,j}$ be orthonormal eigenfunctions corresponding to $\nu_{1,j},\nu_{2,j}$. By \eqref{eq:two-eigenvalue-window}, for $k=1,2$,
\[
\|v_{k,j}\|_2=1,\qquad
F_j[v_{k,j}]
=\frac{\nu_{k,j}-c_j}{\rho_j}\le C''.
\]
The compactness assertion of Theorem~\ref{thm:strip-limit} therefore gives, after a further common subsequence,
\[
v_{k,j}\longrightarrow\psi_kq_{s,H}
\quad\text{strongly in }L^2(\R^2),
\qquad k=1,2,
\]
with $\psi_k\in\mathcal D(F)$. Since
$\|q_{s,H}\|_2=1$, orthonormality passes to the limit:
\[
\begin{aligned}
\langle\psi_k,\psi_l\rangle_{L^2(\R)}
&=\lim_{j\to\infty}
\langle v_{k,j},v_{l,j}\rangle_{L^2(\R^2)}
=\delta_{kl}.
\end{aligned}
\]
Thus the limit functions are nonzero. If
$J:=\operatorname{int}\dom W$ were empty,
then $\mathcal D(F)=\{0\}$, a contradiction.
Consequently, $J\ne\varnothing,\,
0<L:=|J|\le b-a\le1.$

Corollary~\ref{cor:strip-spectrum} now applies and yields
\[
\frac{\nu_{2,j}-\nu_{1,j}}{\rho_j}
\longrightarrow\mu_2-\mu_1,
\]
where $\mu_k$ are the eigenvalues of the effective form \eqref{eq:effective-form}. Combining this with \eqref{eq:near-minimizers} and the scaling identity in Step~1 gives
\begin{equation}\label{eq:optimized-limit-gap}
\liminf_{\theta\downarrow0}
\frac{G_s(\theta)}{f_s(\theta)}
=
\lim_{j\to\infty}
\frac{\gap_s(\Omega_j)}{f_s(\theta_j)}
=\mu_2-\mu_1.
\end{equation}

\step{Step 3. The one-dimensional lower bounds.}
The limiting potential $W$ is finite and convex on the bounded interval $J$ of length $L\le1$. We apply Lemma~\ref{lem:convex-gap-bound}, which also allows convex potentials diverging at the endpoints.

If $0<s<1/2$, set $\sigma=s+1/2$. Since $F=\kappa_sH\,Q_{\sigma,1,W/(\kappa_sH)},$
Lemma~\ref{lem:convex-gap-bound} gives
\[
\begin{aligned}
\mu_2-\mu_1
&\ge\kappa_sH L^{-2\sigma}\gs_\sigma\ge\kappa_s\gs_{s+1/2}.
\end{aligned}
\]
If $s=1/2$, the effective form has kinetic coefficient $\gamma H$, and hence
\[
\mu_2-\mu_1
\ge3\pi^2\gamma\frac{H}{L^2}
\ge3\pi^2\gamma.
\]
Finally, if $1/2<s<1$, its kinetic coefficient is $\beta_sH^{2-2s}$, giving
\[
\mu_2-\mu_1
\ge3\pi^2\beta_s\frac{H^{2-2s}}{L^2}
\ge3\pi^2\beta_s.
\]
Here we used $H\ge1$, $L\le1$, and the positivity of $\gs_\sigma$ and the coefficients.

Together with \eqref{eq:optimized-limit-gap}, these inequalities establish the lower bounds in \eqref{eq:main-asymptotics}. The upper bounds proved in the preceding subsection give the corresponding equalities. Finally, \eqref{eq:rectangle-normalized-gap} shows that thin rectangles attain the leading constant for $s\ge1/2$; for $0<s<1/2$, they do so precisely when $\gs_{s+1/2}=\gap_{s+1/2}(I,0)$. This completes the proof.
\end{proof}

\subsection{The shape of asymptotically minimizing domains}

We show that asymptotically minimizing convex domains converge, after minimal-width normalization, to a rectangle, although their limiting effective potentials need not vanish.

\begin{proof}[\textbf{Proof of Theorem~\ref{thm:shape}}.]
\step{Step 1. The limiting height.}
Consider an arbitrary subsequence of the normalized domains. By asymptotic optimality and Theorem~\ref{thm:main}, their normalized gaps are bounded. The compactness and spectral convergence argument in the proof of Theorem~\ref{thm:main} therefore gives, after passing to a further subsequence,
${H_j\to H\in[1,2],\enspace 0<L:=|J|\le1}$,
with the normalized gaps converging to $\mu_2-\mu_1$.

We now use asymptotic optimality and the one-dimensional gap bounds to determine $H$ and $L$. For $0<s<1/2$, we have
\[
\kappa_s\mathfrak g_{s+1/2}
=\mu_2-\mu_1
\ge \kappa_s H L^{-1-2s}\mathfrak g_{s+1/2}.
\]
At the critical exponent $s=1/2$, the classical fundamental gap bound gives
\[
3\pi^2\gamma
=\mu_2-\mu_1
\ge 3\pi^2\gamma\frac{H}{L^2}.
\]
Similarly, for $1/2<s<1$,
\[
3\pi^2\beta_s
=\mu_2-\mu_1
\ge 3\pi^2\beta_s\frac{H^{2-2s}}{L^2}.
\]
Since $H\ge1$, $L\le1$, and all the coefficients are strictly positive, each inequality forces $H=L=1$.
Thus every subsequence of $(H_j)$ has a further subsequence converging to $1$. Consequently, $H_j\rightarrow1$
along the full sequence.

\step{Step 2. Convergence of the transverse boundaries.}
Since $|K_j|=1$, $h_j\le H_j$, and $|I_j|\le1$,
\[
1=\int_{I_j}h_j(x)\,\dd x
\le H_j|I_j|\le H_j.
\]
Together with $H_j\to1$, this gives
\begin{equation}\label{eq:height-deficit}
|I_j|\longrightarrow1,
\qquad
\int_{I_j}(H_j-h_j(x))\,\dd x
=H_j|I_j|-1\longrightarrow0.
\end{equation}

After a horizontal translation, assume $I_j=\left(-\frac{|I_j|}{2},\frac{|I_j|}{2}\right).$
Then every compact $A\Subset I$ is contained in $I_j$ for all sufficiently large $j$.

Set $d_j(x):=H_j-h_j(x)$. Since $h_j$ is concave and $h_j\le H_j$, the function $d_j$ is nonnegative and convex. Fix $A\Subset I$ and choose $r>0$ such that $\operatorname{dist}(A,\partial I)>2r.$
For all sufficiently large $j$, the interval $[x-r,x+r]$ is contained in $I_j$ for every $x\in A$. Jensen's inequality gives
\[
\begin{aligned}
0\le d_j(x)
&\le\frac1{2r}\int_{x-r}^{x+r}d_j(y)\,\dd y\le\frac1{2r}\int_{I_j}d_j(y)\,\dd y.
\end{aligned}
\]
Hence, by \eqref{eq:height-deficit},
\[
\sup_{x\in A}d_j(x)
\le\frac1{2r}\int_{I_j}d_j(y)\,\dd y
\longrightarrow0.
\]
Thus $d_j\to0$ locally uniformly on $I$. Since
\[
d_j(x)=\ell_j(x)+(H_j-u_j(x)),
\qquad
\ell_j(x)\ge0,\quad H_j-u_j(x)\ge0,
\]
we have
$0\le\ell_j(x)\le d_j(x), \enspace |u_j(x)-1|\le d_j(x)+|H_j-1|$.
Consequently,
\[
\ell_j\longrightarrow0,
\qquad
u_j\longrightarrow1
\quad\text{locally uniformly on }I.
\]

\step{Step 3. Hausdorff convergence.}
Let $R:=\overline I\times[0,1]$. Since the intervals $I_j$ are centered and $|I_j|\le1$, $\overline K_j\subset\overline I\times[0,H_j].$
As $H_j\ge1$, it follows that
\[
\sup_{p\in\overline K_j}\operatorname{dist}(p,R)
\le H_j-1\longrightarrow0.
\]

For the reverse estimate, fix $r\in(0,1/2)$ and set
${I_r:=\left[-\frac12+r,\frac12-r\right]}$.
For any $(x,z)\in R$, let $x_r$ be the projection of $x$ onto $I_r$ and define
$z_j:=\max\{\ell_j(x_r),\min\{z,u_j(x_r)\}\}$.
For all sufficiently large $j$, we have $I_r\subset I_j$, and $(x_r,z_j)\in\overline K_j.$
Moreover,
\[
|x-x_r|\le r,
\qquad
|z-z_j|
\le\max\{\ell_j(x_r),|u_j(x_r)-1|\}.
\]
Therefore,
\[
\begin{aligned}
\sup_{p\in R}\operatorname{dist}(p,\overline K_j)
&\le r+\max\left\{
\sup_{x\in I_r}\ell_j(x),
\sup_{x\in I_r}|u_j(x)-1|
\right\}.
\end{aligned}
\]
For each fixed $r$, the second term tends to zero by Step~2. Hence
\[
\limsup_{j\to\infty}
\sup_{p\in R}\operatorname{dist}(p,\overline K_j)
\le r.
\]
Letting $r\downarrow0$ proves the reverse Hausdorff bound, thus, we complete the proof.
\end{proof}

\section{AI disclosure}

This work was motivated by the thin-rectangle gap estimates of Dyda and Kulczycki \cite{DydaKulczycki2007} and the fractional fundamental gap problems studied by Bao, Ruan, Shen and Sheng \cite{BaoRuanShenSheng2019}. The author initially sought to determine the sharp asymptotic coefficients for collapsing rectangles and their dependence on the relative transverse width. This led to the formulation of the optimization problem $G_s(\theta)$, with $\theta=|\Omega|/(\diam\Omega)^2$, over general planar convex domains. During the research and preparation of this manuscript, the author used OpenAI's GPT-6 model for the following purposes:
\begin{enumerate}
\item[(i)] Assisting with mathematical calculations, checking the accuracy of my arguments and consistency of notation, and improving the language and presentation of the manuscript.

\item[(ii)] Exploring approaches to the construction of recovery sequences in part~(iii) of Theorem~\ref{thm:strip-limit} and to the second longitudinal rescaling argument in Section~\ref{sec:window}. The suggested ideas were critically examined, refined, and developed into the final proofs by the author.
\end{enumerate}

The author's principal contributions are as follows:
\begin{enumerate}
\item[(i)] Formulating the research questions and the four principal results, reconstructing and independently verifying the mathematical proofs, and preparing their final exposition.

\item[(ii)] Writing and organizing the Introduction, including the historical background, the main results, the key ideas of the proofs, and the discussion in Subsection~\ref{remark}, with the aim of making the mathematical arguments and their connections clear to readers.

\item[(iii)] Organizing the manuscript and its proofs, selecting the order of results and arguments, and refining the mathematical presentation to ensure rigor, clarity, and readability.
\end{enumerate}

The author takes full responsibility for all mathematical statements, proofs, calculations, and references in the manuscript.

\medskip
\noindent\textbf{Conflict of interest.}
The author declares no conflict of interest.

\medskip
\noindent\textbf{Data availability.}
No datasets were generated or analyzed during this study.

	\bibliographystyle{amsplain}

\begin{thebibliography}{99}

\bibitem{AmatoBucurFragala2025}
Vincenzo Amato, Dorin Bucur, and Ilaria Fragal\`a,
\emph{The geometric size of the fundamental gap},
2024, revised 2025, arXiv:2407.01341v2.

\bibitem{AmbrosioFreddiMusina2020}
Vincenzo Ambrosio, Lorenzo Freddi, and Roberta Musina,
\emph{Asymptotic analysis of the Dirichlet fractional Laplacian
in domains becoming unbounded},
J. Math. Anal. Appl. \textbf{485} (2020), no.~2, 123845.

\bibitem{AndrewsClutterbuck2011}
Ben Andrews and Julie Clutterbuck,
\emph{Proof of the fundamental gap conjecture},
J. Amer. Math. Soc. \textbf{24} (2011), no.~3, 899--916.

\bibitem{AndrewsClutterbuckHauer2021}
Ben Andrews, Julie Clutterbuck, and Daniel Hauer,
\emph{The fundamental gap for a one-dimensional
Schr\"odinger operator with Robin boundary conditions},
Proc. Amer. Math. Soc. \textbf{149} (2021), no.~4, 1481--1493.

\bibitem{AshbaughBenguria1989}
Mark S. Ashbaugh and Rafael Benguria,
\emph{Optimal lower bound for the gap between the first two
eigenvalues of one-dimensional Schr\"odinger operators
with symmetric single-well potentials},
Proc. Amer. Math. Soc. \textbf{105} (1989), no.~2, 419--424.

\bibitem{BanuelosKulczycki2004}
Rodrigo Ba\~nuelos and Tadeusz Kulczycki,
\emph{The Cauchy process and the Steklov problem},
J. Funct. Anal. \textbf{211} (2004), no.~2, 355--423.

\bibitem{BanuelosKulczycki2006JFA}
Rodrigo Ba\~nuelos and Tadeusz Kulczycki,
\emph{Eigenvalue gaps for the Cauchy process and a
Poincar\'e inequality},
J. Funct. Anal. \textbf{234} (2006), no.~1, 199--225.

\bibitem{BanuelosKulczycki2006CPDE}
Rodrigo Ba\~nuelos and Tadeusz Kulczycki,
\emph{Spectral gap for the Cauchy process on convex,
symmetric domains},
Comm. Partial Differential Equations
\textbf{31} (2006), no.~10--12, 1841--1878.

\bibitem{BanuelosKulczyckiMendez2006}
Rodrigo Ba\~nuelos, Tadeusz Kulczycki,
and Pedro J. M\'endez-Hern\'andez,
\emph{On the shape of the ground state eigenfunction
for stable processes},
Potential Anal. \textbf{24} (2006), no.~3, 205--221.

\bibitem{BaoChenJiangMa2020}
Weizhu Bao, Lizhen Chen, Xiaoyun Jiang, and Ying Ma,
\emph{A Jacobi spectral method for computing eigenvalue gaps
and their distribution statistics of the fractional
Schr\"odinger operator},
J. Comput. Phys. \textbf{421} (2020), 109733.

\bibitem{BaoRuanShenSheng2019}
Weizhu Bao, Xinran Ruan, Jie Shen, and Changtao Sheng,
\emph{Fundamental gaps of the fractional Schr\"odinger operator},
Commun. Math. Sci. \textbf{17} (2019), no.~2, 447--471.

\bibitem{BorisovFreitas2009}
Denis Borisov and Pedro Freitas,
\emph{Singular asymptotic expansions for Dirichlet eigenvalues
and eigenfunctions of the Laplacian on thin planar domains},
Ann. Inst. H. Poincar\'e C Anal. Non Lin\'eaire
\textbf{26} (2009), no.~2, 547--560.

\bibitem{ChandlerWildeHewettMoiola2017}
Simon N. Chandler-Wilde, David P. Hewett, and Andrea Moiola,
\emph{Sobolev spaces on non-Lipschitz subsets of
$\mathbb{R}^n$ with application to boundary integral
equations on fractal screens},
Integral Equations Operator Theory
\textbf{87} (2017), no.~2, 179--224.

\bibitem{ChenHauer2026Gap}
Rui Chen and Daniel Hauer,
\emph{Fundamental gaps for the Dirichlet $p$-Laplacian
on convex domains},
2026, arXiv:2608.13443v3.

\bibitem{ChenCheng2020}
Wei-Chuan Chen and Yan-Hsiou Cheng,
\emph{Remarks on the one-dimensional sloshing problem
involving the $p$-Laplacian operator},
Turkish J. Math. \textbf{44} (2020), no.~4, 1376--1387.

\bibitem{ChengLianWang2014}
Yan-Hsiou Cheng, Wei-Cheng Lian, and Wei-Chuan Wang,
\emph{The dual eigenvalue problems for $p$-Laplacian},
Acta Math. Hungar. \textbf{142} (2014), no.~1, 132--151.

\bibitem{DelPezzoQuaas2017}
L. M. Del Pezzo and A. Quaas,
\emph{A Hopf's lemma and a strong minimum principle
for the fractional $p$-Laplacian},
J. Differential Equations
\textbf{263} (2017), no.~1, 765--778.

\bibitem{DiNezzaPalatucciValdinoci2012}
Eleonora Di Nezza, Giampiero Palatucci, and Enrico Valdinoci,
\emph{Hitchhiker's guide to the fractional Sobolev spaces},
Bull. Sci. Math. \textbf{136} (2012), no.~5, 521--573.

\bibitem{DydaKulczycki2007}
Bart\l omiej Dyda and Tadeusz Kulczycki,
\emph{Spectral gap for stable process on convex planar
double symmetric domains},
Potential Anal. \textbf{27} (2007), no.~2, 101--132.

\bibitem{FallGhimentiMichelettiPistoia2023}
Mouhamed Moustapha Fall, Marco Ghimenti,
Anna Maria Micheletti, and Angela Pistoia,
\emph{Generic properties of eigenvalues of the
fractional Laplacian},
Calc. Var. Partial Differential Equations
\textbf{62} (2023), Paper No.~233.

\bibitem{FleckingerHarrellDeThelin2007}
Jacqueline Fleckinger, Evans M. Harrell II,
and Fran\c{c}ois de Th\'elin,
\emph{On the fundamental eigenvalue ratio of the
$p$-Laplacian},
Bull. Sci. Math. \textbf{131} (2007), no.~7, 613--619.

\bibitem{FriedlanderSolomyak2009}
Leonid Friedlander and Michael Solomyak,
\emph{On the spectrum of the Dirichlet Laplacian
in a narrow strip},
Israel J. Math. \textbf{170} (2009), 337--354.

\bibitem{GongLiLuo2016}
Fuzhou Gong, Huaiqian Li, and Dejun Luo,
\emph{A probabilistic proof of the fundamental gap
conjecture via the coupling by reflection},
Potential Anal. \textbf{44} (2016), no.~3, 423--442.

\bibitem{GongLiLuo2016Erratum}
Fuzhou Gong, Huaiqian Li, and Dejun Luo,
\emph{Erratum to: A probabilistic proof of the
fundamental gap conjecture via the coupling by reflection},
Potential Anal. \textbf{44} (2016), no.~3, 443--447.

\bibitem{Grubb2015}
Gerd Grubb,
\emph{Fractional Laplacians on domains, a development
of H\"ormander's theory of $\mu$-transmission
pseudodifferential operators},
Adv. Math. \textbf{268} (2015), 478--528.

\bibitem{Kaleta2012}
Kamil Kaleta,
\emph{Spectral gap lower bound for the one-dimensional
fractional Schr\"odinger operator in the interval},
Studia Math. \textbf{209} (2012), no.~3, 267--287.

\bibitem{Kato1995}
Tosio Kato,
\emph{Perturbation theory for linear operators},
Classics in Mathematics, Springer, Berlin, 1995,
reprint of the 1980 edition.

\bibitem{KulczyckiKwasnickiMaleckiStos2010}
Tadeusz Kulczycki, Mateusz Kwa\'snicki,
Jacek Ma\l ecki, and Andrzej St\'os,
\emph{Spectral properties of the Cauchy process
on half-line and interval},
Proc. Lond. Math. Soc. (3)
\textbf{101} (2010), no.~2, 589--622.

\bibitem{KuwaeShioya2003}
Kazuhiro Kuwae and Takashi Shioya,
\emph{Convergence of spectral structures:
a functional analytic theory and its applications
to spectral geometry},
Comm. Anal. Geom. \textbf{11} (2003), no.~4, 599--673.

\bibitem{Kwasnicki2008}
Mateusz Kwa\'snicki,
\emph{Spectral gap estimate for stable processes
on arbitrary bounded open sets},
Probab. Math. Statist. \textbf{28} (2008), no.~1, 163--167.

\bibitem{Kwasnicki2012}
Mateusz Kwa\'snicki,
\emph{Eigenvalues of the fractional Laplace operator
in the interval},
J. Funct. Anal. \textbf{262} (2012), no.~5, 2379--2402.

\bibitem{Kwasnicki2017}
Mateusz Kwa\'snicki,
\emph{Ten equivalent definitions of the
fractional Laplace operator},
Fract. Calc. Appl. Anal.
\textbf{20} (2017), no.~1, 7--51.

\bibitem{Kwasnicki2023Simplicity}
Mateusz Kwa\'snicki,
\emph{Simplicity of eigenvalues of the fractional
Laplace operator in an interval},
2023, revised 2026, arXiv:2311.00713v2.

\bibitem{KwasnickiWszola2026}
Mateusz Kwa\'snicki and Jacek Wszo\l a,
\emph{Dirichlet--Neumann bracketing for nonlocal operators},
2026, arXiv:2609.24892v1.

\bibitem{Lavine1994}
Richard Lavine,
\emph{The eigenvalue gap for one-dimensional convex potentials},
Proc. Amer. Math. Soc. \textbf{121} (1994), no.~3, 815--821.


\bibitem{LuRowlett2013Simplices}
Zhiqin Lu and Julie Rowlett,
The fundamental gap of simplices,
Comm. Math. Phys. 319 (2013), no.~1, 111--145.


\bibitem{LuRowlett2014}
Zhiqin Lu and Julie Rowlett,
\emph{The fundamental gap and one-dimensional collapse},
in \emph{Geometric and spectral analysis},
Contemp. Math., vol.~630,
Amer. Math. Soc., Providence, RI, 2014, pp.~223--246.



\bibitem{Ni2013}
Lei Ni,
\emph{Estimates on the modulus of expansion for
vector fields solving nonlinear equations},
J. Math. Pures Appl. (9)
\textbf{99} (2013), no.~1, 1--16.

\bibitem{ReedSimonIV}
Michael Reed and Barry Simon,
\emph{Methods of modern mathematical physics.
IV. Analysis of operators},
Academic Press, New York, 1978.

\bibitem{RockafellarWets1998}
R. Tyrrell Rockafellar and Roger J.-B. Wets,
\emph{Variational analysis},
Grundlehren der mathematischen Wissenschaften,
vol.~317, Springer, Berlin, 1998.

\bibitem{RosOtonSerra2014}
Xavier Ros-Oton and Joaquim Serra,
\emph{The Dirichlet problem for the fractional Laplacian:
regularity up to the boundary},
J. Math. Pures Appl. (9)
\textbf{101} (2014), no.~3, 275--302.

\bibitem{SingerEtAl1985}
I. M. Singer, Bun Wong, Shing-Tung Yau, and Stephen S.-T. Yau,
\emph{An estimate of the gap of the first two
eigenvalues in the Schr\"odinger operator},
Ann. Scuola Norm. Sup. Pisa Cl. Sci. (4)
\textbf{12} (1985), no.~2, 319--333.

\bibitem{Smits1996}
Robert G. Smits,
\emph{Spectral gaps and rates to equilibrium
for diffusions in convex domains},
Michigan Math. J. \textbf{43} (1996), no.~1, 141--157.

\bibitem{Soltan2005}
Valeriu Soltan,
\emph{Affine diameters of convex bodies---a survey},
Expo. Math. \textbf{23} (2005), no.~1, 47--63.

\bibitem{vandenBerg1983}
M. van den Berg,
\emph{On condensation in the free-boson gas
and the spectrum of the Laplacian},
J. Statist. Phys. \textbf{31} (1983), no.~3, 623--637.

\bibitem{WenZhou2026}
Zhiyuan Wen and Lijuan Zhou,
\emph{Maximal gaps of eigenvalues of one-dimensional
$p$-Laplacian},
Mediterr. J. Math. \textbf{23} (2026), no.~1,
Paper No.~47.

\bibitem{Yau1986}
Shing-Tung Yau,
\emph{Nonlinear analysis in geometry},
Monographies de L'Enseignement Math\'ematique,
vol.~33, L'Enseignement Math\'ematique, Geneva, 1986.

\bibitem{YessirkegenovZhangirbayev2026}
Nurgissa Yessirkegenov and Amir Zhangirbayev,
\emph{Explicit sharp gradient stability of the
$L^p$-Poincar\'e inequality and fundamental gaps
for the Dirichlet $p$-Laplacian operator},
2026, arXiv:2610.27281v1.

\bibitem{YuZhong1986}
Qi Huang Yu and Jia Qing Zhong,
\emph{Lower bounds of the gap between the first
and second eigenvalues of the Schr\"odinger operator},
Trans. Amer. Math. Soc.
\textbf{294} (1986), no.~1, 341--349.

\bibitem{Zhang2026Interval}
Cheng Zhang,
\emph{Eigenvalues and eigenfunctions of the
fractional Laplacian on the interval},
2026, arXiv:2608.23457v2.

\end{thebibliography}
	\providecommand{\bysame}{\leavevmode\hbox to3em{\hrulefill}\thinspace}
	\providecommand{\MR}{\relax\ifhmode\unskip\space\fi MR }
	% \MRhref is called by the amsart/book/proc definition of \MR.
	\providecommand{\MRhref}[2]{%
		\href{http://www.ams.org/mathscinet-getitem?mr=#1}{#2}
	}
	\providecommand{\href}[2]{#2}

\end{document}